\documentclass[12pt,a4paper]{amsart}
\usepackage{comment}
\usepackage[foot]{amsaddr}

\usepackage{mathrsfs}

\usepackage{mathtools}
\makeatletter
\tagsleft@false
\makeatother

\usepackage{mathrsfs}
\usepackage{tikz}
\usepackage{amsfonts}
\usepackage{amsmath}
\usepackage{amssymb}
\usepackage{enumerate}
\usepackage{hyperref}
\usepackage{latexsym}
\usepackage{xcolor}
\usepackage[shortlabels]{enumitem}
\usepackage[numbers,sort&compress]
{natbib}
\setlist[enumerate,1]{label={\upshape(\roman*)}}

 \makeatletter
    
    \newcommand{\Rmnum}[1]
    {\expandafter\@slowromancap\romannumeral #1@}
    \makeatother

\newtheorem{thm}{Theorem}[section]

\newtheorem{prop}[thm]{Proposition}
\newtheorem{lemma}[thm]{Lemma}

\newtheorem{cor}[thm]{Corollary}

\newtheorem{example}[thm]{Example}
\newtheorem{defin}[thm]{Definition}

\theoremstyle{definition}
\newtheorem{remark}[thm]{Remark}

\def\wz{\tilde}

\title[$s$-quasi-transitive weakly distance-regular digraphs]{$s$-quasi-transitive weakly distance-regular digraphs}

\date{}
\thanks{*Corresponding author}

\author[Li]{Shuang Li$^{\lowercase{A,B}}$}
\address{$^{a}$School of Mathematics and Systems Science\\Wuhan University of Science and Technology\\Wuhan 430065\\China}
\address{$^{b}$Laboratory of Mathematics and Complex Systems (MOE),~School of Mathematical Sciences\\Beijing Normal University\\Beijing 100875\\China} 
\email{lishuangyx@mail.bnu.edu.cn}

\author[Yang]{Yuefeng Yang$^{\lowercase{C},*}$}
\address{$^{c}$School of Science\\China University of Geosciences\\Beijing 100083\\China}
\email{yangyf@cugb.edu.cn}

\begin{document}

\begin{abstract}
	A digraph is $s$-quasi-transitive if,  for every path $(x_0,x_1,\ldots,x_s)$ of length $s$, 
	either $(x_0,x_s)\in A\Gamma$ or $(x_s,x_0)\in A\Gamma$ (or both). 
	In this paper,  we give a characterization of $s$-quasi-transitive commutative weakly distance-regular digraphs, where  $s=4$, or $s\geq 5$ and the cycle polar diameter is at least $s$.


\end{abstract}
\maketitle
\section{Introduction}
A \emph{digraph} $\Gamma$ is a pair $(V\Gamma,A\Gamma)$, where $V\Gamma$ is a finite nonempty set of vertices and $A\Gamma$ is a set of ordered pairs ({\em arcs}) $(x,y)$ with distinct vertices $x$ and $y$. A subdigraph of $\Gamma$ induced by a subset $U\subseteq V\Gamma$ is denoted by $\Gamma[U]$. For any arc $(x,y)\in A\Gamma$, if $A\Gamma$ also contains an arc $(y,x)$, then $\{(x,y),(y,x)\}$ can be viewed as an {\em edge}. We say that $\Gamma$ is an \emph{undirected graph} or a {\em graph} if $A\Gamma$ is a symmetric relation.
A vertex $x$ is {\em adjacent} to $y$ if $(x,y)\in A\Gamma$.
In this case, we also call $y$ an \emph{out-neighbour} of $x$, and $x$ an \emph{in-neighbour} of $y$. 
A digraph is said to be {\em regular of valency} $k$ if the number of in-neighbour and out-neighbour of all vertices are equal to $k$.

For the digraphs $\Gamma$ and $\Sigma$, the {\em lexicographic product} $\Gamma\circ\Sigma$ is the digraph with vertex set $V\Gamma\times V\Sigma$ such that $((u_1,u_2),(v_1,v_2))$ is an arc if and only if either $(u_1,v_1)\in A\Gamma$, or
$u_1=v_1$ and $(u_2,v_2)\in A\Sigma$. 
If $\Sigma$ is the edgeless digraph on $m$ vertices, then $\Gamma\circ \Sigma$ is called an {\em$m$-coclique extension} of $\Gamma$.

A \emph{walk} of length $r$ from $x$ to $y$ in the digraph $\Gamma$ is a finite sequence of vertices $W:=(x=w_{0},w_{1},\ldots,w_{r}=y)$ such that $(w_{t-1}, w_{t})\in A\Gamma$ for $1\leq t\leq r$. If all vertices of  $W$ are distinct, then $W$ is a \emph{path} of length $r$. If $(y,x)\in A\Gamma$, $W$ is a \emph{circuit} of length $r+1$.
A digraph (resp. graph) is said to be \emph{strongly connected} (resp. \emph{connected}) if, for any vertices $x$ and $y$, there is a path from $x$ to $y$.  The \emph{girth} of $\Gamma$ is the length of a shortest circuit in $\Gamma$.

Quasi-transitive digraphs were originally introduced by Ghouila-Houri in 1962, this class of digraphs became very popular after Bang-Jensen and Huang proved the so-called Canonical Decomposition Theorem in 1995 \cite{JBJ95}. A digraph is \emph{quasi-transitive} if,  for every path $(x_0,x_1,x_2)$ of length $2$, either $(x_0,x_2)\in A\Gamma$ or $(x_2,x_0)\in A\Gamma$ (or both).
Quasi-transitive digraphs are a natural generalization of both semicomplete digraphs and transitive digraphs. Many nice results for these classes can be generalized for quasi-transitive digraphs: the book chapter \cite{Gal17} by Galeana-S\'{a}nchez and Hern\'{a}ndez-Cruz surveys many of these results.

The concept of $s$-quasi-transitive digraphs was introduced in \cite{Her12} as a generalization of quasi-transitive digraphs. A digraph is \emph{$s$-quasi-transitive} if, for every path $(x_0,x_1,\ldots,x_s)$ of length $s$, either $(x_0,x_s)\in A\Gamma$ or $(x_s,x_0)\in A\Gamma$ (or both).  The $s$-quasi-transitive digraphs have been studied in \cite{Gal13,Her12,RW15,RWM15}.
In \cite{JBJ95}, the authors characterized quasi-transitive digraphs in a recursive fashion. In \cite{Gal10}, Galeana-S\'{a}nchez et al.  proved that a strongly connected $3$-quasi-transitive digraph is either a semicomplete bipartite digraph, a semicomplete digraph, or isomorphic to a digraph in an infinite (but very restricted) family.

Let $\Gamma$ be a strongly connected digraph.  The length of a shortest path from $x$ to $y$ is called the \emph{distance} from $x$ to $y$ in $\Gamma$, denoted by $\partial_\Gamma(x,y)$. The maximum value of the distance function in $\Gamma$ is called the \emph{diameter} of $\Gamma$. For $0\leq i\leq d$,  we define $\Gamma_{i}$ to be the set of ordered pairs $(x,y)$, 
where $d$ is the diameter of $\Gamma$. 
 Let $\wz{\partial}_{\Gamma}(x,y):=(\partial_{\Gamma}(x,y),\partial_{\Gamma}(y,x))$ be the \emph{two-way distance} from $x$ to $y$, and $\wz{\partial}(\Gamma):=\{\wz{\partial}_{\Gamma}(x,y)\mid x,y\in V\Gamma\}$ the \emph{two-way distance set} of $\Gamma$. If no confusion occurs, we write $\partial(x,y)$ (resp. $\tilde{\partial}(x,y)$) instead of $\partial_\Gamma(x,y)$ (resp. ${\tilde{\partial}}_\Gamma(x,y)$).  
 For any $\tilde{i}\in\tilde{\partial}(\Gamma)$, we define $\Gamma_{\tilde{i}}$ to be the set of ordered pairs $(x,y)$ with $\tilde{\partial}_{\Gamma}(x,y)=\tilde{i}$.
  An arc $(x,y)$ of $\Gamma$ is of {\em type} $(1,r)$ if $\partial(y,x)=r$. Let $T$ be the set of integers $q$ such that $(1,q-1)$ is a type of $\Gamma$. That is, $T=\{q\mid (1,q-1)\in \tilde{\partial}(\Gamma)\}$, and the maximum element of $T$ is called  the \emph{cycle polar diameter} of $\Gamma$.


As a natural generalization of distance-regular graphs (see \cite{NB74,AEB98,DKT16} for a background of the theory of distance-regular graphs),  Wang and Suzuki \cite{KSW03} introduced the concept of weakly distance-regular digraphs. A strongly connected digraph $\Gamma$ is said to be \emph{weakly distance-regular} if the configuration $\mathfrak{X}(\Gamma)=(V\Gamma,\{\Gamma_{\wz{i}}\}_{\wz{i}\in \wz{\partial}(\Gamma)})$ is an association scheme (see the definition of association scheme in Section 2). We call $\mathfrak{X}(\Gamma)$ the \emph{attached scheme} of $\Gamma$. A weakly distance-regular digraph is \emph{commutative} if its attached scheme is commutative. 
For more information about weakly distance-regular digraphs, see \cite{YF22,lishuang2,HS04,AM,KSW03,KSW04,YYF16,YYF18,YYF20,YYF22,YYF,QZ23,QZ24}. Motivated by the structural results for quasi-transitive  digraphs and $3$-quasi-transitive  digraphs \cite{JBJ95,Gal10}, we focus on the case $s\geq 4$.




 In this paper,  we give a characterization of $s$-quasi-transitive commutative weakly distance-regular digraphs, where  $s=4$, or $s\geq 5$ and the cycle polar diameter is at least $s$. Our main results are as follows.

\begin{thm}\label{main thm}
	 Let $\Gamma$ be a  commutative weakly distance-regular digraph of valency greater than $3$ and with  cycle polar diameter at least $s$, where  $s\geq 4$.  Then $\Gamma$ is  $s$-quasi-transitive  if and only if it is isomorphic to one of the following digraphs:
	\begin{enumerate}

		 	\item\label{main thm-4}
		 $\Sigma\circ C_4$;

		
		\item\label{main thm-3}
	a $t$-coclique extension of  $C_{s+1}$.
	\end{enumerate}
Here,  $t>3$ and $\Sigma$ is a  semicomplete weakly distance-regular digraph with girth $3$.
\end{thm}

	
\begin{thm}\label{main thm1}
 Let $\Gamma$ be a  commutative weakly distance-regular digraph  of valency greater than $3$.  Then $\Gamma$ is  $4$-quasi-transitive  if and only if it is isomorphic to one of the following digraphs:
	\begin{enumerate}
		\item\label{main thm1-2}
		$	\Lambda$;
		
		
		\item\label{main thm1-3}
		$\Sigma\circ C_4$;

		
		
		\item\label{main thm1-7}
		a $t$-coclique extension of  $\Sigma$;
		
	\item\label{main thm1-6}
	an $l$-coclique extension of  $C_{5}$.
	\end{enumerate}
	Here,  $t\geq 1$, $l>3$ and $\Lambda$ (resp. $\Sigma$) is a  semicomplete weakly distance-regular digraph with girth $2$ (resp. $3$).
\end{thm}

The remainder of this paper is organized as follows. In Section 2, we provide the required notation and concepts for association schemes.
In Section 3, we establish preliminary results. 
In Sections 4 and 5, we treat the cases in which the  cycle polar diameter is $s$ or $s+1$.
In Section 6,  we give the proofs of  Theorems \ref{main thm} and \ref{main thm1}.


\section{Association schemes}

In this section, we present some notation and concepts for association schemes that we shall use in the remainder of the paper.

A \emph{$d$-class association scheme} $\mathfrak{X}$ is a pair $(X,\{R_{i}\}_{i=0}^{d})$, where $X$ is a finite nonempty set, and each $R_{i}$ is a
nonempty subset of $X\times X$ satisfying the following axioms (see \cite{EB21,EB84,PHZ96,PHZ05} for a background of the theory of association schemes):
\begin{enumerate}
\item\label{as-1} $R_{0}=\{(x,x)\mid x\in X\}$ is the diagonal relation;

\item\label{as-2} $X\times X=R_{0}\cup R_{1}\cup\cdots\cup R_{d}$, $R_{i}\cap R_{j}=\emptyset~(i\neq j)$;

\item\label{as-3} for each $i$, $R_{i}^{\top}=R_{i^{*}}$ for some $0\leq i^{*}\leq d$, where $R_{i}^{\top}=\{(y,x)\mid(x,y)\in R_{i}\}$;

\item\label{as-4} for all $i,j,l$, the cardinality of the set $$P_{i,j}(x,y):=R_{i}(x)\cap R_{j^{*}}(y)$$ is constant whenever $(x,y)\in R_{l}$, where $R(x)=\{y\mid (x,y)\in R\}$ for $R\subseteq X\times X$ and $x\in X$. This constant is denoted by $p_{i,j}^{l}$.
\end{enumerate}
A $d$-class association scheme is also called an \emph{association scheme} with $d$ classes (or even simply a scheme). The integers $p_{i,j}^{l}$ are called the \emph{intersection numbers} of $\mathfrak{X}$. We say that $\mathfrak{X}$ is \emph{commutative} if $p_{i,j}^{l}=p_{j,i}^{l}$ for all $i,j,l$. The subsets $R_{i}$ are called the \emph{relations} of $\mathfrak{X}$. For each $i$, the integer $k_{i}:=p_{i,i^{*}}^{0}$ is called the \emph{valency} of $R_{i}$. 

Let $\mathfrak{X}=(X,\{R_{i}\}_{i=0}^{d})$ be an association scheme. For two nonempty subsets $E$ and $F$ of $\{R_{i}\}_{i=0}^{d}$, define
\[EF:=\{R_{l}\mid\sum_{R_{i}\in E}\sum_{R_{j}\in F}p_{i,j}^{l}\neq0\}.\]
We write $R_{i}R_{j}$ instead of $\{R_{i}\}\{R_{j}\}$,
and $R_i^2$ instead of $\{R_{i}\}\{R_{i}\}$.
If $R_{i^{*}}R_{j}\subseteq F$ for any $R_{i},R_{j}\in F$, we say that $F$ is {\em closed}.
For each $R_i$ with $0\leq i\leq d$, define $\left\langle R_i\right\rangle$ to be the smallest closed subset containing $R_i$. We say that $R_{j}$ {\em generates} $\mathfrak{X}$ if $\left\langle R_j\right\rangle=\{R_{i}\}_{i=0}^{d}$.


In this section, we always assume that $(X,\{R_{i}\}_{i=0}^{d})$ is an association scheme. We list some basic properties of intersection numbers which are used frequently in this paper.

\begin{lemma}\label{jb}
	{\rm (\cite[Chapter \Rmnum{2}, Proposition 2.2]{EB84} and \cite[Proposition 5.1]{ZA99})} Let $i,j,h,l\in\{0,1,\ldots,d\}$. The following hold:
	\begin{enumerate}
\item\label{jb-0} $p_{i,j}^{l}=p_{i^{*},j^{*}}^{l^{*}}$;

		\item\label{jb-1} $k_{i}k_{j}=\sum_{r=0}^{d}p_{i,j}^{r}k_{r}$;
		
		\item\label{jb-2}  $p_{i,j}^{l}k_{l}=p_{l,j^{*}}^{i}k_{i}=p_{i^{*},l}^{j}k_{j}$;
		
		\item\label{jb-3}  $\sum_{r=0}^{d}p_{i,r}^{l}=k_{i}$;
		
		\item\label{jb-4} $\sum_{r=0}^{d}p_{i,l}^{r}p_{h,r}^{j}=\sum_{t=0}^{d}p_{h,i}^{t}p_{t,l}^{j}$;

\item\label{jb-5}
		$|R_iR_j|\leq \emph{gcd}(k_{i},k_{j})$.
	\end{enumerate}
\end{lemma}

\begin{lemma}\label{comm}
Let $i,j,h\in\{0,1,\ldots,d\}$. Then $R_{h}R_{j^*}=\{R_i\}$ if and only if $p_{i,j}^{h}=k_{j}$.
\end{lemma}
\begin{proof}
	Since 	$k_{j}=k_{j^*}$ from Lemma \ref{jb} \ref{jb-0}, by  Lemma \ref{jb} \ref{jb-1}, one has $k_{h}k_{j}=k_{h}k_{j^*}=\sum_{r=0}^{d}p_{h,j^*}^{r}k_{r}.$ This implies that $R_{h}R_{j^*}=\{R_i\}$ if and only if $k_{h}k_{j}=p_{h,j^*}^{i}k_{i}$. In view of Lemma \ref{jb} \ref{jb-2}, one gets $p_{h,j^*}^{i}k_{i}=p_{i,j}^{h}k_{h}$. It follows that $R_{h}R_{j^*}=\{R_i\}$ if and only if $k_{h}k_{j}=p_{i,j}^{h}k_{h}$. Thus, the desired result follows.
\end{proof}


 \begin{lemma}\label{one}
Let $i\in\{1,2,\ldots,d\}$. Then $R_i^2\ne \{R_i\}$.
\end{lemma}
\begin{proof}
	Suppose for the contrary that  $R_i^2= \{R_i\}$. By Lemma \ref{comm}, we get $p_{i,i^*}^{i}=k_{i^*}=k_i$. In view of Lemma \ref{jb} \ref{jb-3}, one has $k_i\geq p_{i,i^*}^{i}+p_{i,0}^{i}=k_i+1$, a contradiction.
\end{proof}	

\begin{lemma}\label{comm0}
Let $k_{j}=1$ and $R_i\in R_hR_j$ with $i,j,h\in\{0,1,\ldots,d\}$. Then $R_hR_j=\{R_i\}$ and $R_iR_{j^*}=\{R_h\}$.
\end{lemma}
\begin{proof}
	Since $p_{h,j}^{i}\ne 0$, from Lemma \ref{jb} \ref{jb-2}, one has $p_{i,j^*}^{h}\ne 0$, and so $R_h\in R_iR_{j^*}$. Since $k_{j}=1$, from Lemma \ref{jb} \ref{jb-5}, the desired result follows.
\end{proof}

\section{Preliminary results on $s$-quasi-transitive digraphs}
In this section, concerning $s$-quasi-transitive  digraph, we present some basic results  that we shall use in the remainder of the paper.

\begin{lemma}\label{even}
	{\rm (\cite{Her12})} Let  $\Gamma$ be a strongly connected $s$-quasi-transitive  digraph with $s\geq 2$. Suppose $u,v\in V\Gamma$. Then  the following hold:
	\begin{enumerate}
		\item\label{even-1} If $\partial(u,v)=s$, then $\partial(v,u)=1$;
		
		\item\label{even-2}  If $\partial(u,v)=s+1$, then $\partial(v,u)\leq s+1$;
		
		\item\label{even-3} Assume $\partial(u,v)=n\geq s+2$. If $s$ is even, or $s$ and $n$ are both odd,  then $\partial(v,u)=1$; if $s$ is odd and $n$ is even, then $\partial(v,u)\leq 2$.
	\end{enumerate}
\end{lemma}

In the remainder of this section, $\Gamma$ always denotes   a $s$-quasi-transitive   weakly distance-regular digraph  with  $s\geq 4$.

\begin{lemma}\label{maxt}
	If $q\in T$, then $q\leq s+1$.
\end{lemma}
\begin{proof}
	Suppose for the contrary that there exists $q\in T$ such that $q\geq s+2$. Let $(x_0,x_1,\ldots,x_{m-1})$ be a circuit  consisting of arcs of type $(1,q-1)$. Then $m\geq q\geq s+2$.

Let $t$ and $r$ be integers such that $m-1-s=t(s-1)+r$ with $0\leq r<s-1$. We claim that $(x_0,x_{s+k(s-1)})\in A\Gamma$ for $-1\leq k\leq t$. We prove it by induction on $k$. If $k=-1$, then the desired result follows. Assume that $(x_0,x_{s+i(s-1)})\in A\Gamma$ with $-1\leq i<t$. Since $\Gamma$ is $s$-quasi-transitive and $(x_0,x_{s+i(s-1)},x_{s+i(s-1)+1},\ldots,x_{s+(i+1)(s-1)})$ is a path of length $s$, we obtain either $(x_0,x_{s+(i+1)(s-1)})\in A\Gamma$ or $(x_{s+(i+1)(s-1)},x_0)\in A\Gamma$.
	If $(x_{s+(i+1)(s-1)},x_0)\in A\Gamma$, then $(x_{s+i(s-1)},x_{s+i(s-1)+1},\ldots,x_{s+(i+1)(s-1)},x_0)$ is a circuit of length $s+1$ containing an arc of type $(1,q-1)$, contrary to the fact that $q\geq s+2$. Thus, our claim is valid.

	Since $m-1-s=t(s-1)+r$ and $(x_0,x_{s+t(s-1)})\in A\Gamma$ from the claim, one has $$\partial(x_0,x_{m-1})\leq \partial(x_0,x_{s+t(s-1)})+\partial(x_{s+t(s-1)},x_{m-1})=1+\partial(x_{s+t(s-1)},x_{s+t(s-1)+r}),$$ and so $q-1\leq 1+r<s$, a contradiction.
\end{proof}

For each $\wz{i}:=(a,b)\in \wz{\partial}(\Gamma)$, we write $k_{a,b}$ instead of $k_{(a,b)}$, and  $\Gamma_{a,b}$ instead of $\Gamma_{(a,b)}$.

\begin{lemma}\label{rem}
	$\Gamma$ has diameter at most $s+1$.
	\end{lemma}
\begin{proof}
Assume the contrary, namely, there exist vertices $u,v$ such that $\partial(u,v)=s+2$. If $s$ is even, from Lemma \ref{even} \ref{even-3}, then $(v,u)\in \Gamma_{1,s+2}$; if  $s$ is odd, then  $\partial(u,v)=s+2$ is odd, which implies that  $(v,u)\in \Gamma_{1,s+2}$ from Lemma \ref{even} \ref{even-3} again, both contrary to Lemma \ref{maxt}.
\end{proof}

\begin{lemma}\label{1009}
	Let $(v,v_1,\ldots,v_{l-1},u)$ be a shortest path. Then $|\Gamma_{p}(v)\cap \{v_1,\ldots,v_{l-1}\}|\leq 1$ for any $1\leq p\leq l-1$.
\end{lemma}
\begin{proof}
		If $w \in \Gamma_{p}(v)$ and $w\in \{v_1,\ldots,v_{l-1}\}$, then $w=v_{p}$ since $\partial(v,w)=p$. Thus, the desired result follows.
\end{proof}


\begin{lemma}\label{jichu0}
Let $t+1\in T$. Suppose $p_{(r,s-1),(t,1)}^{(i,j)}>0$ with $r\geq 1$ and $(i,j)\ne (0,0)$. Then at least one of the following holds: $i=1$, $j=1$, or $j=s-t-1$.
\end{lemma}
\begin{proof}
	Let $(x,z)\in \Gamma_{i,j}$. Then there exists $y\in P_{(r,s-1),(t,1)}(x,z)$. Pick a shortest path $(y,y_1,\ldots,y_{s-2},x)$. Suppose  $z\notin \{y_1,\ldots,y_{s-2}\}$. Since $\Gamma$ is $s$-quasi-transitive and $(z,y,y_1,\ldots,y_{s-2},x)$ is a path of length $s$, one has $i=1$ or $j=1$.
	Suppose $z\in \{y_1,\ldots,y_{s-2}\}$. Since $\partial(y,z)=t$ and $(y,y_1,\ldots,y_{s-2},x)$ is a shortest path, we get
	$z=y_t$, which implies $j=\partial(z,x)=\partial(y_t,x)=s-t-1$.
\end{proof}

\begin{lemma}\label{k_{1,s-1}=1}
If $\Gamma_{1,s-1}\notin \Gamma_{1,s-1}\Gamma_{1,q}$ for any $q+1\in T$, then $\Gamma_{1,s-1}\Gamma_{s-1,1}=\{\Gamma_{0,0}\}$ and  $k_{1,s-1}=1$.
\end{lemma}
\begin{proof}
In view of Lemma \ref{jb} \ref{jb-0} and \ref{jb-2}, we obtain $\Gamma_{1,q},\Gamma_{q,1}\notin \Gamma_{1,s-1}\Gamma_{s-1,1}$ for any $q+1\in T$. By setting $r=1$ and $t=s-1$ in   Lemma \ref{jichu0}, one has $\Gamma_{1,s-1}\Gamma_{s-1,1}=\{\Gamma_{0,0}\}$, and so $k_{1,s-1}=1$.
\end{proof}



\begin{lemma}\label{1207}
	Suppose  $(x,y)\in \Gamma_{2,s+(-1)^l}$ and $z\in  P_{(1,m),(1,n)}(x,y)$ with $l\in \{0,1\}$, $m+n>s+(-1)^l$ and $m\geq 4+(-1)^l$. Let $(y=y_0,y_1,\ldots,y_{s+(-1)^l-1},x)$  be a shortest path. If  $z_0\in  P_{(1,p-1),(1,q-1)}(z,y)$ for some $p,q\in T$,
	then $z_0\in \{y_1,\ldots,y_{s+(-1)^l-1}\}$.
\end{lemma}
\begin{proof}
	Assume the contrary, namely, $z_0\notin \{y_1,\ldots,y_{s+(-1)^l-1}\}$. Since $(x,y)\in \Gamma_{2,s+(-1)^l}$ and $(z_0,y)\in \Gamma_{1,q-1}$, one has $z_0\notin \{x,y\}$. If $z\in \{y_1,\ldots,y_{s+(-1)^l-1}\}$,  then $z=y_n$ since $\partial(y,z)=n$, which implies $m=\partial(z,x)=\partial(y_n,x)=s+(-1)^l-n$, contrary to the fact that $m+n>s+(-1)^l$. Thus,  $z\notin \{y_1,\ldots,y_{s+(-1)^l-1}\}$.
	Since $(y,y_1,\ldots,y_{s+(-1)^l-1},x)$  is a shortest path, one has 
	\begin{align}
4+(-1)^l\leq m=\partial(z,x)\leq \partial(z,y_{s-2})+(-1)^l+2.\label{261016}\tag{3.1}
	\end{align}
	Since $(x,z,z_0,y,y_1,\ldots,y_{s-3})$ and $(z,z_0,y,y_1,\ldots,y_{s-2})$ are two paths of length $s$,  one has $(x,y_{s-3})\in A\Gamma$, and $(z,y_{s-2})$ or $(y_{s-2},z)\in A\Gamma$. Since $l=0$ or $l=1$, from \eqref{261016}, we get $5\leq \partial(z,y_{s-2})+3$ or  $3\leq \partial(z,y_{s-2})+1$, which implies $(y_{s-2},z)\in A\Gamma$.

	We prove $(x,y_{s-i})\in A\Gamma$  for any  $3\leq i\leq s$ by induction on $i$. If $i=3$,
	then  $(x,y_{s-3})\in A\Gamma$. Now suppose $(x,y_{s-k})\in A\Gamma$ with $3\leq k<s$.
	Note that $(y,y_1,\ldots,y_{s+(-1)^l-1},x)$  is a shortest path.
	Since
	$(x,y_{s-k},\ldots,y_{s-2},z,y,y_1,\ldots,y_{s-k-1})$ is a path of length $s$, we have $(x,y_{s-k-1})\in A\Gamma$. Then $(x,y_{s-i})\in A\Gamma$  for any  $3\leq i\leq s$,  contrary to the fact that $(x,y_0)=(x,y)\in \Gamma_{2,s+(-1)^l}$.
\end{proof}

\begin{lemma}\label{2s1}
Let $p+q>s+1$.	Suppose  $(x,y)\in \Gamma_{2,s+1}$ and $z\in P_{(1,p),(1,q)}(x,y)$. Assume that $(y=y_0,y_1,\ldots,y_{s},y_{s+1}=x)$ is a shortest path. Then 	the following hold:
	\begin{enumerate}
		\item\label{2s1-2} $(x,y_1),(y_s,y)\in \Gamma_{1,s}$;
		
		\item\label{2s1-1}
		$(y_i,y_{i+1})\in \Gamma_{1,1}\cup \Gamma_{1,2}\cup \Gamma_{1,3}$, $(y_j,y_{j+2})\in \Gamma_{2,1}\cup\Gamma_{2,2}$, and $(y_{l+3},y_{l})\in \Gamma_{1,3}$ for any $0\leq i\leq s$, $0\leq j\leq s-1$ and $0\leq l\leq s-2$;

		\item\label{2s1-3} If $p\geq 4$, then $s$ is odd, $(y_{s-1},z)\in A\Gamma$ and $(x,y_3)\in \Gamma_{1,s-2}$.
	\end{enumerate}	
\end{lemma}
\begin{proof}
	If $z\in \{y_{1},\ldots,y_s\}$, then $z=y_q$ since $\partial(y,z)=q$, which implies $p=\partial(z,x)=\partial(y_q,x)=s+1-q<p$, a contradiction. Thus,  $z\notin \{y_{1},\ldots,y_s\}$.
	
	 (i) Note that $(y,y_1,\ldots,y_{s},x)$ is a shortest path. Since $\Gamma$ is $s$-quasi-transitive and $(y_1,\ldots,y_s,x)$, $(y,y_1,\ldots,y_s)$ are two paths of length $s$, \ref{2s1-2} holds.
	
 (ii) 	
Note that $(y_{l+3},y_{l+4},\ldots,y_s,y_{s+1},z,y_{0},y_1,\ldots,y_{l})$ is a path of length $s$ for any $0\leq l\leq s-2$. Since $(y_{l},y_{l+1},y_{l+2},y_{l+3})$ is a shortest path, we have $(y_{l+3},y_{l})\in \Gamma_{1,3}$. The fact that $(y_{l+3},y_{l},y_{l+1},y_{l+2})$ is a circuit implies $(y_i,y_{i+1})\in \Gamma_{1,1}\cup \Gamma_{1,2}\cup \Gamma_{1,3}$ for any $0\leq i\leq s$ and  $\partial(y_{j+2},y_j)\leq 2$ for any $0\leq j\leq s-1$. Since $(y_{0},y_1,\ldots,y_{s},y_{s+1})$ is a shortest path, one has $\partial(y_j,y_{j+2})=2$ for any $0\leq j\leq s-1$. Thus, \ref{2s1-1}  holds.

(iii) Since $4\leq p=\partial(z,x)\leq \partial(z,y_{s-1})+2$ and  $(z,y,y_1,\ldots,y_{s-1})$ is a path of length $s$, we get $(y_{s-1},z)\in A\Gamma$.
We prove $(x,y_{s-2m-2})\in A\Gamma$ by induction on $2m$ with  $0\leq 2m\leq s-2$.  If $2m=0$,
	from \ref{2s1-1}, then $(x,y_{s-2})=(y_{s+1},y_{s-2})\in \Gamma_{1,3}$. Assume that $(x,y_{s-2k-2})\in A\Gamma$ with $0\leq 2k<s-2$. Since $(x,y_{s-2k-2},\ldots,y_{s-1},z,y,y_1,\ldots,y_{s-2k-4})$ is a path of length $s$, one gets $(x,y_{s-2k-4})\in A\Gamma$. 
	If $s$ is even, then $(x,y)=(x,y_0)\in A\Gamma$, a contradiction. Thus, $s$ is odd. It follows that $(x,y_3)=(x,y_{s-2(\frac{s-5}{2})-2})\in A\Gamma$. Since $(y,y_1,\ldots,y_{s},x)$ is a shortest path,   one obtains $(x,y_3)\in \Gamma_{1,s-2}$.	
\end{proof}

\begin{lemma}\label{jichu3}
	Let $l,m\geq 2$ and $t\geq 1$. Suppose $\Gamma_{2,m}\in \Gamma_{1,l}\Gamma_{1,t}$. Then $\Gamma_{1,l}\notin \Gamma_{1,l}\Gamma_{1,p}$ for any  $p\geq s-2$.
\end{lemma}
\begin{proof}
	Assume the contrary, namely,  $\Gamma_{1,l}\in \Gamma_{1,l}\Gamma_{1,p}$ for some  $p\geq s-2$.
	Let $(x,y)\in \Gamma_{2,m}$ and $z_0\in P_{(1,l),(1,t)}(x,y)$. Since $\Gamma_{1,l}\in \Gamma_{1,l}\Gamma_{1,p}$, there exist distinct vertices $z_i\in P_{(1,l),(1,p)}(x,z_{i-1})$ for $1\leq i\leq p$.  Since $\partial(x,y)=2$, one has $y\notin \{z_0,z_1,\ldots,z_p\}$. Since $p\geq s-2$ and $(z_{p},z_{p-1},\ldots,z_1,z_0)$ is a walk consisting of arcs of type $(1,p)$,  $(x,z_{s-2},z_{s-3},\ldots,z_1,z_0,y)$ is a path of length $s$, which implies that $(x,y)$ or $(y,x)\in A\Gamma$, a contradiction.
\end{proof}

\begin{lemma}\label{oegy6}
	Let $\Gamma_{1,s-2}\Gamma_{1,s}=\{\Gamma_{2,s-1}\}$. Suppose $\Gamma_{l,1}\in \Gamma_{2,s-1}\Gamma_{s-2,1}$ for some $l>1$. Then $\Gamma_{2,s-1}\Gamma_{1,s}\subseteq \{\Gamma_{1,s-2},\Gamma_{1,s-1},\Gamma_{1,s}\}$.
\end{lemma}
\begin{proof}
	Let $(x,y)\in \Gamma_{2,s-1}$ and $(y,z)\in \Gamma_{1,s}$. Since $\Gamma_{l,1}\in \Gamma_{2,s-1}\Gamma_{s-2,1}$, from Lemma \ref{jb} \ref{jb-2}, we get $\Gamma_{2,s-1}\in \Gamma_{l,1}\Gamma_{1,s-2}$. Then there exists $w\in P_{(l,1),(1,s-2)}(x,y)$. Since $y\in P_{(1,s-2),(1,s)}(w,z)$ and $\Gamma_{1,s-2}\Gamma_{1,s}=\{\Gamma_{2,s-1}\}$, one has $(w,z)\in \Gamma_{2,s-1}$. The fact $s-1=\partial(y,x)\leq 1+\partial(z,x)$ implies that $\partial(z,x)\geq s-2$. Since $w\in P_{(l,1),(2,s-1)}(x,z)$,  by  the commutativity of $\Gamma$ and setting $r=2$ and $t=l$ in Lemma \ref{jichu0}, we get $\partial(x,z)=1$, $\partial(z,x)=1$ or $\partial(z,x)=s-l-1$. Since $l>1$ and $\partial(z,x)\geq s-2$, from Lemma \ref{maxt}, one has $(x,z)\in \Gamma_{1,s-2}\cup \Gamma_{1,s-1}\cup \Gamma_{1,s}$.
\end{proof}

The commutativity of $\Gamma$ will be used frequently in the sequel, so we no longer refer to it for the sake of simplicity.

\begin{lemma}\label{oegy5}
	Let $t$ and  $s$ be  positive integers  such that $t<s$ and $2\mid (t+s)$. Suppose $\{t+1\}\subseteq T\subseteq \{t+1,t+3,\ldots,s-1,s+1\}$.	 If $s$ is even, or  $s$ is odd and $t\geq 5$, then $\Gamma_{1,t}\notin \Gamma_{2,s-1}\Gamma_{s-2,1}$.
\end{lemma}
\begin{proof}
		Assume the contrary, namely, $\Gamma_{1,t}\in \Gamma_{2,s-1}\Gamma_{s-2,1}$. Lemma \ref{jb} \ref{jb-2} implies that
	$\Gamma_{2,s-1}\in \Gamma_{1,t}\Gamma_{1,s-2}$. Let $(x,z)\in \Gamma_{2,s-1}$ and $y\in P_{(1,t),(1,s-2)}(x,z)$. Pick two shortest paths $(z=z_0,z_1,\ldots,z_{s-3},y)$ and  $(y,y_1,\ldots,y_{t-1},x)$.
	Since $\partial(z,x)=s-1$, one has $x\notin \{z_1,\ldots,z_{s-3}\}$.
	 Since $\min T=t+1$ and $(y=y_0,y_1,\ldots,y_{t-1},y_{t}=x)$ is a circuit, we get $(y,y_{1})\in \Gamma_{1,t}$. 
     
     Note that $(z=z_0,z_1,\ldots,z_{s-3},y)$ is a shortest path. If $y_1\in \{z=z_0,z_1,\ldots,z_{s-3}\}$, then $y_1=z_{s-t-2}$ since $\partial(y_1,y)=t$, which implies that $s-1=\partial(z,x)\leq \partial(z,z_{s-t-2})+\partial(y_1,x)=s-3$, a contradiction. Thus, $y_1\notin \{z=z_0,z_1,\ldots,z_{s-3}\}$.
Suppose  $y_j=z_i$ for some $0\leq i\leq s-3$ and $2\leq j\leq t$. Since $t=\partial(y_1,y)\leq \partial(y_1,y_j)+\partial(z_i,y)=j-1+s-2-i$, one has $i-j\leq s-t-3$, and so $s-1=\partial(z,x)\leq  \partial(z,z_i)+\partial(y_j,x)=i+t-j\leq s-3$, a contradiction. Thus, $y_j\notin \{z=z_0,z_1,\ldots,z_{s-3}\}$ for any $1\leq j\leq t$.
	
	Suppose $t<4$. Then $s$ is even and $t=2$. Since $(z,z_1,\ldots,z_{s-3},y,y_1=y_{t-1},x)$ is a path of length $s$, we obtain either $(z,x)\in A\Gamma$ or $(x,z)\in A\Gamma$, a contradiction.
	Thus, $t\geq 4$.

 Since $(z,z_1,\ldots,z_{s-3},y,y_1,y_2)$ is a path of length $s$, one obtains either $(z,y_2)\in A\Gamma$ or $(y_2,z)\in A\Gamma$. If $(z,y_2)\in A\Gamma$, then $s-1=\partial(z,x)\leq \partial(z,y_2)+\partial(y_2,x)=t-1$ since $(y_2,\ldots,y_{t-1},x)$ is a shortest path, contrary to the fact that $t<s$. Then $(y_2,z)\in A\Gamma$. Since $(x,y,y_1,y_2,z,z_1,\ldots,z_{s-4})$  is a path of length $s$, we have $(x,z_{s-4})$ or $(z_{s-4},x)\in A\Gamma$.
	Since $s-1=\partial(z,x)\leq\partial(z,z_{s-4})+\partial(z_{s-4},x)$, one gets $(x,z_{s-4})\in A\Gamma$.

	We claim  that $(x,z_{s-2m})\in A\Gamma$  for any  $4\leq 2m\leq s$, where $z_0=z$. We prove it by induction on $2m$. If $2m=4$,
	 then  our claim is valid. Assume that $(x,z_{s-2k})\in A\Gamma$ with $4\leq 2k<s$. Since $s-1=\partial(z,x)\leq\partial(z,z_{s-2k-2})+\partial(z_{s-2k-2},x)=s-2k-2+\partial(z_{s-2k-2},x)$ and
	$(x,z_{s-2k},\ldots,z_{s-3},y,y_1,y_2,z,z_1,\ldots,z_{s-2k-2})$ is a path of length $s$, we get $(x,z_{s-2k-2})\in A\Gamma$. Thus, our claim is valid.

	If $s$ is even, from the claim, then $(x,z)=(x,z_0)\in A\Gamma$, contrary to the fact that $(x,z)\in \Gamma_{2,s-1}$.
	Then $s$ is odd and $t\geq 5$. By the claim, we get $(x,z_1)\in A\Gamma$.
	Since $(x,z_1,\ldots,z_{s-3},y,y_1,y_2)$ is a path of length $s$, we have either $(x,y_2)\in A\Gamma$ or $(y_2,x)\in A\Gamma$. Since $t\geq 5$ and $(y,y_1,\ldots,y_{t-1},y_{t}=x)$ is a circuit, one obtains $(x,y_2)\in \Gamma_{1,t-2}$, contrary to the fact that $\min T=t+1$.
\end{proof}

\begin{lemma}\label{jichu4}
	Suppose  $(2,s+1)\in \wz{\partial}(\Gamma)$. Then $\Gamma_{2,s-1}\notin \Gamma_{1,s-1}\Gamma_{1,s}$.
\end{lemma}
\begin{proof}
  Let $(x,y)\in \Gamma_{2,s+1}$. Pick a shortest path $(y,y_1,\ldots,y_{s},x)$.
	Since $5\leq s+1=\partial(y,x)\leq \partial(y,y_s)+1$ and $(y,y_1,\ldots,y_{s})$ is a path of length $s$, one has
	$(y_s,y)\in \Gamma_{1,s}$.
	Suppose  $\Gamma_{2,s-1}\in \Gamma_{1,s-1}\Gamma_{1,s}$. By Lemma \ref{jb} \ref{jb-2}, one gets $p_{(s-1,1),(2,s-1)}^{(1,s)}=\frac{k_{2,s-1}}{k_{1,s}}p_{(1,s-1),(1,s)}^{(2,s-1)}\ne 0$, which implies that there exists   $z\in P_{(2,s-1),(s-1,1)}(y_s,y)$.
	The fact $s+1=\partial(y,x)\leq 1+\partial(z,x)$ implies $\partial(z,x)\geq s$. Since $y_s\in P_{(n,1),(2,s-1)}(x,z)$ for some $n+1\in T$, by Lemma \ref{maxt} and setting $r=2$ and $t=n$ in Lemma \ref{jichu0}, we get  $(x,z)\in \Gamma_{1,s}$. Since $z\in P_{(1,s),(s-1,1)}(x,y)$, one obtains $\Gamma_{2,s+1}\in \Gamma_{1,s}\Gamma_{s-1,1}$, which implies $p_{(1,s-1),(s,1)}^{(s+1,2)}>0$ from Lemma \ref{jb} \ref{jb-0}, contrary to Lemma \ref{jichu0}. Thus, $\Gamma_{2,s-1}\notin \Gamma_{1,s-1}\Gamma_{1,s}$.
\end{proof}

\begin{lemma}\label{34lemma00}
	Suppose $\Gamma_{1,m}\Gamma_{1,s}=\{\Gamma_{2,s-1}\}$ for some $m\in \{s-1,s\}$. Then $\Gamma_{1,s}\notin \Gamma_{1,s}\Gamma_{1,t}$ for any $t+1\in T$.
\end{lemma}
\begin{proof}
	Let $(x,y)\in \Gamma_{2,s-1}$ and $z\in  P_{(1,m),(1,s)}(x,y)$.
	Assume the contrary, namely,  $\Gamma_{1,s}\in \Gamma_{1,s}\Gamma_{1,t}$ for some $t+1\in T$. Then there exists $w\in P_{(1,s),(1,t)}(z,y)$. Since
$z\in  P_{(1,m),(1,s)}(x,w)$,  one has $(x,w)\in \Gamma_{2,s-1}$. Since $y\in P_{(2,s-1),(t,1)}(x,w)$, we get $\Gamma_{2,s-1}\in \Gamma_{2,s-1}\Gamma_{t,1}$, contrary to Lemma \ref{jichu0}.
\end{proof}

\begin{lemma}\label{1205}
	Let  $\Gamma_{1,m}\in \Gamma_{1,p}\Gamma_{1,q}$ with $m\geq 2$ and $q\geq 3$. Then the following hold:
	\begin{enumerate}
		\item\label{1205-1} if $s\geq 4$, then $\Gamma_{2,s-1}\notin \Gamma_{1,s}\Gamma_{1,m}$;
		
		\item\label{1205-2}
		if $s\geq 5$, then $\Gamma_{2,s+1}\notin \Gamma_{1,s}\Gamma_{1,m}$.
	\end{enumerate}	
\end{lemma}
\begin{proof}
Let $t\in \{4,5\}$. Suppose for the contrary that $s\geq t$ and $\Gamma_{2,s-(-1)^t}\in \Gamma_{1,s}\Gamma_{1,m}$. Let $(x,y)\in \Gamma_{2,s-(-1)^t}$ and $z\in P_{(1,s),(1,m)}(x,y)$.   Since $\Gamma_{1,m}\in \Gamma_{1,p}\Gamma_{1,q}$, there exists $w\in P_{(1,p),(1,q)}(z,y)$. Pick a shortest path $(y=y_0,y_1,\ldots,y_{s-(-1)^t-1},x)$. 

If $t=4$,  by setting $l=1$  in  Lemma \ref{1207}, then $w\in \{y_1,\ldots,y_{s-(-1)^t-1}\}$ since $s+m\geq s+2>s-1=s+(-1)^l$ and $s\geq t=4>3=4+(-1)^l$. If $t=5$,  by setting $l=0$  in  Lemma \ref{1207}, then $w\in \{y_1,\ldots,y_{s-(-1)^t-1}\}$ since $s+m\geq s+2>s+1=s+(-1)^l$ and $s\geq t=5=4+(-1)^l$. Thus, $w\in \{y_1,\ldots,y_{s-(-1)^t-1}\}$.

 Since   $\partial(y,w)=q$, we get $w=y_q$, which implies that $s=\partial(z,x)\leq \partial(z,w)+\partial(y_q,x)=1+s-(-1)^t-q\leq s-(-1)^t-2\leq s+1-2=s-1$, a contradiction. 
\end{proof}

\begin{lemma}\label{jiougy1}
Suppose $\Gamma_{2,s+1}\notin \Gamma_{1,m}\Gamma_{1,s}$ for any $s-2\leq m\leq s$. Then $s\notin T$.
\end{lemma}
\begin{proof}
Assume the contrary, namely, $s\in T$. We divide the proof into two cases.

	\textbf{Case 1.} $\Gamma_{2,s-1}\in \Gamma_{1,s-1}\Gamma_{1,s}$.

 By Lemma \ref{jichu4}, we have $(2,s+1)\notin \wz{\partial}(\Gamma)$. 	Since $s\geq 4$, from Lemma \ref{jichu3}, we get $\Gamma_{1,s-1},\Gamma_{1,s}\notin \Gamma_{1,s-1}\Gamma_{1,s}$.  In view of Lemma \ref{even} \ref{even-1}  and  Lemmas \ref{maxt}, \ref{rem}, one has $\Gamma_{1,s-1}\Gamma_{1,s}=\{\Gamma_{2,s-1}\}$.

By  Lemma \ref{34lemma00},  one obtains $\Gamma_{1,s}\notin \Gamma_{1,s}\Gamma_{1,t}$ for any $t+1\in T$. Lemma \ref{jb} \ref{jb-2} implies $\Gamma_{1,s}\Gamma_{s,1}\cap \{\Gamma_{1,t},\Gamma_{t,1}\mid t+1\in T\}=\emptyset$. By setting $r=1$ and $t=s$ in  Lemma \ref{jichu0}, one has $\Gamma_{1,s-1} \Gamma_{s,1}\cap \Gamma_{1,s} \Gamma_{s,1}=\emptyset$, and so $\Gamma_{1,s}\Gamma_{s-1,1}\cap \Gamma_{1,s} \Gamma_{s,1}=\emptyset$. By Lemma \ref{jb} \ref{jb-2}, we get $\sum_{\tilde{t}\in\wz{\partial}(\Gamma)}p_{(s,1),(1,s)}^{\tilde{t}}p_{\tilde{t},(1,s-1)}^{(1,s)}=0$.
By setting $i=h^*=(1,s)$,   $l=(1,s-1)$ and $j=(1,s)$ in Lemma \ref{jb} \ref{jb-4}, we obtain $\sum_{\tilde{r}\in\wz{\partial}(\Gamma)}p_{(1,s),(1,s-1)}^{\tilde{r}}p_{(s,1),\tilde{r}}^{(1,s)}=0$.
Since $\Gamma_{1,s-1}\Gamma_{1,s}=\{\Gamma_{2,s-1}\}$, one has $p_{(s,1),(2,s-1)}^{(1,s)}=0$. It follows from Lemma \ref{jb} \ref{jb-2} that
$p_{(1,s),(1,s)}^{(2,s-1)}=0$. Since $\Gamma_{1,s}\notin \Gamma_{1,s}\Gamma_{1,t}$ for any $t+1\in T$, from Lemma \ref{even} \ref{even-1}  and  Lemmas \ref{maxt}, \ref{rem}, we get $\Gamma_{1,s}^2= \{\Gamma_{1,s-1}\}$. Since $s\geq 4$ and $\Gamma_{2,s-1}\in \Gamma_{1,s}\Gamma_{1,s-1}$, from Lemma \ref{1205} \ref{1205-1}, one has $\Gamma_{1,s-1}\notin \Gamma_{1,p}\Gamma_{1,q}$ for any $p\geq 1$ and $q\geq 3$, contrary to the fact that $\Gamma_{1,s-1}\in \Gamma_{1,s}^2$.

	\textbf{Case 2.} $\Gamma_{2,s-1}\notin \Gamma_{1,s-1}\Gamma_{1,s}$.

	Since $\Gamma_{2,s+1}\notin \Gamma_{1,s-1}\Gamma_{1,s}$, from Lemma \ref{even} \ref{even-1} and Lemmas \ref{maxt}, \ref{rem}, one obtains $\Gamma_{1,s-1}\Gamma_{1,s}\subseteq \{\Gamma_{1,s-1},\Gamma_{1,s}\}$. 	By Lemma \ref{jb} \ref{jb-2} and \ref{jb-3}, we obtain \begin{align}
	k_{1,s}=p_{(s,1),(1,s-1)}^{(1,s-1)}+p_{(s,1),(1,s)}^{(1,s-1)}=p_{(1,s),(1,s-1)}^{(1,s-1)}+p_{(1,s),(s,1)}^{(1,s-1)},\nonumber
	\end{align}
	which implies that $\Gamma_{1,s-1}\Gamma_{s,1}\subseteq \{\Gamma_{1,s-1},\Gamma_{s,1}\}$ and  $p_{(1,s),(1,s)}^{(1,s-1)}=p_{(1,s),(s-1,1)}^{(1,s-1)}=0$. By Lemma \ref{jb} \ref{jb-2}, one has $\Gamma_{1,s-1}\notin \Gamma_{1,s}^2$ and $\Gamma_{1,s}\notin \Gamma_{1,s-1}^2$.
	
	Since $\Gamma_{2,s+1}\notin \Gamma_{1,s}^2\cup \Gamma_{1,s-2}\Gamma_{1,s}$, from Lemma \ref{even} \ref{even-1} and Lemmas \ref{maxt}, \ref{rem},  one obtains that $\Gamma_{1,s}^2\subseteq \{\Gamma_{1,s},\Gamma_{2,s-1}\}$, $\Gamma_{1,s-2}\Gamma_{1,s}\subseteq \{\Gamma_{1,s-1},\Gamma_{1,s},\Gamma_{2,s-1}\}$ and $\Gamma_{1,s-1}^2\subseteq \{\Gamma_{1,s-2},\Gamma_{1,s-1},\Gamma_{2,s-2},\Gamma_{2,s-1},\Gamma_{2,s+1}\}.$ By Lemma \ref{one}, one gets $\Gamma_{2,s-1}\in \Gamma_{1,s}^2$. It follows from Lemma \ref{jichu3} that  $\Gamma_{1,s}\notin  \Gamma_{1,s}^2\cup \Gamma_{1,s}\Gamma_{1,s-1}\cup\Gamma_{1,s}\Gamma_{1,s-2}$. Then  $\Gamma_{1,s}^2=\{\Gamma_{2,s-1}\}$, $\Gamma_{1,s-2}\Gamma_{1,s}\subseteq \{\Gamma_{1,s-1},\Gamma_{2,s-1}\}$ and $\Gamma_{1,s-1}\Gamma_{1,s}=\{\Gamma_{1,s-1}\}$. By Lemma \ref{jichu3} again, we obtain $\Gamma_{2,m}\notin \Gamma_{1,s-1}^2$ for any $m\geq 2$, and so $\Gamma_{1,s-1}^2\subseteq\{\Gamma_{1,s-2},\Gamma_{1,s-1}\}.$
	Lemma \ref{one} implies that $\Gamma_{1,s-2}\in \Gamma_{1,s-1}^2$.
 Since $\Gamma_{1,s-1}\Gamma_{s,1}\subseteq \{\Gamma_{1,s-1},\Gamma_{s,1}\}$, one gets
 $p_{(1,s-1),(s,1)}^{(1,s-2)}=0$.  Lemma \ref{jb} \ref{jb-2} implies  $p_{(1,s-2),(1,s)}^{(1,s-1)}=0$, and so $\Gamma_{1,s-2}\Gamma_{1,s}=\{\Gamma_{2,s-1}\}$.
	
	Let $(x,y)\in \Gamma_{2,s-1}$ and $z\in P_{(1,s-2),(1,s)}(x,y)$.  Since $\Gamma_{1,s-2}\in \Gamma_{1,s-1}^2$, there exists  $w\in P_{(1,s-1),(1,s-1)}(x,z)$.
	Since $\Gamma_{1,s-1}\Gamma_{1,s}=\{\Gamma_{1,s-1}\}$ and $z\in P_{(1,s-1),(1,s)}(w,y)$, we get $(w,y)\in \Gamma_{1,s-1}$, and so $w\in  P_{(1,s-1),(1,s-1)}(x,y)$. Then $\Gamma_{2,s-1}\in\Gamma_{1,s-1}^2$, contrary to the fact that $\Gamma_{1,s-1}^2\subseteq\{\Gamma_{1,s-1},\Gamma_{1,s-2}\}.$
\end{proof}

\begin{lemma}\label{jichu6}
	Let $r\in \{0,1\}$ and $1\leq n\leq s-2-r$. Suppose $\{\Gamma_{2,s-r-1},\Gamma_{2,s+(-1)^r}\}\cap \Gamma_{1,n}\Gamma_{1,s-r}\ne \emptyset$. Then $n-r+3\in T$.
\end{lemma}
\begin{proof}
	Suppose $\Gamma_{2,m}\in \Gamma_{1,n}\Gamma_{1,s-r}$ for some $m\in \{s-r-1,s+(-1)^r\}$.  Pick  $(x,y)\in \Gamma_{2,m}$. Then there exist $z_1\in  P_{(1,n),(1,s-r)}(x,y)$ and $z_2\in  P_{(1,s-r),(1,n)}(x,y)$. Pick a shortest path $(y,y_1,\ldots,y_{s-r-1},z_1)$.  Since $\partial(y,x)=m\geq s-r-1$, one has $x\notin \{y_1,\ldots,y_{s-r-2}\}$.
We divide the proof into two cases.

\textbf{Case 1.}  $z_2\in \{y_1,\ldots,y_{s-r-1}\}$.

If $x=y_{s-r-1}$, then $z_2\notin \{y_1,\ldots,y_{s-r-1}\}$  since $s-r=\partial(z_2,x)=\partial(z_2,y_{s-r-1})$,   a contradiction. Thus, $x\ne y_{s-r-1}$, and so $x\notin \{y_1,\ldots,y_{s-r-1}\}$.
 Since $\partial(y,z_2)=n$ and $(y,y_1,\ldots,y_{s-r-1},z_1)$ is a shortest path, one obtains $z_2=y_n$. Since $\Gamma$ is $s$-quasi-transitive and $(x,z_1,y,y_1,\ldots,y_{s-2})$ is a path of length $s$, we have either $(x,y_{s-2})\in A\Gamma$ or $(y_{s-2},x)\in A\Gamma$. If $(y_{s-2},x)\in A\Gamma$, then $s-r=\partial(z_2,x)\leq \partial(y_n,y_{s-2})+1=s-2-n+1\leq s-2$, a contradiction. Thus, $(x,y_{s-2})\in A\Gamma$, and so $\partial(y_{s-2},x)\leq \partial(y_{s-2},z_1)+\partial(z_1,x)=n-r+2$. Since $s-r=\partial(z_2,x)\leq \partial(y_{n},y_{s-2})+\partial(y_{s-2},x)=s-2-n+\partial(y_{s-2},x)$, one has $\partial(y_{s-2},x)\geq n-r+2$, and so $(x,y_{s-2})\in \Gamma_{1,n-r+2}$. Then $n-r+3\in T$.

\textbf{Case 2.}  $z_2\notin \{y_1,\ldots,y_{s-r-1}\}$.

Pick a shortest path $(z_1,x_1,\ldots,x_{n-1},x)$. Since   $\partial(z_2,x)=s-r>n$, one has $z_2\notin \{z_1,x_1,\ldots,x_{n-1}\}$.
	If $x_i\in \{y=y_0,y_1,\ldots,y_{s-n-2}\}$ for some $1\leq i\leq n-1$, then $s-r=\partial(z_2,x)\leq 1+\partial(y,x_i)+n-i\leq s-2$,   a contradiction. Thus, $x_i\notin \{y=y_0,y_1,\ldots,y_{s-n-2}\}$ for any $1\leq i\leq n-1$. Note that $(y,y_1,\ldots,y_{s-r-1},z_1)$ is a shortest path. Since $(z_1,x_1,\ldots,x_{n-1},x,z_2,y,y_1,\ldots,y_{s-n-2})$ is a path of length $s$, we get $(z_1,y_{s-n-2})\in \Gamma_{1,n-r+2}$, and so $n-r+3\in T$.
\end{proof}

\begin{lemma}\label{oegy4}
	Let $m,n\in \{s-1,s+1\}$.  Suppose  $\Gamma_{2,m}\in \Gamma_{1,p}\Gamma_{1,s}$ with $p\geq 4$. Then 		$\Gamma_{2,s-1}\notin \Gamma_{2,m}\Gamma_{2,n}$.
\end{lemma}
\begin{proof}
	Assume the contrary, namely, $\Gamma_{2,s-1}\in \Gamma_{2,m}\Gamma_{2,n}$. Let $(x,y)\in \Gamma_{2,s-1}$ and $z\in P_{(2,m),(2,n)}(x,y)$.  Pick a shortest path $(y,y_1,\ldots,y_{s-2},x)$. Since $\partial(y,z)=n\in \{s-1,s+1\}$, we get $z\notin \{y_1,\ldots,y_{s-2}\}$. 

    Since $\Gamma_{2,m}\in \Gamma_{1,p}\Gamma_{1,s}$,  there exists $z_1\in P_{(1,p),(1,s)}(x,z)$. If $z_1\in \{y=y_0,\ldots,y_{s-2}\}$, then $s=\partial(z,z_1)\leq 2+\partial(y,z_1)\leq s$, which implies $z_1=y_{s-2}$, and so $p=1$, a contradiction. Thus, $z_1\notin \{y=y_0,\ldots,y_{s-2}\}$.
    Since $s\geq 4$, we may assume  $\Gamma_{2,n}\in \Gamma_{1,q_1}\Gamma_{1,q_2}$ with $q_1\geq 1$ and $q_2\geq 2$. Then there exists  $z_2\in P_{(1,q_1),(1,q_2)}(z,y)$. Since $q_2\geq 2$, one has $z_2\ne y_1$.
If $z_2\in \{y_2,\ldots,y_{s-2}\}$, then $s=\partial(z,z_1)\leq 1+\partial(z_2,x)+1\leq s-1$, a contradiction. Thus, $z_2\notin \{y_1,\ldots,y_{s-2}\}$.

Since $4\leq p=\partial(z_1,x)\leq \partial(z_1,y_{s-3})+2$ and $(z_1,z,z_2,y,y_1,\ldots,y_{s-3})$ is a path of length $s$, one gets $(y_{s-3},z_1)\in A\Gamma$. Since $3\leq s-1\leq n=\partial(y,z)\leq1+ \partial(y_1,z)$ and $(y_1,\ldots,y_{s-2},x,z_1,z)$ is a path of length $s$, we have $(z,y_1)\in A\Gamma$. It follows that $s=\partial(z,z_1)\leq \partial(z,y_1)+\partial(y_1,y_{s-3})+\partial(y_{s-3},z_1)=s-2$, a contradiction.
\end{proof}

\section{Cycle polar diameter $s$}

In  this section, $\Gamma$ always denotes  an $s$-quasi-transitive  commutative weakly distance-regular digraph of cycle polar diameter  $s$, where $s\geq 4$. If $\Gamma$ has diameter at least $s$, from Lemma \ref{even} \ref{even-1}, then  $s+1\in T$, a contradiction.
Thus, $\Gamma$ has diameter at most $s-1$.

The main result of this section is as follows.

\begin{prop}\label{bc4}
 If  $\Gamma$ has cycle polar diameter $s$, then   $\Gamma$ is isomorphic to the digraph in Theorem \ref{main thm} \ref{main thm-4}.
\end{prop}

To prove Proposition \ref{bc4}, we first establish several auxiliary lemmas.

\begin{lemma}\label{sec5-x-2}
	Let $m+1,n+1\in T$ with $m\geq 3$ and  $m+n>s-1$. Suppose $\Gamma_{1,n}\in \Gamma_{1,a}\Gamma_{1,b}$ for some $a+1,b+1\in T$ with $b\in \{s-2,s-1\}$. Then $\Gamma_{2,s-1}\notin \Gamma_{1,m}\Gamma_{1,n}$.
\end{lemma}
\begin{proof}
	Assume the contrary, namely, $\Gamma_{2,s-1}\in \Gamma_{1,m}\Gamma_{1,n}$. Let $(x,y)\in \Gamma_{2,s-1}$ and $z\in  P_{(1,m),(1,n)}(x,y)$.
	Pick a shortest path $(y=y_0,y_1,\ldots,y_{s-2},x)$. Since $\Gamma_{1,n}\in \Gamma_{1,a}\Gamma_{1,b}$, there exists $z_0\in P_{(1,a),(1,b)}(z,y)$. Since $m\geq 3$ and $m+n>s-1$, by setting $l=1$ in Lemma \ref{1207}, one has $z_0\in \{y_1,\ldots,y_{s-2}\}$. Since $s-2\leq b=\partial(y,z_0)\leq \partial(y,y_{s-2})=s-2$, we get $z_0=y_{s-2}$. It follows that $m=\partial(z,x)\leq\partial(z,z_0)+\partial(y_{s-2},x)=2$, a contradiction.
\end{proof}

\begin{lemma}\label{sec5-x-1}
	Let $a+1\in T$. Suppose $\Gamma_{1,s-2}\in \Gamma_{1,a}\Gamma_{1,s-1}$. Then  $\Gamma_{2,s-1}\notin \Gamma_{1,s-2}^2$.
\end{lemma}
\begin{proof}
	If $s\geq 5$, from Lemma \ref{sec5-x-2},  then $\Gamma_{2,s-1}\notin \Gamma_{1,s-2}^2$. Now we consider the case that $s=4$.
	Assume the contrary, namely,  $\Gamma_{2,3}\in \Gamma_{1,2}^2$. Let $(x,y)\in \Gamma_{2,3}$ and $z\in  P_{(1,2),(1,2)}(x,y)$. Since $\Gamma_{1,2}\in \Gamma_{1,a}\Gamma_{1,3}$, there exist $z_1\in P_{(1,a),(1,3)}(x,z)$ and $z_2\in P_{(1,a),(1,3)}(z,y)$. Since $(x,z_1,z,z_2,y)$ is a path of length $4$, we get $(x,y)$ or $(y,x)\in A\Gamma$, a contradiction.
\end{proof}

A family of paths in $\Gamma$
is said to be \emph{internally-disjoint} if no vertex in $\Gamma$ is an internal vertex of more than one path of the family.

\begin{lemma}\label{sec5-1}
	Let  $(x,y_0)\in \Gamma_{2,s-2}$. Suppose that $(y_0,y_1,\ldots,y_{s-3},x)$ and $(x,z,z_1,y_0)$ are two internally-disjoint paths. Assume $\partial(z,x)\geq 3$. If  $z_2\in  P_{(1,p),(1,q)}(z,z_1)\cup P_{(1,q),(1,p)}(z_1,y_0)$ for some $p+1,q+1\in T$ with $q>1$, then $z_2\in \{y_1,\ldots,y_{s-3}\}$.
\end{lemma}
\begin{proof}
	Assume the contrary, namely,  $z_2\notin \{y_1,\ldots,y_{s-3}\}$. Since $q>1$, we get $z_2\notin \{z,z_1,y_0\}$. Since  $3\leq\partial(z,x)\leq \partial(z,y_{s-3})+1$, and $(z,z_2,z_1,y_0,y_1,\ldots,y_{s-3})$ or $(z,z_1,z_2,y_0,y_1,\ldots,y_{s-3})$ is a path of length $s$,  one obtains $(y_{s-3},z)\in A\Gamma$.

	If $x=z_2\in P_{(1,p),(1,q)}(z,z_1)\cup P_{(1,q),(1,p)}(z_1,y_0)$, then $\partial(z,x)=1$ or $\partial(x,y_0)=1$, a contradiction. Then $x\ne z_2$. Observe that $(y_0,y_1,\ldots,y_{s-3},x)$ is a shortest path. Since $(x,z,z_2,z_1,y_0,y_1,\ldots,y_{s-4})$ or $(x,z,z_1,z_2,y_0,y_1,\ldots,y_{s-4})$ is a path of length $s$, one has $(x,y_{s-4})\in A\Gamma$.
	
	We prove $(x,y_{s-i})\in A\Gamma$ for any  $4\leq i\leq s$ by induction on $i$.
	If $i=4$,
	then  $(x,y_{s-4})\in A\Gamma$. Now suppose  $(x,y_{s-k})\in A\Gamma$ with $4\leq k<s$. Note that $(y_0,y_1,\ldots,y_{s-3},x)$ is a shortest path. Since
	$(x,y_{s-k},\ldots,y_{s-3},z,z_1,y_0,y_1,\ldots,y_{s-k-1})$ is a path of length $s$, we have $(x,y_{s-k-1})\in A\Gamma$. Then $(x,y_{s-i})\in A\Gamma$ for any  $4\leq i\leq s$, contrary to the fact that $(x,y_0)\in  \Gamma_{2,s-2}$.
\end{proof}

\begin{lemma}\label{sec5-1'}
	If $\Gamma_{2,s-2}\in \Gamma_{1,s-1}^2\cup \Gamma_{1,s-2}\Gamma_{1,s-1}$, then $k_{1,s-1}=1$ and $\Gamma_{1,s-1}\notin \Gamma_{1,s-1}\Gamma_{1,q}$ for any $q+1\in T$.
\end{lemma}
\begin{proof}
	Assume the contrary, namely, $\Gamma_{1,s-1}\in \Gamma_{1,s-1}\Gamma_{1,q}$ for some $q+1\in T$. Let $(x,y_0)\in \Gamma_{2,s-2}$ and $z\in P_{(1,m),(1,s-1)}(x,y_0)$ for some $m\in \{s-1,s-2\}$. Pick a shortest path $(y_0,y_1,\ldots,y_{s-3},x)$.
	Since $\partial(y_0,z)=s-1$, one has $z\notin \{y_1,\ldots,y_{s-3}\}$.
	Since $\Gamma_{1,s-1}\in \Gamma_{1,s-1}\Gamma_{1,q}$,  there exist $z_1\in P_{(1,s-1),(1,q)}(z,y_0)$ and $z_2\in P_{(1,q),(1,s-1)}(z,z_1)$. Since $\partial(z_1,z_2)=s-1$ and $\partial(z_1,y_0)=1$, we get $y_0\notin\{z_1,z_2\}$.

	Suppose $z_1$ or $z_2\in \{y_1,\ldots,y_{s-3}\}$. Since $s-2\leq m=\partial(z,x)\leq \min\{1+\partial(z_1,x),1+\partial(z_2,x)\}$, one has $\partial(z_1,x),\partial(z_2,x)\geq s-3$, which implies $z_1=y_1$ or $z_2=y_1$. It follows that $s-1=\partial(z_1,z)\leq \partial(y_1,x)+1=s-2$ or $s-1=\partial(z_1,z_2)\leq \partial(z_1,y_0)+\partial(y_0,y_1)=2$, a contradiction. Thus, $z_1,z_2\notin \{y_1,\ldots,y_{s-3}\}$.
	
	Note that $(y_0,y_1,\ldots,y_{s-3},x)$ and $(x,z,z_1,y_0)$ are two internally-disjoint paths. Since $z_2\in P_{(1,q),(1,s-1)}(z,z_1)$ and  $z_2\notin \{y_1,\ldots,y_{s-3}\}$, from Lemma \ref{sec5-1}, one has $\partial(z,x)\leq 2$, which implies  $s-2\leq m=\partial(z,x)\leq 2$, and so $s=4$. Since  $(x,z,z_2,z_1,y_0)$ is a path of length $4$, one obtains either $(x,y_0)\in A\Gamma$ or $(y_0,x)\in A\Gamma$, a contradiction.
	Thus, $\Gamma_{1,s-1}\notin \Gamma_{1,s-1}\Gamma_{1,q}$ for any $q+1\in T$. It follows from Lemma \ref{k_{1,s-1}=1} that  $k_{1,s-1}=1$.
\end{proof}

\begin{lemma}\label{bc1}
	We have  $s-1\in T$. 
\end{lemma}
\begin{proof}
		Assume the contrary, namely,  $s-1\notin T$. 
	Since $\Gamma$ has diameter at most $s-1$,  we get
	\begin{align}
		\Gamma_{1,a}\Gamma_{1,s-1}&\subseteq \{\Gamma_{1,s-1},\Gamma_{2,s-2},\Gamma_{2,s-1}\},\label{buc-2}\tag{4.1}
	\end{align}
	where $a+1\in T$.
	By Lemma \ref{one}, one has $\Gamma_{1,s-1}^2 \cap \{\Gamma_{2,s-2},\Gamma_{2,s-1}\}\ne \emptyset$. Since $s\geq 4$, from Lemma \ref{jichu3}, we obtain $\Gamma_{1,s-1}\notin \Gamma_{1,s-1}^2$, and so $\Gamma_{1,s-1}^2\subseteq \{\Gamma_{2,s-2},\Gamma_{2,s-1}\}$.
	
	If $\Gamma_{2,s-2}\in \Gamma_{1,s-1}^2$, from Lemma \ref{sec5-1'}, then   $\Gamma_{1,s-1}\notin \Gamma_{1,s-1}\Gamma_{1,a}$ for any $a+1\in T$. If $\Gamma_{2,s-1}\in \Gamma_{1,s-1}^2$, from Lemma  \ref{sec5-x-2}, then $\Gamma_{1,s-1}\notin \Gamma_{1,a}\Gamma_{1,s-1}$ for any $a+1\in T$.
	Thus,  $\Gamma_{1,s-1}\notin \Gamma_{1,s-1}\Gamma_{1,a}$ for any $a+1\in T$. It follows from Lemma \ref{k_{1,s-1}=1} that  $k_{1,s-1}=1$.
	In view of  \eqref{buc-2}, one has $\Gamma_{1,a}\Gamma_{1,s-1}\subseteq \{\Gamma_{2,s-2},\Gamma_{2,s-1}\}$ for any $a+1\in T$.

	If $a+1\in T$ for some positive  integer $a<s-1$, by setting $r=1$ in Lemma \ref{jichu6}, then $s-1\in T$ since $\Gamma_{2,s-2}$ or $\Gamma_{2,s-1}\in \Gamma_{1,a}\Gamma_{1,s-1}$, a contradiction. The fact that $\Gamma$ has cycle polar diameter $s$  implies
	$T=\{s\}$.
	
	Since $k_{1,s-1}=1$, $\Gamma$ has valency $1$, a contradiction. 
\end{proof}

 Since $\Gamma$ has diameter at most $s-1$, from Lemma \ref{bc1},  one has
\begin{align}
	\Gamma_{1,q}\Gamma_{1,s-1}&\subseteq \{\Gamma_{1,s-2},\Gamma_{1,s-1},\Gamma_{2,s-2},\Gamma_{2,s-1}\},\label{buc-1}\tag{4.2}\\
	\Gamma_{1,s-2}^2&\subseteq \{\Gamma_{1,s-3},\Gamma_{1,s-2},\Gamma_{1,s-1},\Gamma_{2,s-3},\Gamma_{2,s-2},\Gamma_{2,s-1}\},\label{buc-x}\tag{4.3}
\end{align}
where $q+1\in T$.

\begin{lemma}\label{bc2}
	We have $\Gamma_{2,s-1}\notin \Gamma_{1,s-2}\Gamma_{1,s-1}$.
\end{lemma}
\begin{proof}
	Assume the contrary, namely, $\Gamma_{2,s-1}\in \Gamma_{1,s-2}\Gamma_{1,s-1}$. Let $(x,y)\in \Gamma_{2,s-1}$, $z\in  P_{(1,s-2),(1,s-1)}(x,y)$ and $z'\in P_{(1,s-1),(1,s-2)}(x,y)$. Pick a shortest path $(y,y_1,\ldots,y_{s-2},x)$.
	Since $\Gamma_{2,s-1}\in \Gamma_{1,s-2}\Gamma_{1,s-1}$ and $s\geq 4$, from Lemma \ref{jichu3}, we have $\Gamma_{1,s-1},\Gamma_{1,s-2}\notin  \Gamma_{1,s-2}\Gamma_{1,s-1}$. \eqref{buc-1} implies that $\Gamma_{1,s-2}\Gamma_{1,s-1}\subseteq \{\Gamma_{2,s-2},\Gamma_{2,s-1}\}$.
	
	Suppose  $\Gamma_{1,s-1}\in \Gamma_{1,s-1}\Gamma_{1,q}$ for some $q+1\in T$. By Lemma \ref{sec5-1'}, we get $\Gamma_{1,s-2}\Gamma_{1,s-1}=\{\Gamma_{2,s-1}\}$. Pick a vertex $w\in P_{(1,s-1),(1,q)}(z,y)$. Since $z\in P_{(1,s-2),(1,s-1)}(x,w)$, we have $(x,w)\in\Gamma_{2,s-1}$. Since $y\in P_{(2,s-1),(q,1)}(x,w)$ and $s\geq4$, by setting $r=2$ and $t=q$ in Lemma \ref{jichu0}, one obtains $s-1=\partial(w,x)=s-q-1$, and so $q=0$, a contradiction. Then $\Gamma_{1,s-1}\notin \Gamma_{1,s-1}\Gamma_{1,q}$ for any $q+1\in T$. By Lemma \ref{k_{1,s-1}=1}, we get $\Gamma_{1,s-1}\Gamma_{s-1,1}=\{\Gamma_{0,0}\}$ and  $k_{1,s-1}=1$.

	Since $\Gamma_{2,s-1}\in \Gamma_{1,s-2}\Gamma_{1,s-1}$, from Lemma \ref{comm0}, one has $\Gamma_{1,s-2}\Gamma_{1,s-1}=\{\Gamma_{2,s-1}\}$ and $\Gamma_{2,s-1}\Gamma_{s-1,1}=\{\Gamma_{1,s-2}\}$. By Lemma \ref{jb} \ref{jb-2}, we get $p_{(1,q),(1,s-1)}^{(2,s-1)}=\frac{k_{1,q}}{k_{2,s-1}}p_{(2,s-1),(s-1,1)}^{(1,q)}=0$ for any  $q\ne s-2$. Since $\Gamma_{1,s-1}\notin \Gamma_{1,s-1}\Gamma_{1,q}$ for any $q+1\in T\setminus\{s-1\}$, from \eqref{buc-1}, one obtains
	\begin{align}
		\Gamma_{1,q}\Gamma_{1,s-1}&\subseteq \{\Gamma_{1,s-2},\Gamma_{2,s-2}\}~{\rm for~ any}~q\ne s-2.\label{buc-y}\tag{4.4}
	\end{align}
	
	Suppose $\Gamma_{1,s-2}\in \Gamma_{1,q}\Gamma_{1,s-1}$ for some $q+1\in T$.  Let $w'\in P_{(1,q),(1,s-1)}(z',y)$. Since $z'\in P_{(1,s-1),(1,s-2)}(x,y)$ and $s\geq 4$, by setting $l=1$ in Lemma \ref{1207},  we get $w'\in \{y_1,\ldots,y_{s-2}\}$, which implies that $s-1=\partial(y,w')\leq \partial(y,y_{s-2})=s-2$, a contradiction. Thus,  $\Gamma_{1,s-2}\notin \Gamma_{1,q}\Gamma_{1,s-1}$ for any $q+1\in T$. By \eqref{buc-y}, we get $\Gamma_{1,q}\Gamma_{1,s-1}=\{\Gamma_{2,s-2}\}$  for any $q+1\in T$ and $q\ne s-2.$  Since $s\in T$, we get $\Gamma_{1,s-1}^2= \{\Gamma_{2,s-2}\}$. Since $k_{1,s-1}=1$, from Lemma \ref{comm0}, one has $\Gamma_{2,s-2}\Gamma_{s-1,1}=\{\Gamma_{1,s-1}\}$.  If $q+1\in T$ for some $1\leq q\leq s-3$, from Lemma \ref{jb} \ref{jb-2}, then  $\Gamma_{1,q}\in \Gamma_{2,s-2}\Gamma_{s-1,1}$ since $\Gamma_{1,q}\Gamma_{1,s-1}=\{\Gamma_{2,s-2}\}$, a contradiction. Thus,  $q+1\notin T$ for any $1\leq q\leq s-3$, and so $T=\{s-1,s\}$.

	Since $\Gamma_{1,s-1}^2=\{\Gamma_{2,s-2}\}$, $\Gamma_{1,s-1}\Gamma_{s-1,1}=\{\Gamma_{0,0}\}$ and $\Gamma_{1,s-1}\Gamma_{1,s-2}=\{\Gamma_{2,s-1}\}$, we get $p_{(1,s-1),(1,s-1)}^{(1,s-2)}=p_{(1,s-1),(s-1,1)}^{(1,s-2)}=p_{(1,s-1),(1,s-2)}^{(1,s-2)}=0$. Lemma \ref{jb} \ref{jb-2} implies that $p_{(1,s-2),(s-1,1)}^{(1,s-1)}=p_{(s-1,1),(1,s-2)}^{(s-1,1)}=p_{(s-1,1),(1,s-2)}^{(1,s-2)}=0$, and so
	$\Gamma_{1,s-1},\Gamma_{s-1,1},\Gamma_{1,s-2}\notin \Gamma_{1,s-2}\Gamma_{s-1,1}$. Since $T=\{s-1,s\}$, by Lemma \ref{jb} \ref{jb-0} and setting $r=1$ and $t=s-2$ in  Lemma \ref{jichu0}, one has $\Gamma_{1,s-2}\Gamma_{s-1,1}=\{\Gamma_{s-2,1}\}$.
	
	Lemma \ref{jb} \ref{jb-0} and \ref{jb-2} imply
	$p_{(1,s-2),(1,s-2)}^{(1,s-1)}=p_{(s-2,1),(s-2,1)}^{(s-1,1)}=\frac{k_{1,s-2}}{k_{1,s-1}}p_{(1,s-2),(s-1,1)}^{(s-2,1)}\ne 0$, and so
	$\Gamma_{1,s-1}\in \Gamma_{1,s-2}^2$.
	Then there exists $w_0'\in P_{(1,s-2),(1,s-2)}(z,y)$. If $s\geq 5$, by setting $l=1$ in Lemma \ref{1207}, then $w_0'=y_{s-2}$ since $\partial(y,w_0')=s-2$ and  $(y,y_1,\ldots,y_{s-2},x)$ is a shortest path,   which implies that $s-2=\partial(z,x)\leq 1+\partial(w_0',x)=2$, a contradiction. Thus,  $s=4$.
	
	Note that $\Gamma_{1,3}^2=\{\Gamma_{2,2}\}$. Let $(x_0,y_0)\in \Gamma_{2,2}$ and $x_1\in  P_{(1,3),(1,3)}(x_0,y_0)$.
	Since 	$\Gamma_{1,3}\in \Gamma_{1,2}^2$, there exist $x_2\in P_{(1,2),(1,2)}(x_0,x_1)$ and $x_3\in P_{(1,2),(1,2)}(x_1,y_0)$. Since $(x_0,x_2,x_1,x_3,y_0)$ is a path of length $4$, one has either $(x_0,y_0)\in A\Gamma$ or $(y_0,x_0)\in A\Gamma$, a contradiction.
\end{proof}

\begin{lemma}\label{1213}
	Suppose  $T=\{q+1,s-1,s\}$ with $q\leq s-3$. Let $(u,v)\in \Gamma_{1,s-1}$ and  $(v,v_1,\ldots,v_{s-2},u)$ be a shortest path. If $k_{1,s-1}=1$ and $k_{1,q}\geq 2$, then there exists  $u_0 \in \Gamma_{q,1}(v)$ such that $(u,u_0)\in \Gamma_{1,s-2}\cup\Gamma_{1,q}$. 
\end{lemma}
\begin{proof}
	Since $k_{1,q}\geq 2$, by setting $l=s-1$ and $p=q$ in  Lemma \ref{1009}, there exists $u_0 \in \Gamma_{q,1}(v)$ such that $u_0\notin \{v_1,\ldots,v_{s-2}\}$.  Since $(u_0,v,v_1,\ldots,v_{s-2},u)$ is a path of length $s$ and $s-1=\partial(v,u)\leq q+\partial(u_0,u)\leq s-3+\partial(u_0,u)$, one has $(u,u_0)\in A\Gamma$. Since $v\in \Gamma_{1,s-1}(u)$  and $k_{1,s-1}=1$, one obtains $(u,u_0)\notin \Gamma_{1,s-1}$. The fact that $T=\{q+1,s-1,s\}$ implies $(u,u_0)\in \Gamma_{1,s-2}\cup\Gamma_{1,q}$.
\end{proof}

\begin{lemma}\label{bc3}
	We have $\Gamma_{2,s-1}\notin \Gamma_{1,s-1}^2$.
\end{lemma}
\begin{proof}
	Assume the contrary, namely, $\Gamma_{2,s-1}\in \Gamma_{1,s-1}^2$. Let $(x,y)\in \Gamma_{2,s-1}$ and $z\in  P_{(1,s-1),(1,s-1)}(x,y)$. Pick a shortest path $(y,y_1,\ldots,y_{s-2},x)$.

	Since $s\geq 4$, from Lemma  \ref{sec5-x-2}, we get  $\Gamma_{1,s-1}\notin \Gamma_{1,a}\Gamma_{1,s-1}\cup \Gamma_{1,a}\Gamma_{1,s-2}$ for any $a+1\in T$, and so $\Gamma_{1,s-1}\notin \Gamma_{1,s-2}^2$.
	By Lemma \ref{k_{1,s-1}=1}, one gets $\Gamma_{1,s-1}\Gamma_{s-1,1}=\{\Gamma_{0,0}\}$ and  $k_{1,s-1}=1$. Since $\Gamma_{2,s-1}\in \Gamma_{1,s-1}^2$, from Lemma \ref{comm0}, we get
	$\Gamma_{1,s-1}^2=\{\Gamma_{2,s-1}\}$ and $\Gamma_{2,s-1}\Gamma_{s-1,1}=\{\Gamma_{1,s-1}\}$, which implies  $p_{(1,q),(1,s-1)}^{(2,s-1)}=\frac{k_{1,q}}{k_{2,s-1}}p_{(2,s-1),(s-1,1)}^{(1,q)}=0$ for any $q\leq s-2$. Since $\Gamma_{1,s-1}\notin \Gamma_{1,q}\Gamma_{1,s-1}$ for any $q+1\in T$ with $q\leq s-2$, from
	\eqref{buc-1}, one has
	\begin{align}
		\Gamma_{1,q}\Gamma_{1,s-1}&\subseteq \{\Gamma_{1,s-2},\Gamma_{2,s-2}\}.\label{buc-q}\tag{4.5}
	\end{align}
	Since $k_{1,s-1}=1$, from Lemma \ref{jb} \ref{jb-5}, we obtain $\Gamma_{1,s-2}\Gamma_{1,s-1}=\{\Gamma_{1,s-2}\}$ or $\{\Gamma_{2,s-2}\}$. We divide the proof into two cases.
	
	\textbf{Case 1.} $\Gamma_{1,s-2}\Gamma_{1,s-1}=\{\Gamma_{1,s-2}\}$.
	
	Since  $k_{1,s-1}=1$, from Lemma \ref{comm0}, we have $\Gamma_{1,s-2}\Gamma_{s-1,1}=\{\Gamma_{1,s-2}\}$, which implies $p_{(1,q),(1,s-1)}^{(1,s-2)}=\frac{k_{1,q}}{k_{1,s-2}}p_{(1,s-2),(s-1,1)}^{(1,q)}=0$ for any $q\leq s-3$. In view of \eqref{buc-q}, one gets
	\begin{align}
		\Gamma_{1,q}\Gamma_{1,s-1}&\subseteq \{\Gamma_{2,s-2}\}~{\rm for~ any}~q\leq s-3.\label{buc-w}\tag{4.6}
	\end{align}
	Since $\Gamma_{1,s-2}\Gamma_{1,s-1}=\{\Gamma_{1,s-2}\}$, from  Lemma \ref{jichu3}, we have $\Gamma_{2,m}\notin \Gamma_{1,s-2}^2$ for any $m\geq 2$. Since $\Gamma_{1,s-1}\notin \Gamma_{1,s-2}^2$, from \eqref{buc-x}, one gets $s\geq 5$ and $\Gamma_{1,s-2}^2\subseteq \{\Gamma_{1,s-3},\Gamma_{1,s-2}\}$, or $s=4$ and $\Gamma_{1,2}^2\subseteq \{\Gamma_{1,1},\Gamma_{1,2},\Gamma_{2,1}\}$.
	
	Suppose $\Gamma_{1,s-3}\in \Gamma_{1,s-2}^2$. By \eqref{buc-w}, we get $\Gamma_{1,s-3}\Gamma_{1,s-1}= \{\Gamma_{2,s-2}\}$.
	Let $(x',y')\in \Gamma_{2,s-2}$ and $z'\in  P_{(1,s-3),(1,s-1)}(x',y')$. Since $\Gamma_{1,s-3}\in \Gamma_{1,s-2}^2$, there exists a vertex $w'\in  P_{(1,s-2),(1,s-2)}(x',z')$. Since $\Gamma_{1,s-2}\Gamma_{1,s-1}=\{\Gamma_{1,s-2}\}$, one has $(w',y')\in \Gamma_{1,s-2}$, which implies that $w'\in  P_{(1,s-2),(1,s-2)}(x',y')$, and so $\Gamma_{2,s-2}\in \Gamma_{1,s-2}^2$, a contradiction.
	By Lemma \ref{one}, we get $s=4$ and  $\{\Gamma_{2,1}\}\subseteq \Gamma_{1,2}^2\subseteq \{\Gamma_{2,1},\Gamma_{1,2}\}$.
	
	Since $\Gamma_{1,2}\Gamma_{1,3}=\{\Gamma_{1,2}\}$, from Lemma \ref{jb} \ref{jb-2}, one has $\Gamma_{1,3}\in \Gamma_{1,2}\Gamma_{2,1}$, which implies that there exists $w_1\in P_{(1,2),(2,1)}(x,z)$. Since $\Gamma_{2,1}\in \Gamma_{1,2}^2$, there exists $w_2\in P_{(1,2),(1,2)}(w_1,z)$. Since $(x,w_1,w_2,z,y)$ is a path of length $4$, we get $(x,y)$ or $(y,x)\in A\Gamma$, contrary to the fact that $(x,y)\in \Gamma_{2,3}$.
	
	\textbf{Case 2.} $\Gamma_{1,s-2}\Gamma_{1,s-1}=\{\Gamma_{2,s-2}\}$.
	
	Since  $k_{1,s-1}=1$,  from Lemma \ref{comm0}, one has  $\Gamma_{2,s-2}\Gamma_{s-1,1}=\{\Gamma_{1,s-2}\}$, which implies $p_{(1,q),(1,s-1)}^{(2,s-2)}=\frac{k_{1,q}}{k_{2,s-2}}p_{(2,s-2),(s-1,1)}^{(1,q)}=0$ for any $q\leq s-3$. By \eqref{buc-q}, we get
	\begin{align}
		\Gamma_{1,q}\Gamma_{1,s-1}&\subseteq \{\Gamma_{1,s-2}\}~{\rm for~ any}~q\leq s-3.\label{buc-r}\tag{4.7}
	\end{align}

	Since $\Gamma_{1,s-1}\notin \Gamma_{1,s-2}^2$ and $\Gamma_{1,s-2}\notin \Gamma_{1,s-2}\Gamma_{1,s-1}\cup \Gamma_{1,s-1}\Gamma_{s-1,1}\cup \Gamma_{1,s-1}^2$, 
	we get $p_{(1,s-2),(1,s-2)}^{(1,s-1)}=p_{(1,s-2),(1,s-1)}^{(1,s-2)}=p_{(s-1,1),(1,s-1)}^{(s-2,1)}=p_{(1,s-1),(1,s-1)}^{(1,s-2)}=0$. By Lemma \ref{jb} \ref{jb-2}, one has $p_{(1,s-1),(s-2,1)}^{(1,s-2)}=p_{(1,s-2),(s-1,1)}^{(1,s-2)}=p_{(1,s-1),(s-2,1)}^{(1,s-1)}=p_{(s-1,1),(1,s-2)}^{(1,s-1)}=0$. It follows  that $\Gamma_{1,s-2},\Gamma_{s-2,1},\Gamma_{1,s-1},\Gamma_{s-1,1}\notin \Gamma_{1,s-1}\Gamma_{s-2,1}$.
	By setting $r=1$ and $t=s-2$ in Lemma \ref{jichu0}, one has  $\Gamma_{1,s-1}\Gamma_{s-2,1}\subseteq\{\Gamma_{1,q},\Gamma_{q,1}\mid q+1\in T~{\rm and}~q\leq s-3\}$. Since $k_{1,s-1}=1$, from Lemma \ref{jb} \ref{jb-5}, one gets
	$\Gamma_{1,s-1}\Gamma_{s-2,1}=\{\Gamma_{1,q}\}$ or $\{\Gamma_{q,1}\}$ for some $q \leq s-3$.
	
	Let $p+1\in T$ with $p\leq s-3$. By \eqref{buc-r}, we have $\Gamma_{1,p}\Gamma_{1,s-1}=\{\Gamma_{1,s-2}\}$. It follows from Lemma \ref{jb} \ref{jb-0} and \ref{jb-2} that $\Gamma_{p,1}\in \Gamma_{s-2,1}\Gamma_{1,s-1}$. Since $\Gamma_{1,s-1}\Gamma_{s-2,1}=\{\Gamma_{1,q}\}$ or $\{\Gamma_{q,1}\}$ for some $q \leq s-3$, one gets $p=q$. Since $p+1\in T$ with $p\leq s-3$ was arbitrary, we obtain $T=\{q+1,s-1,s\}$  and $\Gamma_{1,s-1}\Gamma_{s-2,1}=\{\Gamma_{q,1}\}.$
	
	If $k_{1,q}=1$, from Lemma \ref{jb} \ref{jb-1}, then $k_{1,s-2}=1$ since $\Gamma_{1,s-1}\Gamma_{s-2,1}=\{\Gamma_{q,1}\}$, which implies that $\Gamma$ has valency $3$ since $T=\{q+1,s-1,s\}$, a contradiction. Then $k_{1,q}\geq 2$.
	Let $(u,v)\in \Gamma_{1,s-1}$. Pick a shortest path $(v,v_1,\ldots,v_{s-2},u)$.
	By Lemma \ref{1213}, there exists $u_0 \in \Gamma_{q,1}(v)$ such that $(u,u_0)\in \Gamma_{1,s-2}\cup\Gamma_{1,q}$.
	
	Suppose $(u,u_0)\in \Gamma_{1,s-2}$. Since $u\in  P_{(s-2,1),(1,s-1)}(u_0,v)$, we get $\Gamma_{1,q}\in \Gamma_{s-2,1}\Gamma_{1,s-1}$. The fact $\Gamma_{1,s-1}\Gamma_{s-2,1}=\{\Gamma_{q,1}\}$ implies $q=1$. It follows from Lemma \ref{jb} \ref{jb-2} that $\Gamma_{1,s-1}\in \Gamma_{1,1}\Gamma_{1,s-2}$. Pick a vertex $z_0\in P_{(1,1),(1,s-2)}(z,y)$. By setting $l=1$ in Lemma \ref{1207}, one obtains $z_0\in\{y_1,\ldots,y_{s-2}\}$.
	Since $\partial(y,z_0)=s-2$ and  $(y,y_1,\ldots,y_{s-2},x)$ is a shortest path,  we get $z_0=y_{s-2}$, and so $s-1=\partial(z,x)\leq \partial(z,z_0)+\partial(y_{s-2},x)=2$, a contradiction.
	
	Suppose $(u,u_0)\in \Gamma_{1,q}$.
	Since $u_0\in P_{(1,q),(1,q)}(u,v)$, one has $\Gamma_{1,s-1}\in \Gamma_{1,q}^2$. Pick a vertex $z_1\in P_{(1,q),(1,q)}(z,y)$. By  setting $l=1$ in Lemma \ref{1207}, one has $z_1\in\{y_1,\ldots,y_{s-2}\}$. Since $s-1=\partial(z,x)\leq 1+\partial(z_1,x)$ and $(y,y_1,\ldots,y_{s-2},x)$ is a shortest path, we get $z_1=y_1$. Since $(z_1,y)\in\Gamma_{1,q}$, one obtains $q=1$, which implies that $s-1=\partial(y,z)\leq \partial(y,y_1)+\partial(z_1,z)=2$, a contradiction.
\end{proof}

In view of Lemmas \ref{bc2} and  \ref{bc3}, from \eqref{buc-1}, one obtains
\begin{align}
	\Gamma_{1,s-1}^2,\Gamma_{1,s-2}\Gamma_{1,s-1}&\subseteq \{\Gamma_{1,s-2},\Gamma_{1,s-1},\Gamma_{2,s-2}\}.\label{buc-5}\tag{4.8}
\end{align}

\begin{lemma}\label{bc5}
	We have $\Gamma_{2,s-2}\notin \Gamma_{1,s-2}\Gamma_{1,s-1}$.
\end{lemma}
\begin{proof}
	Assume the contrary, namely, $\Gamma_{2,s-2}\in \Gamma_{1,s-2}\Gamma_{1,s-1}$.
	By   Lemma \ref{jichu3}, one has $\Gamma_{1,s-2},\Gamma_{1,s-1}\notin  \Gamma_{1,s-1}\Gamma_{1,s-2}$. It follows from \eqref{buc-5} that    $\Gamma_{1,s-2}\Gamma_{1,s-1}=\{\Gamma_{2,s-2}\}$. By Lemma \ref{sec5-1'},  we get $k_{1,s-1}=1$ and $\Gamma_{1,s-1}\notin \Gamma_{1,s-1}\Gamma_{1,q}$ for any $q+1\in T$. Since $\Gamma_{1,s-2}\Gamma_{1,s-1}=\{\Gamma_{2,s-2}\}$, from Lemma \ref{comm0}, we get $\Gamma_{2,s-2}\Gamma_{s-1,1}=\{\Gamma_{1,s-2}\}$. By Lemma \ref{jb} \ref{jb-2}, one has $p_{(1,s-1),(1,s-1)}^{(2,s-2)}=\frac{k_{1,s-1}}{k_{2,s-2}}p_{(2,s-2),(s-1,1)}^{(1,s-1)}=0$. Since  $\Gamma_{1,s-1}\notin \Gamma_{1,s-1}^2$, from \eqref{buc-5}, one obtains  $\Gamma_{1,s-1}^2=\{\Gamma_{1,s-2}\}$.
	
	Since  $k_{1,s-1}=1$, from Lemma \ref{jb} \ref{jb-1} and Lemma \ref{comm0}, one has $k_{1,s-2}=1$  and $\Gamma_{1,s-2}\Gamma_{s-1,1}=\{\Gamma_{1,s-1}\}$.	
	Note that $\Gamma_{2,s-2}\Gamma_{s-1,1}=\{\Gamma_{1,s-2}\}$. Then $p_{(1,s-2),(s-1,1)}^{(1,q)}=p_{(2,s-2),(s-1,1)}^{(1,q)}=0$ for any $q \leq s-3$.
	By Lemma \ref{jb} \ref{jb-2}, we get $p_{(1,q),(1,s-1)}^{(1,s-2)}=p_{(1,q),(1,s-1)}^{(2,s-2)}=0$ for any $q \leq s-3$. Since $\Gamma_{1,s-1}\notin \Gamma_{1,s-1}\Gamma_{1,q}$ for any $q \leq s-3$, from  \eqref{buc-1}, one obtains $\Gamma_{1,q}\Gamma_{1,s-1}=\{\Gamma_{2,s-1}\}$ for any $q+1\in T$ with $ q\leq s-3$.
	
	Since $\Gamma$ has valency greater than $3$ and $k_{1,s-2}=k_{1,s-1}=1$, there exists $q+1\in T$ with $q\leq s-3$.


	
	For each $q'+1\in T$ with $q'\leq s-3$, since $\Gamma_{1,q'}\Gamma_{1,s-1}=\{\Gamma_{2,s-1}\}$ and $k_{1,s-1}=1$, from Lemma \ref{comm0}, we obtain $\Gamma_{2,s-1}\Gamma_{s-1,1}=\{\Gamma_{1,q'}\}$. Since $q'+1\in T$ with $q'\leq s-3$ was arbitrary, one gets $T=\{q+1,s-1,s\}$, $\Gamma_{1,q}\Gamma_{1,s-1}=\{\Gamma_{2,s-1}\}$ and $\Gamma_{2,s-1}\Gamma_{s-1,1}=\{\Gamma_{1,q}\}$.
	
	By setting $r=1$ in Lemma \ref{jichu6}, we get $q+2\in T$, which implies $q+2\in \{s-1,s\}$. Since $q\leq s-3$, we get $q=s-3$, and so $T=\{s-2,s-1,s\}$.
	If $k_{1,s-3}=1$, then  $\Gamma$ has valency $3$, a contradiction. Then $k_{1,s-3}\geq 2$.

	Let $(u,v)\in \Gamma_{1,s-1}$ and $(v,v_1,\ldots,v_{s-2},u)$ be a shortest path.  Since $T=\{s-2,s-1,s\}$, from  Lemma \ref{1213}, there exists $u_0\in \Gamma_{s-3,1}(v)$ such that   $(u,u_0)\in \Gamma_{1,s-2}\cup\Gamma_{1,s-3}$.
	If $(u,u_0)\in \Gamma_{1,s-2}$, then $\Gamma_{s-3,1}\in \Gamma_{s-1,1}\Gamma_{1,s-2}$ since $u\in P_{(s-1,1),(1,s-2)}(v,u_0)$, contrary to the fact that $\Gamma_{1,s-2}\Gamma_{s-1,1}=\{\Gamma_{1,s-1}\}$. Then  $(u,u_0)\in \Gamma_{1,s-3}$.
	
	Since $u_0\in P_{(1,s-3),(1,s-3)}(u,v)$, we get $\Gamma_{1,s-1}\in \Gamma_{1,s-3}^2$.
	Observe that $\Gamma_{2,s-2}\in \Gamma_{1,s-2}\Gamma_{1,s-1}$. Pick $(x,y)\in \Gamma_{2,s-2}$ and  $z\in  P_{(1,s-1),(1,s-2)}(x,y)$. Let  $(y,y_1,\ldots,y_{s-3},x)$ be a shortest path.
	Since $\Gamma_{1,s-2}\in \Gamma_{1,s-1}^2$ and $\Gamma_{1,s-1}\in \Gamma_{1,s-3}^2$, there exist $z_1\in P_{(1,s-1),(1,s-1)}(z,y)$ and $z_2\in P_{(1,s-3),(1,s-3)}(z_1,y)$.
	Since $\partial(z,x)=\partial(y,z_1)=s-1$, we get $z,z_1\notin \{y_1,\ldots,y_{s-3}\}$, which implies  that $(y,y_1,\ldots,y_{s-3},x)$ and $(x,z,z_1,y)$ are two internally-disjoint paths.
	By Lemma \ref{sec5-1}, one has $z_2\in \{y_1,\ldots,y_{s-3}\}$. The fact $\partial(y,z_2)=s-3$ implies $z_2=y_{s-3}$, and so $s-1=\partial(z,x)\leq \partial(z,z_2)+\partial(y_{s-3},x)\leq 3$. Then $s=4$. Since $(x,z,z_1,z_2,y)$ is a path of length $4$, we obtain either $(x,y)\in A\Gamma$ or $(y,x)\in A\Gamma$, a contradiction.
\end{proof}

For a nonempty closed subset $F$ of $\{\Gamma_{\wz i}\}_{\wz i\in \tilde{\partial}(\Gamma)}$,
let
$$V\Gamma/F :=\{F(x)\mid x\in V\Gamma\} ~{\rm and}~ \Gamma_{\tilde{i}}^F :=\{(F(x),F(y))\mid y\in F\Gamma_{\tilde{i}}F(x)\},$$
where $F(x):=\{y\in V\Gamma\mid (x,y)\in \underset{f\in F}{\cup}f\}$. The digraph $(V\Gamma/F,\underset{(1,s)\in \tilde{\partial}(\Gamma)}{\cup}\Gamma_{1,s}^F)$ is said to be the  {\em quotient digraph} of $\Gamma$ over $F$, denoted by $\Gamma/F$.

\begin{proof}[Proof of Proposition \ref{bc4}]
	 Lemma \ref{bc5} and \eqref{buc-5} imply $\Gamma_{1,s-2}\Gamma_{1,s-1}\subseteq \{\Gamma_{1,s-2},\Gamma_{1,s-1}\}$. In view of Lemma \ref{jb} \ref{jb-2} and  \ref{jb-3}, one gets
	$k_{1,s-1}=p_{(s-1,1),(1,s-2)}^{(1,s-2)}+p_{(s-1,1),(1,s-1)}^{(1,s-2)}=p_{(1,s-1),(1,s-2)}^{(1,s-2)}+p_{(1,s-1),(s-1,1)}^{(1,s-2)}$, which implies  $p_{(1,s-1),(1,s-1)}^{(1,s-2)}=p_{(1,s-1),(s-2,1)}^{(1,s-2)}=0$.
	By 	Lemma \ref{jb} \ref{jb-2} again, we get $p_{(1,s-2),(1,s-2)}^{(1,s-1)}=0$.
	\eqref{buc-x} and \eqref{buc-5} imply	$\Gamma_{1,s-2}^2\subseteq \{\Gamma_{1,s-3},\Gamma_{1,s-2},\Gamma_{2,s-3},\Gamma_{2,s-2},\Gamma_{2,s-1}\}$ and $\Gamma_{1,s-1}^2\subseteq \{\Gamma_{1,s-1},\Gamma_{2,s-2}\}$. By Lemma \ref{one}, we obtain $\Gamma_{2,s-2}\in \Gamma_{1,s-1}^2$.
	By Lemma \ref{sec5-1'}, one has  $k_{1,s-1}=1$ and $\Gamma_{1,s-1}\notin \Gamma_{1,s-1}\Gamma_{1,q}$ for any $q+1\in T$, which imply
	$\Gamma_{1,s-1}^2=\{\Gamma_{2,s-2}\}~{\rm and}~ \Gamma_{1,s-2}\Gamma_{1,s-1}=\{\Gamma_{1,s-2}\}$. \mbox{Since} $k_{1,s-1}=1$, from  Lemma \ref{comm0}, we get
	$\Gamma_{2,s-2}\Gamma_{s-1,1}=\{\Gamma_{1,s-1}\}~{\rm and}~\Gamma_{1,s-2}\Gamma_{s-1,1}=\{\Gamma_{1,s-2}\}$,
	which imply $p_{(2,s-2),(s-1,1)}^{(1,q)}=p_{(1,s-2),(s-1,1)}^{(1,q)}=0$ for any $q\leq s-3$. By Lemma \ref{jb} \ref{jb-2}, one has  $p_{(1,q),(1,s-1)}^{(2,s-2)}=p_{(1,q),(1,s-1)}^{(1,s-2)}=0$. Since	$\Gamma_{1,s-1}\notin \Gamma_{1,s-1}\Gamma_{1,q}$  for any $q+1\in T$ with $q\leq s-3$,  from \eqref{buc-1}, we get  $\Gamma_{1,q}\Gamma_{1,s-1}=\{\Gamma_{2,s-1}\}$.
	
	Since $\Gamma_{1,s-2}\Gamma_{1,s-1}=\{\Gamma_{1,s-2}\}$, from  Lemma \ref{jichu3}, one has $\Gamma_{2,m}\notin \Gamma_{1,s-2}\Gamma_{1,q}$ for any $q\geq 1$ and $m\geq 2$, which implies that $s\geq 5$ and $\Gamma_{1,s-2}^2\subseteq \{\Gamma_{1,s-3},\Gamma_{1,s-2}\}$, or $s=4$ and $\Gamma_{1,2}^2\subseteq \{\Gamma_{1,1},\Gamma_{1,2},\Gamma_{2,1}\}$.
	If $\Gamma_{1,s-3}\in \Gamma_{1,s-2}^2$, then $\Gamma_{1,s-1}\Gamma_{1,s-3}=\{\Gamma_{2,s-1}\}$,  contrary to Lemma \ref{sec5-x-2}. 
	By Lemma \ref{one}, we get $s=4$ and $\{\Gamma_{2,1}\}\subseteq  \Gamma_{1,2}^2\subseteq \{\Gamma_{2,1},\Gamma_{1,2}\}$.
	
	Suppose $2\in T$. Then $\Gamma_{1,1}\Gamma_{1,3}=\{\Gamma_{2,3}\}$.  Let $(u,v)\in \Gamma_{2,3}$ and $u_0\in  P_{(1,1),(1,3)}(u,v)$. Pick a shortest path $(v,v_1,v_2,u_0)$. If $u\in \{v_1,v_2\}$, then $u=v_2$ since $\partial(u,u_0)=1$, which implies  $3=\partial(v,u)=\partial(v,v_2)=2$, a contradiction. If $u\notin \{v_1,v_2\}$, then $(u,v)$ or $(v,u)\in A\Gamma$ since $(v,v_1,v_2,u_0,u)$ is a path of length $4$, a contradiction. Thus, $2\notin T$, and so $T=\{3,4\}$.
	
	Note that $\Gamma_{1,3}^2=\{\Gamma_{2,2}\}$ and $\Gamma_{2,2}\Gamma_{3,1}=\{\Gamma_{1,3}\}$. By Lemma \ref{jb} \ref{jb-0}, one gets $\Gamma_{1,3}^3=\Gamma_{2,2}\Gamma_{1,3}=\{\Gamma_{3,1}\}$. The fact  $k_{1,3}=1$ implies $\Gamma_{1,3}^4=\Gamma_{1,3}\Gamma_{3,1}=\{\Gamma_{0,0}\}$. It follows that $\Gamma_{1,3}^{4i}=\{\Gamma_{0,0}\}$ and $\Gamma_{1,3}^{4i+r}=\{\Gamma_{r,4-r}\}$ for $i\geq0$ and $r\in\{1,2,3\}$. Then $F:=\langle\Gamma_{1,3}\rangle=\{\Gamma_{0,0},\Gamma_{1,3},\Gamma_{3,1},\Gamma_{2,2}\}$. Since $k_{1,3}=1$, $(F(x),(F(x))^2\cap\Gamma_{1,3})$ is isomorphic to $C_4$ for any $x\in V\Gamma$.
	
	Let $(x,y)\in \Gamma_{2,2}$. Pick a vertex $z\in \Gamma_{1,2}(y)$. Since $\Gamma_{1,3}^2=\{\Gamma_{2,2}\}$, there exists $w\in P_{(1,3),(1,3)}(x,y)$.  Since $\Gamma_{1,2}\Gamma_{1,3}=\{\Gamma_{1,2}\}$, one has $(w,z)\in \Gamma_{1,2}$, and so $(x,z)\in \Gamma_{1,2}$. Since $(x,y)\in \Gamma_{2,2}$ and $z\in \Gamma_{1,2}(y)$ were arbitrary,  we get $\Gamma_{2,2}\Gamma_{1,2}=\{\Gamma_{1,2}\}$.
	
	Suppose $\Gamma_{3,3}\in \Gamma_{1,2}\Gamma_{2,1}$. Then there exist $(x',y')\in \Gamma_{3,3}$ and $z'\in P_{(1,2),(2,1)}(x',y')$. Since $\Gamma_{1,2}\Gamma_{1,3}=\{\Gamma_{1,2}\}$ and $\Gamma_{2,1}\in \Gamma_{1,2}^2$, there exist  $x_1\in P_{(1,2),(1,3)}(x',z')$ and  $x_2\in P_{(1,2),(1,2)}(z',y')$, which imply that $(x',x_1,z',x_2,y')$ is a path of length $4$, and so $(x',y')$ or $(y',x')\in A\Gamma$, a contradiction. Thus, $\Gamma_{3,3}\notin \Gamma_{1,2}\Gamma_{2,1}$.
	
	Since  $\Gamma_{1,3}^2=\{\Gamma_{2,2}\}$, $\Gamma_{1,2}\Gamma_{1,3}=\{\Gamma_{1,2}\}$ and $\Gamma_{1,2}^2\subseteq \{\Gamma_{1,2},\Gamma_{2,1}\}$, one has $(2,3)\notin \tilde{\partial}(\Gamma)$. It follows that $\Gamma_{1,2}\Gamma_{2,1}\subseteq \{\Gamma_{0,0},\Gamma_{1,2},\Gamma_{2,1},\Gamma_{1,3},\Gamma_{3,1},\Gamma_{2,2}\}$.
	Combining this with $\Gamma_{1,2}\Gamma_{3,1}=\Gamma_{2,2}\Gamma_{1,2}=\{\Gamma_{1,2}\}$, we obtain $\langle\Gamma_{1,2}\rangle\subseteq\{\Gamma_{0,0},\Gamma_{1,2},\Gamma_{2,1},\Gamma_{1,3},\Gamma_{3,1},\Gamma_{2,2}\}$. Since $T=\{3,4\}$, one has $\tilde{\partial}(\Gamma)=\{(0,0),(1,2),(2,1),(1,3),(3,1),(2,2)\}$.

	We claim that $(u',v')\in\Gamma_{1,2}$ for $u'\in F(u)$ and $v'\in F(v)$ when $(u,v)\in\Gamma_{1,2}$. Without loss of generality, we may assume $u'\neq u$. Since $F=\{\Gamma_{0,0},\Gamma_{1,3},\Gamma_{3,1},\Gamma_{2,2}\}$, we have $u\in P_{(1,3),(1,2)}(u',v)\cup  P_{(3,1),(1,2)}(u',v)\cup  P_{(2,2),(1,2)}(u',v)$. The fact $\Gamma_{1,3}\Gamma_{1,2}=\Gamma_{1,2}\Gamma_{3,1}=\Gamma_{2,2}\Gamma_{1,2}=\{\Gamma_{1,2}\}$ implies $(u',v)\in\Gamma_{1,2}$. If $v=v'$, then our claim is valid. If $v\neq v'$, then $v\in P_{(1,2),(1,3)}(u',v')\cup  P_{(1,2),(3,1)}(u',v')\cup  P_{(1,2),(2,2)}(u',v')$, and so $(u',v')\in\Gamma_{1,2}$. Thus, our claim is valid.
	
	Set $\Sigma:=\Gamma/F$. By  the claim, $\Gamma$ is isomorphic to $\Sigma\circ C_4$, which implies $\tilde{\partial}_{\Gamma}(x_0,y_0)=\tilde{\partial}_{\Sigma}(F(x_0),F(y_0))$ for vertices $x_0,y_0$ with  $y_0\notin F(x_0)$.
	Since $\tilde{\partial}(\Gamma)=\{(0,0),(1,2),(2,1),(1,3),(3,1),(2,2)\}$, we get  $\tilde{\partial}(\Sigma)=\{(0,0),(1,2),(2,1)\}$.
	Let $\wz{h}\in\wz{\partial}(\Sigma)$ and $(F(u_1),F(v_1))\in\Sigma_{\wz{h}}$. Without loss of generality, we may assume that $(u_1,v_1)\in \Gamma_{\wz{h}}$. Let $\wz{i},\wz{j}\in\wz{\partial}(\Sigma)$. If $\wz{i}$ or $\wz{j}=(0,0)$, then $|\Sigma_{\wz{i}}(F(u_1))\cap\Sigma_{\wz{j}^*}(F(v_1))|=\delta_{\wz{h},\wz{j}}$ or $\delta_{\wz{h},\wz{i}}$.
	If $\wz{i},\wz{j}\ne(0,0)$, then
	$$P_{\wz{i},\wz{j}}(u_1,v_1)=\cup_{w\in P_{\wz{i},\wz{j}}(u_1,v_1)}F(w)
	=\cup_{F(w)\in \Sigma_{\wz{i}}(F(u_1))\cap\Sigma_{\wz{j}^*}(F(v_1))}F(w).$$
	Since $(F(x),(F(x))^2\cap\Gamma_{1,3})$ is isomorphic to $C_4$ for any $x\in V\Gamma$,  one gets $|F(w)|=4$, which implies $|\Sigma_{\wz{i}}(F(u_1))\cap\Sigma_{\wz{j}^*}(F(v_1))|=\frac{p_{\wz{i},\wz{j}}^{\wz{h}}}{4}.$ Since $\Gamma$ is a commutative weakly distance-regular digraph, $\Sigma$ is a  commutative weakly distance-regular digraph. Thus, the desired result follows.
\end{proof}

\section{Cycle polar diameter $s+1$}

In  this section, $\Gamma$ always denotes  $\Gamma$ always denotes  an $s$-quasi-transitive  commutative weakly distance-regular digraph of cycle polar diameter  $s+1$, where $s\geq 4$. The main result of this section is as follows.

\begin{prop}\label{propxj}
	If $\Gamma$  has  cycle polar diameter  $s+1$, then $\Gamma$ is isomorphic to the digraph in Theorem \ref{main thm} \ref{main thm-3}.
\end{prop}

We first prove that $T=\{s+1\}$, considering the cases of even and odd $s$ separately.


\subsection{$s$ is even}

The main result of this subsection is as follows.

\begin{prop}\label{prop4}
If $s$ is even, then	$T=\{s+1\}$.
\end{prop}

To prove Proposition \ref{prop4}, we need two auxiliary lemmas. In the remainder of this subsection, we always assume that $s$ is even.

\begin{lemma}\label{7.3xj}
If $\Gamma_{2,s+1}\in \Gamma_{1,m}\Gamma_{1,n}$ for some $m+1,n+1\in T$, then $m+n=s+1$.
\end{lemma}
\begin{proof}
Since $\Gamma_{2,s+1}\in \Gamma_{1,m}\Gamma_{1,n}$, we have $m+n\geq s+1$. Assume the contrary, namely, $m+n> s+1$. Let $(x,y)\in \Gamma_{2,s+1}$. Pick a shortest path $(y=y_0,y_1,\ldots,y_{s},y_{s+1}=x)$.
	Suppose $s=4$. By Lemma \ref{2s1} \ref{2s1-2} and \ref{2s1-1}, we have $(x,y_1),(y_3,y)=(y_3,y_0)\in A\Gamma$, which imply that $(x,y_1,y_2,y_3,y)$ is a path of length $4$, and so $(x,y)$ or $(y,x)\in A\Gamma$, a contradiction.
	 Since $s$ is even, one gets  $s\geq 6$. By Lemma \ref{2s1} \ref{2s1-3} and the commutativity of $\Gamma$, one has    $m,n\leq 3$, and so $s+1<m+n\leq 6$, a contradiction.
\end{proof}

In view of Lemma \ref{even} \ref{even-1} and Lemmas \ref{maxt}, \ref{rem},  \ref{7.3xj}, one obtains
\begin{align}
\Gamma_{1,1}\Gamma_{1,s}&\subseteq \{\Gamma_{1,s-1},\Gamma_{1,s},\Gamma_{2,s-1},\Gamma_{2,s+1}\},\label{34gs-0}\tag{5.1}\\
\Gamma_{1,q}\Gamma_{1,s}&\subseteq \{\Gamma_{1,s-1},\Gamma_{1,s},\Gamma_{2,s-1}\},\label{34gs-1}\tag{5.2}\\
\Gamma_{1,s-1}^2&\subseteq \{\Gamma_{1,s-2},\Gamma_{1,s-1},\Gamma_{1,s},\Gamma_{2,s-2},\Gamma_{2,s-1}\},\label{34gs-2}\tag{5.3}
\end{align}
where $q+1\in T$ and $q\geq 2$.

\begin{lemma}\label{prop1}
The following hold:
	\begin{enumerate}
		\item\label{prop1-1}
			$T\subseteq \{3,5,\ldots,s-1,s+1\}$;
		
		\item\label{prop1-2} $\Gamma_{1,q}\Gamma_{1,s}=\{\Gamma_{2,s-1}\}$ for any $q+1\in T$.
	\end{enumerate}	
\end{lemma}
\begin{proof}
	Since $s\geq 4$, from Lemma \ref{7.3xj}, one has $\Gamma_{2,s+1}\notin \Gamma_{1,m}\Gamma_{1,s}$ for $s-2\leq m\leq s$. It follows from  Lemma \ref{jiougy1} that $s\notin T$. \eqref{34gs-1} implies $\Gamma_{1,s}^2\subseteq \{\Gamma_{1,s},\Gamma_{2,s-1}\}$. By Lemma \ref{one}, we get $\Gamma_{2,s-1}\in \Gamma_{1,s}^2$. Lemma \ref{jichu3} implies  $\Gamma_{1,s}\notin \Gamma_{1,s}^2$, and so $\Gamma_{1,s}^2=\{\Gamma_{2,s-1}\}$.
	By Lemma \ref{34lemma00}, one has $\Gamma_{1,s}\notin \Gamma_{1,s}\Gamma_{1,t}$ for any $t+1\in T$. Since $s\notin T$, from \eqref{34gs-0} and \eqref{34gs-1},
	we get $\Gamma_{1,1}\Gamma_{1,s}\subseteq \{\Gamma_{2,s-1},\Gamma_{2,s+1}\}$ and  $\Gamma_{1,q}\Gamma_{1,s}=\{\Gamma_{2,s-1}\}$ for any $q+1\in T$ and $q\geq 2$.

If $p+1\in T$ for some positive odd integer $p$, by setting $r=0$ in Lemma \ref{jichu6}, then $s\in T$ since $\Gamma_{2,s-1}$ or $\Gamma_{2,s+1}\in \Gamma_{1,p}\Gamma_{1,s}$, a contradiction. Thus, \ref{prop1-1}  holds. Since $2\notin T$, \ref{prop1-2}  holds.
\end{proof}

\begin{proof}[Proof of Proposition \ref{prop4}]
	Suppose for the contrary that $T\ne \{s+1\}$. Let $t+1=\min T$.
	By Lemma \ref{prop1} \ref{prop1-1},  $t$ is  even  and $t\leq s-2$.  In view of Lemma \ref{prop1} \ref{prop1-2}, we have $\Gamma_{1,t}\Gamma_{1,s}=\{\Gamma_{2,s-1}\}$. By setting $r=0$ in Lemma \ref{jichu6},  one has $t+3\in T$. Repeating this process, since $\max T=s+1$, we get $T=\{t+1,t+3,t+5,\ldots,s-1,s+1\}$.

	 Suppose  $\Gamma_{p,1}\in \Gamma_{2,s-1}\Gamma_{s,1}$ for some $p+1\in T$.  Lemma \ref{jb} \ref{jb-2} implies that $\Gamma_{2,s-1}\in \Gamma_{1,s}\Gamma_{p,1}$. Let $(u,v)\in \Gamma_{2,s-1}$ and $v_0\in  P_{(1,s),(p,1)}(u,v)$. Since $\Gamma_{1,s}^2=\{\Gamma_{2,s-1}\}$ by Lemma \ref{prop1} \ref{prop1-2}, there exists $u_1\in  P_{(1,s),(1,s)}(u,v)$. Then $u_1\ne v_0$. Pick a shortest path $(v_0,v_1,\ldots,v_{s-1},u)$. If $v\in \{v_0,v_1,\ldots,v_{s-1}\}$,  then $v=v_1$ since $\partial(v,u)=s-1$, which implies $(v,v_0)\in \Gamma_{1,1}$, contrary to Lemma \ref{prop1} \ref{prop1-1}. Thus, $v\notin \{v_0,v_1,\ldots,v_{s-1}\}$.
	Since $\partial(u_1,u)=s$ and  $(v_0,v_1,\ldots,v_{s-1},u)$ is a shortest path, one obtains $u_1\notin \{v_0,\ldots,v_{s-1}\}$.
	Since $\Gamma$ is $s$-quasi-transitive and $(u,u_1,v,v_0,v_1,\ldots,v_{s-3})$ is a path of length $s$, we have $(u,v_{s-3})\in \Gamma_{1,3}$, contrary to Lemma \ref{prop1} \ref{prop1-1}. Thus,  $\Gamma_{p,1}\notin \Gamma_{2,s-1}\Gamma_{s,1}$ for any $p+1\in T$.

By setting $r=2$ and $t=s$ in Lemma \ref{jichu0}, one gets $\Gamma_{2,s-1}\Gamma_{s,1}\subseteq\{\Gamma_{1,p}\mid p+1\in T\}$. By  Lemma \ref{comm} and Lemma \ref{prop1} \ref{prop1-2}, we get $p_{(2,s-1),(q,1)}^{(1,s)}=k_{1,q}$ for any $q+1\in T$.
	In view of Lemma \ref{jb} \ref{jb-0} and  \ref{jb-3}, one has
	\begin{align}
	k_{2,s-1}= p_{(2,s-1),(s,1)}^{(1,s)}+\ldots+p_{(2,s-1),(t,1)}^{(1,s)}=k_{1,s}+\dots+k_{1,t}.\label{5gs-4}\tag{5.4}
	\end{align}

	Since $\Gamma_{1,s}^2=\{\Gamma_{2,s-1}\}$ by Lemma \ref{prop1} \ref{prop1-2}, from Lemma \ref{oegy4}, one has  $\Gamma_{2,s-1}\notin \Gamma_{2,s-1}^2$. Then $\sum_{\tilde{t}\in\wz{\partial}(\Gamma)}p_{(1,s),(1,s)}^{\tilde{t}}p_{\tilde{t},(2,s-1)}^{(2,s-1)}=0$. By setting $i=h=(1,s)$ and $l=j=(2,s-1)$ in Lemma \ref{jb} \ref{jb-4}, from Lemma \ref{prop1} \ref{prop1-2}, one obtains
	  $p_{(1,s),(2,s-1)}^{(1,p)}p_{(1,s),(1,p)}^{(2,s-1)}=0$ for any $p+1\in T$.
	Then   $p_{(1,s),(2,s-1)}^{(1,p)}=0$, and so $\Gamma_{1,p}\notin \Gamma_{2,s-1}\Gamma_{1,s}$ for any $p+1\in T$.

	Suppose $\Gamma_{p,1}\in \Gamma_{2,s-1}\Gamma_{s-2,1}$ for some $p+1\in T$. By Lemma \ref{prop1} \ref{prop1-1}, one has $p\geq 2$. Since $\Gamma_{1,s-2}\Gamma_{1,s}=\{\Gamma_{2,s-1}\}$ by Lemma \ref{prop1} \ref{prop1-2},
from Lemma \ref{oegy6}, we have $\Gamma_{2,s-1}\Gamma_{1,s}\subseteq \{\Gamma_{1,s-2},\Gamma_{1,s-1},\Gamma_{1,s}\}$, contrary to the fact that $\Gamma_{1,p}\notin \Gamma_{2,s-1}\Gamma_{1,s}$ for any $p+1\in T$.
	 Thus, $\Gamma_{p,1}\notin \Gamma_{2,s-1}\Gamma_{s-2,1}$ for any $p+1\in T$.

By setting $r=2$ and $t=s-2$ in Lemma \ref{jichu0}, we get $\Gamma_{2,s-1}\Gamma_{s-2,1}\subseteq \{\Gamma_{1,p}\mid p+1\in T\}$.
Since $s$ and $t$ are even, from Lemma \ref{oegy5}, one has  $\Gamma_{2,s-1}\Gamma_{s-2,1}\subseteq \{\Gamma_{1,s},\Gamma_{1,s-2},\ldots,\Gamma_{1,t+2}\}$.
It follows from Lemma \ref{jb} \ref{jb-0} and \ref{jb-3} that
	\begin{align}
	k_{2,s-1}= p_{(2,s-1),(s,1)}^{(1,s-2)}+\ldots+p_{(2,s-1),(t+2,1)}^{(1,s-2)}\leq k_{1,s}+\dots+k_{1,t+2},\nonumber
	\end{align}
	contrary to \eqref{5gs-4}.
\end{proof}

\subsection{$s$ is odd} In this subsection, we always assume that $s$ is odd.

In view of Lemma \ref{even} \ref{even-1} and Lemmas \ref{maxt}, \ref{rem}, one obtains
\begin{align}
\Gamma_{1,q}\Gamma_{1,s}&\subseteq \{\Gamma_{1,s-1},\Gamma_{1,s},\Gamma_{2,s-1},\Gamma_{2,s+1}\},\label{34gso-1}\tag{5.5}\\
\Gamma_{1,s-1}^2&\subseteq \{\Gamma_{1,s-2},\Gamma_{1,s-1},\Gamma_{1,s},\Gamma_{2,s-2},\Gamma_{2,s-1},\Gamma_{2,s+1}\},\label{34gso-2}\tag{5.6}
\end{align}
where $q+1\in T$.


\begin{lemma}\label{odd57}
	Let $T\subseteq \{2,4,\ldots,s-1,s+1\}$ and $\Gamma_{1,q}\Gamma_{1,s}\subseteq \{\Gamma_{2,s-1},\Gamma_{2,s+1}\}$ for any $q+1\in T$. Suppose $\Gamma_{2,s+1}\in \Gamma_{1,q}\Gamma_{1,s}$ for some $q\geq 5$. Then the following hold:
	\begin{enumerate}
		\item\label{odd57-1} 	$\Gamma_{2,s+1}\in \Gamma_{1,s}\Gamma_{1,1}\cup \Gamma_{1,s}\Gamma_{3,1}$;
		
		\item\label{odd57-2}
		$\Gamma_{2,s-1}\in \Gamma_{1,s}\Gamma_{1,1}\cup \Gamma_{1,s}\Gamma_{1,3}$;
		
			\item\label{odd57-5}
			If 	$q=s$, then $\Gamma_{2,s-1}\in  \Gamma_{1,s-2}\Gamma_{1,s}\cap \Gamma_{2,s-1}^2$.
	\end{enumerate}	
\end{lemma}
\begin{proof}
Let $(x,y)\in \Gamma_{2,s+1}$ and $z\in P_{(1,q),(1,s)}(x,y)$.
Pick a shortest path $(y=y_0,y_1,\ldots,y_{s},y_{s+1}=x)$.
If $z\in \{y_0,y_1,\ldots,y_{s+1}\}$,  then $z=y_s$ since $\partial(y,z)=s$, which implies $q=\partial(z,x)=\partial(y_s,x)=1$, a contradiction. Thus,    $z\notin \{y_0,y_1,\ldots,y_{s+1}\}$. Since $\Gamma$ is $s$-quasi-transitive and  $(y_2,\ldots,y_s,x,z)$ is a path of length $s$, we obtain either $(z,y_2)\in A\Gamma$ or $(y_2,z)\in A\Gamma$.
Since $T\subseteq \{2,4,\ldots,s-1,s+1\}$ and $5\leq s=\partial(y,z)\leq 2+\partial(y_2,z)$, we obtain $(z,y_2)\in \Gamma_{1,s-2}\cup \Gamma_{1,s}$. By Lemma \ref{2s1} \ref{2s1-3}, one has  $(y_{s-1},z)\in A\Gamma$. Since $\partial(y_2,z)\leq \partial(y_2,y_{s-1})+\partial(y_{s-1},z)=s-2$, one has $(z,y_2)\in \Gamma_{1,s-2}$.

By Lemma \ref{2s1} \ref{2s1-2}, we obtain $(x,y_1),(y_s,y)\in \Gamma_{1,s}$. Observe that $\Gamma_{1,q}\Gamma_{1,s}\subseteq \{\Gamma_{2,s-1},\Gamma_{2,s+1}\}$ for $q+1\in T$. Since $y_1\in P_{(1,s),(1,q_1-1)}(x,y_2)$, $z\in P_{(1,q_2-1),(1,s)}(y_{s-1},y)$ and $y\in P_{(1,s),(1,q_3-1)}(y_s,y_1)\cap P_{(1,s),(1,q_3-1)}(z,y_1)$ with $q_1,q_2,q_3\in T$, one obtains $(x,y_2),(y_{s-1},y),(y_s,y_1),(z,y_1)\in \Gamma_{2,s-1}\cup \Gamma_{2,s+1}$. Since   $(y,y_1,\ldots,y_{s},x)$ is a shortest path, we have  $(x,y_2),(y_{s-1},y),(y_s,y_1)\in \Gamma_{2,s-1}$. The fact $\partial(y_1,z)\leq 1+\partial(y_2,z)=s-1$ implies  $(z,y_1)\in \Gamma_{2,s-1}$.

By Lemma \ref{2s1} \ref{2s1-1}, one has $(y_i,y_{i+1})\in \Gamma_{1,1}\cup \Gamma_{1,3}$ for any $0\leq i\leq s$.

(i)  Since $y_1\in P_{(1,s),(1,1)}(x,y)\cup P_{(1,s),(3,1)}(x,y)$, \ref{odd57-1} holds.

(ii) Since $y_1\in P_{(1,s),(1,1)}(x,y_2)\cup P_{(1,s),(1,3)}(x,y_2)$,  \ref{odd57-2} holds.



(iii)  Since $(y_{s-1},z)\in A\Gamma$ and $s=q=\partial(z,x)\leq \partial(z,y_{s-1})+2$, one has $(y_{s-1},z)\in \Gamma_{1,s-2}\cup \Gamma_{1,s}$. The fact $\partial(z,y_{s-1})\leq \partial(z,y_2)+\partial(y_2,y_{s-1})=s-2$ implies $(y_{s-1},z)\in \Gamma_{1,s-2}$. Since $x\in P_{(1,1),(1,s)}(y_{s},z)\cup P_{(1,3),(1,s)}(y_{s},z)$ and $\Gamma_{1,q}\Gamma_{1,s}\subseteq \{\Gamma_{2,s-1},\Gamma_{2,s+1}\}$ for any $q+1\in T$, we get $(y_s,z)\in \Gamma_{2,s-1}\cup \Gamma_{2,s+1}$. Since $\partial(z,y_s)\leq 1+\partial(y_2,y_s)=s-1$, one has $(y_s,z)\in \Gamma_{2,s-1}$. The fact  $z\in  P_{(1,s-2),(1,s)}(y_{s-1},y)\cap P_{(2,s-1),(2,s-1)}(y_{s},y_1)$ implies that the desired result follows.
\end{proof}

Next, we divide this subsection into two cases according to the value of $s$.

\subsubsection{$s=5$}

The main result of this subsubsection is as follows.

\begin{prop}\label{txodd5-2}
	Let $s=5$. Then $T=\{6\}$.
\end{prop}

To prove Proposition \ref{txodd5-2}, we need some auxiliary lemmas. In the remainder of this subsubsection, we always assume that $s=5$.

\begin{lemma}\label{txodd5-1}
We have $5\notin T$.
\end{lemma}
\begin{proof}
	Assume the contrary, namely, $5\in T$.
	Suppose $\Gamma_{1,p}\in \Gamma_{1,4}\Gamma_{5,1}$ for some  $p+1\in T$. Lemma \ref{jb} \ref{jb-2} implies $\Gamma_{1,4}\in \Gamma_{1,p}\Gamma_{1,5}$. By Lemma \ref{1205}, one gets  $\Gamma_{2,4},\Gamma_{2,6}\notin \Gamma_{1,5}\Gamma_{1,4}$. It follows from \eqref{34gso-1} that $\Gamma_{1,4}\Gamma_{1,5}\subseteq \{\Gamma_{1,4},\Gamma_{1,5}\}$.
	In view of Lemma \ref{jb}  \ref{jb-2} and \ref{jb-3}, one obtains $k_{1,5}=p_{(5,1),(1,4)}^{(1,4)}+p_{(5,1),(1,5)}^{(1,4)}=p_{(1,5),(1,4)}^{(1,4)}+p_{(1,5),(5,1)}^{(1,4)}$, which  implies  $\Gamma_{1,4}\Gamma_{5,1}\subseteq \{\Gamma_{1,4},\Gamma_{5,1}\}$, and so $p_{(1,4),(5,1)}^{(4,1)}=p_{(1,4),(5,1)}^{(1,3)}=p_{(1,4),(5,1)}^{(1,5)}=0$. By Lemma \ref{jb} \ref{jb-0} and \ref{jb-2}, we get $\Gamma_{1,5}\notin \Gamma_{1,4}^2$ and $\Gamma_{1,4}\notin \Gamma_{1,3}\Gamma_{1,5}\cup \Gamma_{1,5}^2$.
	Since $\Gamma_{1,p}\in \Gamma_{1,4}\Gamma_{5,1}$ and $\Gamma_{1,4}\Gamma_{5,1}\subseteq \{\Gamma_{1,4},\Gamma_{5,1}\}$, one has $p=4$.  Lemma \ref{jb} \ref{jb-2} implies that $\Gamma_{1,4}\in \Gamma_{1,4}\Gamma_{1,5}$. 
    By Lemma \ref{jichu3}, we get $\Gamma_{2,m}\notin \Gamma_{1,4}^2$ for any $m\geq 2$. Since $\Gamma_{1,5}\notin \Gamma_{1,4}^2$, from \eqref{34gso-2}, one has 
	$\Gamma_{1,4}^2\subseteq \{\Gamma_{1,3},\Gamma_{1,4}\}$. Lemma \ref{one} implies  $\Gamma_{1,3}\in \Gamma_{1,4}^2$. By Lemma \ref{1205}, we have $\Gamma_{2,4},\Gamma_{2,6}\notin \Gamma_{1,5}\Gamma_{1,3}$. \eqref{34gso-1} implies $\Gamma_{1,3}\Gamma_{1,5}=\{\Gamma_{1,5}\}$.  By Lemma \ref{jichu3} and \eqref{34gso-1}, we get
	$\Gamma_{1,5}^2=\{\Gamma_{1,5}\}$, contrary to Lemma \ref{one}.  Thus, $\Gamma_{1,p}\notin \Gamma_{1,4}\Gamma_{5,1}$ for any  $p+1\in T$.
	
	By setting $r=1$ and $t=s=5$ in Lemma \ref{jichu0}, there exists $p_0 \geq 2$ such that $\Gamma_{p_0,1}\in \Gamma_{1,4}\Gamma_{5,1}$.  By Lemma \ref{jb} \ref{jb-0} and \ref{jb-2}, we have $\Gamma_{1,5}\in \Gamma_{1,p_0}\Gamma_{1,4}$. By \eqref{34gso-1} and Lemma \ref{1205}, one gets $\Gamma_{1,5}^2\subseteq \{\Gamma_{1,4},\Gamma_{1,5}\}$.  Lemma \ref{jb} \ref{jb-2} and Lemma \ref{one} imply that  $p_{(1,4),(5,1)}^{(1,5)}=\frac{k_{1,4}}{k_{1,5}}p_{(1,5),(1,5)}^{(1,4)}\ne 0$, a contradiction.
\end{proof}

By \eqref{34gso-1} and Lemma \ref{txodd5-1}, we get $\Gamma_{1,q}\Gamma_{1,5}\subseteq \{\Gamma_{1,5},\Gamma_{2,4},\Gamma_{2,6}\}$ for $q+1\in T$. In view of Lemma \ref{one}, one has $\Gamma_{2,4}$ or $\Gamma_{2,6}\in \Gamma_{1,5}^2$. Lemma \ref{1205} implies $\Gamma_{1,5}\notin\Gamma_{1,q}\Gamma_{1,5}$, and so
\begin{align}
\Gamma_{1,q}\Gamma_{1,5}&\subseteq \{\Gamma_{2,4},\Gamma_{2,6}\},\label{odd5-3}\tag{5.7}
\end{align}
where $q+1\in T$. If $3\in T$, by setting $r=0$ in Lemma \ref{jichu6}, then $5\in T$ from  \eqref{odd5-3}, a contradiction. Since $\max T=6$, we obtain $T\subseteq \{2,4,6\}$.

\begin{lemma}\label{txodd5-2'}
	If $\Gamma_{2,4}\in \Gamma_{1,3}\Gamma_{1,5}\cap \Gamma_{2,4}^2$, then $\Gamma_{2,4}\in \Gamma_{1,3}\Gamma_{3,1}$.
\end{lemma}
\begin{proof}
	Let $(x,y)\in \Gamma_{2,4}$ and $z\in P_{(2,4),(2,4)}(x,y)$. Since $\Gamma_{2,4}\in \Gamma_{1,3}\Gamma_{1,5}$, there exist vertices $z_1\in P_{(1,3),(1,5)}(x,z)$ and $z_2\in P_{(1,3),(1,5)}(z,y)$. Pick a shortest path $(y,y_1,y_2,y_3,x)$.  Since $\partial(y,z)=4$ and $\partial(y,z_2)=5$, one has $z,z_2\notin \{y_1,y_2,y_3\}$. If $z_1\in \{y_1,y_2,y_3\}$, then $z_1=y_1$ since $\partial(z_1,x)=3$, which implies $4=\partial(y,z)\leq \partial(y,y_1)+\partial(z_1,z)=2$, a contradiction. Thus, $z_1\notin \{y_1,y_2,y_3\}$.

Since $5=\partial(z,z_1)\leq \partial(z,y_3)+\partial(y_3,x)+1=\partial(z,y_3)+2$ and $5=\partial(y,z_2)\leq 1+\partial(y_1,z)+1$, we have $\partial(z,y_3),\partial(y_1,z)\geq3$. Note that $(z,z_2,y,y_1,y_2,y_3)$ and $(y_1,y_2,y_3,x,z_1,z)$ are two paths of length $5$. Since $T\subseteq \{2,4,6\}$ and $\Gamma$ is $5$-quasi-transitive, one has $(y_3,z),(z,y_1)\in \Gamma_{1,3}\cup \Gamma_{1,5}$. Since $(y_3,z,y_1,y_2)$ is a circuit, one obtains  $(z,y_1)\in\Gamma_{1,3}$. Since $\partial(y_1,z_2)\leq 2+\partial(y_3,z)+1=4$ and $y\in P_{(1,5),(1,q)}(z_2,y_1)$ for some $q+1\in T$, from \eqref{odd5-3}, one has $(z_2,y_1)\in \Gamma_{2,4}$. Since $z\in P_{(3,1),(1,3)}(z_2,y_1)$,  the desired result follows.
\end{proof}

\begin{lemma}\label{txodd5-2''}
	If $\Gamma_{2,6}\in \Gamma_{1,5}^2$, then $\Gamma_{2,6}\Gamma_{3,1}\subseteq \{\Gamma_{1,3},\Gamma_{1,5}\}$.
\end{lemma}
\begin{proof}
	Let $(x,y)\in \Gamma_{2,6}$ and $(y,w)\in \Gamma_{3,1}$. Since $\Gamma_{2,6}\in \Gamma_{1,5}^2$, there exists $z\in P_{(1,5),(1,5)}(x,y)$. Pick a shortest path $(y,y_1,y_2,w)$.
	Since $\partial(y,z)=5$ and $\partial(y,x)=6$, we get $x,z\notin \{y_1,y_2\}$.
	 Since $\Gamma$ is $5$-quasi-transitive and $(x,z,y,y_1,y_2,w)$ is a path of length $5$, one has either $(x,w)\in A\Gamma$ or $(w,x)\in A\Gamma$. Since $T\subseteq \{2,4,6\}$ and $6=\partial(y,x)\leq 3+\partial(w,x)$, we get  $(x,w)\in \Gamma_{1,3}\cup \Gamma_{1,5}$. Thus, the desired result follows.
\end{proof}

\begin{lemma}\label{txodd5-3''}
 If $\Gamma_{2,6}\in \Gamma_{1,3}\Gamma_{1,5}\cap \Gamma_{1,5}^2$, then $\Gamma_{2,4}\in \Gamma_{3,1}\Gamma_{1,5}$.
\end{lemma}
\begin{proof}
	Let $(x,y)\in \Gamma_{2,6}$, $w_1\in P_{(1,5),(1,5)}(x,y)$ and $w_2\in P_{(1,3),(1,5)}(x,y)$. Then $w_1\ne w_2$. Pick a shortest path $(y,y_1,y_2,y_3,y_4,w_2)$.	Since $y\in P_{(1,5),(1,q)}(w_2,y_1)$ for some $q+1\in T$, from \eqref{odd5-3}, one has $(w_2,y_1)\in \Gamma_{2,4}$.
	Since $\partial(y,x)=6$ and $\partial(y,w_1)=5$, one gets $x,w_1\notin \{y,y_1,y_2,y_3,y_4,w_2\}$.
	Since $\Gamma$ is $5$-quasi-transitive and $(x,w_1,y,y_1,y_2,y_3)$ is a path of length $5$, we have either $(x,y_3)\in A\Gamma$ or $(y_3,x)\in A\Gamma$. The fact $6=\partial(y,x)\leq 3+\partial(y_3,x)$ implies $(x,y_3)\in A\Gamma$. Since $6=\partial(y,x)\leq 1+\partial(y_1,x)$ and $(x,y_3,y_4,w_2,y,y_1)$ is a path of length $5$, one obtains $(x,y_1)\in \Gamma_{1,5}$, and so $x\in P_{(3,1),(1,5)}(w_2,y_1)$. Then the desired result follows.
\end{proof}	

\begin{lemma}\label{lem411fl}
If $\Gamma_{2,4}\notin \Gamma_{1,5}\Gamma_{1,1}$ and $\Gamma_{1,5}^2=\{\Gamma_{2,4}\}$, then $\Gamma_{2,6}\notin \Gamma_{1,3}\Gamma_{1,5}$.
\end{lemma}
\begin{proof}
	Suppose for the contrary that $\Gamma_{2,6}\in \Gamma_{1,3}\Gamma_{1,5}$.	Let $(x,y)\in \Gamma_{2,6}$ and $z\in P_{(1,5),(1,3)}(x,y)$. Pick a shortest path $(y=y_0,y_1,y_2,y_3,y_4,y_{5},y_6=x)$.
	Since $T\subseteq \{2,4,6\}$, by Lemma \ref{2s1} \ref{2s1-2} and \ref{2s1-1}, one has $(x,y_{1}),(y_{5},y)\in \Gamma_{1,5}$ and $(y,y_1)=(y_0,y_1)\in \Gamma_{1,1}\cup \Gamma_{1,3}$. Since $\partial(y_1,y_5)=4$ and $y\in P_{(1,5),(1,3)}(y_5,y_1)\cup P_{(1,5),(1,1)}(y_5,y_1)$, from \eqref{odd5-3},  we get $(y_5,y_1)\in \Gamma_{2,4}$.
	Since $\Gamma_{2,4}\notin \Gamma_{1,5}\Gamma_{1,1}$, one has
	$(y,y_1)\in \Gamma_{1,3}$, which implies that $y_1\in P_{(1,5),(3,1)}(x,y),$ and so  $\Gamma_{2,6}\in \Gamma_{1,5}\Gamma_{3,1}$.
	
	 Pick a vertex $w\in \Gamma_{1,5}(y)$.  Since $\Gamma_{2,6}\in \Gamma_{1,5}\Gamma_{3,1}$, there exists $u\in P_{(3,1),(1,5)}(x,y)$. Since $\Gamma_{1,5}^2=\{\Gamma_{2,4}\}$, we get $(u,w)\in \Gamma_{2,4}$.
	The fact $6=\partial(y,x)\leq 1+\partial(w,x)$ implies that $\partial(w,x)\geq 5$.
	Since $u\in P_{(3,1),(2,4)}(x,w)$, by setting $r=2$ and $t=s-2=3$ in Lemma \ref{jichu0}, one has $(x,w)\in \Gamma_{1,5}$. Since $(x,y)\in \Gamma_{2,6}$ and $w\in \Gamma_{1,5}(y)$ were arbitrary, we obtain $\Gamma_{2,6}\Gamma_{1,5}=\{\Gamma_{1,5}\}$.
	
	By Lemma \ref{jb} \ref{jb-1} and \ref{jb-2}, we get $k_{2,6}=p_{(2,6),(1,5)}^{(1,5)}=p_{(6,2),(1,5)}^{(1,5)}$, which implies   $\Gamma_{1,5}\Gamma_{6,2}=\{\Gamma_{1,5}\}$. In view of Lemma \ref{jb} \ref{jb-0}, we have $\Gamma_{2,6}\Gamma_{5,1}=\{\Gamma_{5,1}\}$.
	Since $\Gamma_{2,6}\in \Gamma_{1,3}\Gamma_{1,5}$, from Lemma \ref{jb} \ref{jb-2}, one has $\Gamma_{1,3}\in \Gamma_{2,6}\Gamma_{5,1}$,
	a contradiction.
\end{proof}

	\begin{lemma}\label{lem411f2}
	 If $\Gamma_{1,5}\Gamma_{1,1}=\{\Gamma_{2,6}\}$, then $\Gamma_{2,6}\in \Gamma_{1,3}^2$.
	\end{lemma}
	\begin{proof}
		Let $(x,y)\in \Gamma_{2,6}$, $z_1\in P_{(1,5),(1,1)}(x,y)$ and $z_2\in P_{(1,1),(1,5)}(x,y)$. Pick a shortest path $(y,y_1,y_2,y_3,y_4,y_{5},x)$.
		Suppose $z_1,z_2\in \{y_1,y_2,y_3,y_4,y_{5}\}$. Since $\partial(y,z_2)=\partial(z_1,x)=5$, one has $z_2=y_5$ and $z_1=y_1$, which imply  $x\in P_{(1,1),(1,5)}(y_5,y_1)$. Since  $\Gamma_{1,5}\Gamma_{1,1}=\{\Gamma_{2,6}\}$, we get $(y_5,y_1)\in \Gamma_{2,6}$, contrary to the fact  $\partial(y_1,y_5)=4$. Thus, $z_i\notin \{y_1,y_2,y_3,y_4,y_{5}\}$  for some $i\in \{1,2\}$. Note that  $(y,y_1,y_2,y_3,y_4,y_{5},x)$ is a shortest path. Since $(x,z_i,y,y_1,y_2,y_3)$ and $(y_3,y_4,y_5,x,z_i,y)$ are two paths of length $5$, one has
		 $(x,y_{3}),(y_{3},y)\in \Gamma_{1,3}$, which implies $y_3\in P_{(1,3),(1,3)}(x,y)$, and so $\Gamma_{2,6}\in \Gamma_{1,3}^2$.
	\end{proof}

	\begin{lemma}\label{lem411f3}
 If $\Gamma_{2,4}\in  \Gamma_{2,4}^2$ and $\Gamma_{1,3}\notin \Gamma_{2,4}\Gamma_{1,5}$, then $\Gamma_{2,6}\in \Gamma_{1,3}\Gamma_{1,5}$.
\end{lemma}
\begin{proof}
	Suppose not. By \eqref{odd5-3}, we get $\Gamma_{1,3}\Gamma_{1,5}=\{\Gamma_{2,4}\}$.   Let $(x,y)\in \Gamma_{2,4}$ and $(y,w)\in \Gamma_{1,5}$. By Lemma \ref{txodd5-2'}, we have $\Gamma_{2,4}\in \Gamma_{1,3}\Gamma_{3,1}$. Then there exists $z\in P_{(3,1),(1,3)}(x,y)$. Since $\Gamma_{1,3}\Gamma_{1,5}=\{\Gamma_{2,4}\}$ and $y\in P_{(1,3),(1,5)}(z,w)$, one gets $(z,w)\in \Gamma_{2,4}$. Since $4=\partial(y,x)\leq 1+\partial(w,x)$, one obtains $\partial(w,x)\geq 3$. Since $y\in P_{(2,4),(1,5)}(x,w)$ and $\Gamma_{1,3}\notin \Gamma_{2,4}\Gamma_{1,5}$,  we get $(x,w)\notin \Gamma_{1,3}$.
Since $z\in P_{(3,1),(2,4)}(x,w)$ and $T\subseteq \{2,4,6\}$, by setting $r=2$ and $t=s-2=3$ in Lemma \ref{jichu0}, we have $(x,w)\in \Gamma_{1,5}$. Since $(x,y)\in \Gamma_{2,4}$ and $(y,w)\in \Gamma_{1,5}$ were arbitrary, one gets $\Gamma_{2,4}\Gamma_{1,5}=\{\Gamma_{1,5}\}$.
	
In view of Lemma \ref{jb} \ref{jb-1} and \ref{jb-2}, one obtains
	$k_{2,4}=p_{(2,4),(1,5)}^{(1,5)}=p_{(4,2),(1,5)}^{(1,5)}$,
	which implies that $\Gamma_{1,5}\Gamma_{4,2}=\{\Gamma_{1,5}\}$. Since $\Gamma_{2,4}\in \Gamma_{1,3}\Gamma_{1,5}$, from  Lemma \ref{jb} \ref{jb-0} and \ref{jb-2}, one has 
	$\Gamma_{3,1}\in \Gamma_{1,5}\Gamma_{4,2}$, a contradiction.
		\end{proof}
	
Now we are ready to give a proof of Proposition \ref{txodd5-2}.

\begin{proof}[Proof of Proposition \ref{txodd5-2}]
	Suppose not. Since $T\subseteq \{2,4,6\}$, we get  $2$ or $4\in T$. If $2\in T$, from \eqref{odd5-3}, then $\Gamma_{1,1}\Gamma_{1,5}\subseteq \{\Gamma_{2,4},\Gamma_{2,6}\}$, which implies $4\in T$ by setting $r=0$ in Lemma \ref{jichu6}. Thus, $T=\{4,6\}$ or $\{2,4,6\}$.

By \eqref{odd5-3}, we have $\Gamma_{2,b}\in\Gamma_{1,5}^2$ for some $b\in\{4,6\}$.  Lemma \ref{oegy4} implies $\Gamma_{2,4}\notin \Gamma_{2,b}\Gamma_{2,4}\cup \Gamma_{2,b}\Gamma_{2,6}$. In view of \eqref{odd5-3}, one gets $\sum_{\tilde{r}\in\wz{\partial}(\Gamma)}p_{(1,a),(1,5)}^{\tilde{r}}p_{(2,b),\tilde{r}}^{(2,4)}=0$  for $a\in \{1,3,5\}$. By setting $i=(1,a)$, $l=(1,5)$, $j=(2,4)$ and $h=(2,b)$ in Lemma \ref{jb} \ref{jb-4},   one obtains
	\begin{align}
	p_{(1,a),(2,b)}^{(1,3)}p_{(1,3),(1,5)}^{(2,4)}=p_{(1,a),(2,b)}^{(1,5)}p_{(1,5),(1,5)}^{(2,4)}=0.\label{odd5-4}\tag{5.8}
	\end{align}
We divide the proof into two cases.

	\textbf{Case 1.} $\Gamma_{2,6}\notin \Gamma_{1,5}^2$.
	
By \eqref{odd5-3}, we have  $\Gamma_{1,5}^2=\{\Gamma_{2,4}\}$, and so $b=4$. By \eqref{odd5-4}, we obtain $p_{(1,a),(2,4)}^{(1,5)}=0$ for $a\in\{1,3,5\}$.  Lemma \ref{jb}  \ref{jb-2} implies $\Gamma_{2,4}\notin \Gamma_{1,5}\Gamma_{a,1}$ and $\Gamma_{a,1}\notin \Gamma_{2,4}\Gamma_{5,1}$. 

Since  $\Gamma_{2,4}\notin \Gamma_{1,5}\Gamma_{1,1}$ and $\Gamma_{1,5}^2=\{\Gamma_{2,4}\}$, from Lemma \ref{lem411fl}, one has $\Gamma_{2,6}\notin\Gamma_{1,3}\Gamma_{1,5}$, which implies $\Gamma_{1,3}\Gamma_{1,5}=\{\Gamma_{2,4}\}$ from \eqref{odd5-3}.  By Lemma \ref{jb} \ref{jb-2} and setting $r=2$ and $t=s=5$ in Lemma \ref{jichu0}, one has $\Gamma_{2,4}\Gamma_{5,1}=\{\Gamma_{1,3},\Gamma_{1,5}\}$.  By \eqref{odd5-4}, we get $p_{(1,a),(2,4)}^{(1,3)}=0$ for $a\in\{1,3,5\}$.    By Lemma \ref{jb} \ref{jb-0} and \ref{jb-2}, one obtains $\Gamma_{1,1},\Gamma_{3,1},\Gamma_{5,1}\notin \Gamma_{2,4}\Gamma_{3,1}$.
By setting $r=2$ and $t=s-2=3$ in Lemma \ref{jichu0}, one has  $\Gamma_{2,4}\Gamma_{3,1}\subseteq \{\Gamma_{1,3},\Gamma_{1,5}\}$.	
	Since $\Gamma_{1,q}\Gamma_{1,5}=\{\Gamma_{2,4}\}$, from Lemma \ref{comm}, we get $p_{(2,4),(q,1)}^{(1,5)}=k_{1,q}$ for $q\in \{3,5\}$. Since $\Gamma_{2,4}\Gamma_{3,1}\subseteq\{\Gamma_{1,3},\Gamma_{1,5}\}=\Gamma_{2,4}\Gamma_{5,1}$, from Lemma \ref{jb} \ref{jb-2} and \ref{jb-3}, one obtains
	\begin{align}
	k_{2,4}&= p_{(2,4),(3,1)}^{(1,5)}+p_{(2,4),(5,1)}^{(1,5)}=k_{1,3}+k_{1,5},\nonumber\\
	k_{2,4}&\leq p_{(2,4),(3,1)}^{(1,3)}+p_{(2,4),(5,1)}^{(1,3)}\leq k_{1,3}+k_{1,5},\nonumber
	\end{align}
	which imply that $p_{(2,4),(q,1)}^{(1,3)}=k_{1,q}$ for $q\in \{3,5\}$. By Lemma \ref{comm}, we get $\Gamma_{1,3}^2=\{\Gamma_{2,4}\}$. If $2\notin T$, then $\Gamma_{2,2}\in \Gamma_{1,3}^2$, a contradiction; if $2\in T$, from \eqref{odd5-3}, then $\Gamma_{1,5}\Gamma_{1,1}=\{\Gamma_{2,6}\}$ since $\Gamma_{2,4}\notin\Gamma_{1,5}\Gamma_{1,1}$, contrary to Lemma \ref{lem411f2}.
	
	\textbf{Case 2.} $\Gamma_{2,6}\in \Gamma_{1,5}^2$.

By \eqref{odd5-3} and Lemma \ref{odd57} \ref{odd57-1}, we have $\Gamma_{2,6}\in \Gamma_{1,5}\Gamma_{a',1}$ for some $a'\in \{1,3\}$. By Lemma \ref{jb} \ref{jb-2}, one gets $p_{(2,6),(1,a')}^{(1,5)}=\frac{k_{2,6}}{k_{1,5}}p_{(1,5),(a',1)}^{(2,6)}\ne 0$. By setting $b=6$ and $a=a'$ in \eqref{odd5-4}, we obtain
 $p_{(1,5),(1,5)}^{(2,4)}=0$, and so $\Gamma_{1,5}^2=\{\Gamma_{2,6}\}$.

By \eqref{odd5-3} and  Lemma \ref{odd57} \ref{odd57-5}, one has   $\Gamma_{2,4}\in \Gamma_{1,3}\Gamma_{1,5}\cap \Gamma_{2,4}^2$.
 Since $\Gamma_{2,6}\notin \Gamma_{1,3}\Gamma_{3,1}$, from Lemma \ref{jb} \ref{jb-2}, we get
 $\sum_{\tilde{r}\in\wz{\partial}(\Gamma)}p_{(1,5),(1,5)}^{\tilde{r}}p_{(1,3),\tilde{r}}^{(1,3)}=0$.
 By setting $i=l=(1,5)$  and  $j=h=(1,3)$  in Lemma \ref{jb} \ref{jb-4},   one has
	$p_{(1,5),(1,3)}^{(2,4)}p_{(2,4),(1,5)}^{(1,3)}=0$,  which implies that $p_{(2,4),(1,5)}^{(1,3)}=0$, and so $\Gamma_{1,3}\notin \Gamma_{2,4}\Gamma_{1,5}$.

	Let $(x,y)\in \Gamma_{2,4}$ and $(y,z)\in \Gamma_{1,5}$. Since $\Gamma_{2,4}\in  \Gamma_{2,4}^2$, from Lemma \ref{lem411f3}, we get $\Gamma_{2,6}\in \Gamma_{1,3}\Gamma_{1,5}$. It follows from Lemma \ref{txodd5-3''} that there exists $x_1\in P_{(3,1),(1,5)}(x,y)$. Since $\Gamma_{1,5}^2=\{\Gamma_{2,6}\}$, one has $(x_1,z)\in \Gamma_{2,6}$. Since $x_1\in P_{(3,1),(2,6)}(x,z)$,
	from Lemma \ref{txodd5-2''},
	we get $(x,z)\in \Gamma_{1,3}\cup \Gamma_{1,5}$.
	Since  $\Gamma_{1,3}\notin \Gamma_{2,4}\Gamma_{1,5}$ and $y\in P_{(2,4),(1,5)}(x,z)$, one has 	$(x,z)\notin \Gamma_{1,3}$, and so
	$(x,z)\in \Gamma_{1,5}$. Since $(x,y)\in \Gamma_{2,4}$ and $(y,z)\in \Gamma_{1,5}$ were arbitrary, we get $\Gamma_{2,4}\Gamma_{1,5}=\{\Gamma_{1,5}\}$. 	

By Lemma \ref{jb} \ref{jb-1} and \ref{jb-2}, we obtain
	$k_{2,4}=p_{(2,4),(1,5)}^{(1,5)}=p_{(4,2),(1,5)}^{(1,5)}$, which implies  $\Gamma_{1,5}\Gamma_{4,2}=\{\Gamma_{1,5}\}$. Since  $\Gamma_{2,4}\in \Gamma_{1,3}\Gamma_{1,5}$, from Lemma \ref{jb} \ref{jb-0} and \ref{jb-2}, one has 
	$\Gamma_{3,1}\in \Gamma_{1,5}\Gamma_{4,2}$, a contradiction.
\end{proof}

\subsubsection{$s$ is odd and $s\geq 7$}

The main result of this subsubsection is as follows.

\begin{prop}\label{oprop4}
	Let $s$ be odd with $s\geq 7$. Then $T=\{s+1\}$.
\end{prop}

To prove Proposition \ref{oprop4}, we need some auxiliary lemmas. In this subsubsection, we  assume that  $s$ is odd and $s\geq 7$.

\begin{lemma}\label{34lemma0}
	 Suppose $\Gamma_{2,s+1}\in \Gamma_{1,s}\Gamma_{1,s-1}$. Then
	$k_{1,n}\leq1$ for any $n\in \{1,2,3\}$.
\end{lemma}
\begin{proof}
	Assume the contrary, namely, $k_{1,n}>1$ for some $n\in \{1,2,3\}$. Let $(x,y)\in \Gamma_{1,s-1}$. Pick a shortest path $(y,y_1,\ldots,y_{s-2},x)$.   Since $k_{1,n}>1$, by setting $l=s-1$ and $p=n$ in Lemma \ref{1009}, there exists $z\in \Gamma_{n,1}(y)$ such that $z\notin \{y_1,\ldots,y_{s-2}\}$, which implies that $(z,y,y_1,\ldots,y_{s-2},x)$ is a path of length $s$, and so $(x,z)$ or $(z,x)\in A\Gamma$. Since $s-1=\partial(y,x)\leq n+\partial(z,x)$, we get $\partial(z,x)\geq s-n-1\geq s-4>1$, and so $(x,z)\in \Gamma_{1,q}$ for some $q\geq s-n-1$. Note that $z\in P_{(1,q),(1,n)}(x,y)$. By   the commutativity of $\Gamma$,  one has $\Gamma_{1,s-1}\in \Gamma_{1,n}\Gamma_{1,q}$ for some $q+1\in T$ and $q\geq s-n-1$. Since $s\geq 7$ and $n\leq 3$, we obtain $q\geq 3$.

    Since $s\geq 7$ and $\Gamma_{2,s+1}\in \Gamma_{1,s}\Gamma_{1,s-1}$, from Lemma \ref{1205} \ref{1205-2}, one gets $\Gamma_{1,s-1}\notin \Gamma_{1,p}\Gamma_{1,q}$ for any $p\geq1$ and $q\geq 3$, a contradiction.
\end{proof}

\begin{lemma}\label{34lemma1chai}
	If $\Gamma_{2,s+1}\in \Gamma_{1,s}\Gamma_{1,t}$ for some $t\geq 5$, then $\Gamma_{1,s}\Gamma_{1,s-1}\neq\{\Gamma_{1,s-1}\}$.
\end{lemma}
\begin{proof}
Assume the contrary, namely, $\Gamma_{1,s}\Gamma_{1,s-1}=\{\Gamma_{1,s-1}\}$. By Lemma \ref{jichu3}, one has $\Gamma_{2,m}\notin \Gamma_{1,s-1}^2$ for any $m\geq 2$.

Suppose $\Gamma_{1,s}\in \Gamma_{1,s-1}^2$. Since $\Gamma_{2,s+1}\in \Gamma_{1,s}\Gamma_{1,t}$, 
there exist  $(x,y)\in \Gamma_{2,s+1}$ and $z\in  P_{(1,t),(1,s)}(x,y)$.  Pick a shortest path $(y=y_0,y_1,\ldots,y_{s},y_{s+1}=x)$. Since $\Gamma_{1,s}\in \Gamma_{1,s-1}^2$, there exists $z_1\in  P_{(1,s-1),(1,s-1)}(z,y)$. Since $t\geq 5$, by setting $l=0$ in Lemma \ref{1207}, we get
$z_1\in \{y_1,\ldots,y_{s}\}$. Since $\partial(y,z_1)=s-1$, one has $z_1=y_{s-1}$, which implies  $t=\partial(z,x)\leq \partial(z,z_1)+\partial(y_{s-1},x)=3$, a contradiction. Since $s\geq 7$, from
\eqref{34gso-2}, we obtain $\Gamma_{1,s-1}^2\subseteq \{\Gamma_{1,s-2},\Gamma_{1,s-1}\}$. In view of Lemma \ref{one}, one has $\Gamma_{1,s-2}\in \Gamma_{1,s-1}^2$.

 Suppose that $\Gamma_{2,m}\in \Gamma_{1,s-2}\Gamma_{1,s}$ for some $m\geq 2$. Let $(u,v)\in \Gamma_{2,m}$ and  $u_0\in  P_{(1,s-2),(1,s)}(u,v)$. Since $\Gamma_{1,s-2}\in \Gamma_{1,s-1}^2$, there exists $w\in P_{(1,s-1),(1,s-1)}(u,u_0)$. Since $\Gamma_{1,s}\Gamma_{1,s-1}=\{\Gamma_{1,s-1}\}$, we get $(w,v)\in \Gamma_{1,s-1}$, and so $w\in P_{(1,s-1),(1,s-1)}(u,v)$, contrary to the fact that $\Gamma_{1,s-1}^2\subseteq \{\Gamma_{1,s-2},\Gamma_{1,s-1}\}$. Then $\Gamma_{2,m}\notin \Gamma_{1,s-2}\Gamma_{1,s}$ for any $m\geq 2$.  Since $s\geq 7$, from \eqref{34gso-1}, one obtains $\Gamma_{1,s-2}\Gamma_{1,s}\subseteq \{\Gamma_{1,s-1},\Gamma_{1,s}\}$.
Since $\Gamma_{2,s+1}\in \Gamma_{1,s}\Gamma_{1,t}$, from Lemma \ref{jichu3}, we have $\Gamma_{1,s}\notin \Gamma_{1,s}\Gamma_{1,s-2}$. Thus,
$\Gamma_{1,s-2}\Gamma_{1,s}= \{\Gamma_{1,s-1}\}$.

By Lemma \ref{jb} \ref{jb-2}, we get $p_{(1,s-1),(s,1)}^{(1,s-2)}=\frac{k_{1,s-1}}{k_{1,s-2}}p_{(1,s-2),(1,s)}^{(1,s-1)}\ne 0$, and so $\Gamma_{1,s-2}\in \Gamma_{1,s-1}\Gamma_{s,1}$. Since $\Gamma_{1,s-1}\Gamma_{1,s}=\{\Gamma_{1,s-1}\}$, from Lemma \ref{jb} \ref{jb-1} and  \ref{jb-2}, one has
$k_{1,s}=p_{(1,s),(1,s-1)}^{(1,s-1)}=p_{(s,1),(1,s-1)}^{(1,s-1)}$, which implies  $\Gamma_{1,s-1}\Gamma_{s,1}=\{\Gamma_{1,s-1}\}$,
contrary to the fact that $\Gamma_{1,s-2}\in \Gamma_{1,s-1}\Gamma_{s,1}$.
\end{proof}

\begin{lemma}\label{34lemma1}
	Let  $s\in T$. Then $\Gamma_{2,s+1}\notin \Gamma_{1,s}\Gamma_{1,t}$ for any $t\geq 5$.
\end{lemma}
\begin{proof}
	Assume the contrary, namely, $\Gamma_{2,s+1}\in \Gamma_{1,s}\Gamma_{1,t}$ for some $t\geq 5$. 	Let $(x,y)\in \Gamma_{2,s+1}$.  Pick a shortest path $(y=y_0,y_1,\ldots,y_{s},y_{s+1}=x)$.
	By Lemmas \ref{jichu3} and \ref{jichu4}, one has $\Gamma_{1,s},\Gamma_{2,s-1}\notin \Gamma_{1,s}\Gamma_{1,s-1}$.  \eqref{34gso-1} implies that $\Gamma_{1,s}\Gamma_{1,s-1}\subseteq \{\Gamma_{1,s-1},\Gamma_{2,s+1}\}$.

	By Lemma \ref{34lemma1chai}, we get  $\Gamma_{2,s+1}\in \Gamma_{1,s}\Gamma_{1,s-1}$. Pick a vertex  $z\in P_{(1,s),(1,s-1)}(x,y)$.  By Lemma \ref{2s1} \ref{2s1-1}, one has $(y_i,y_{i+1})\in \Gamma_{1,1}\cup\Gamma_{1,2}\cup\Gamma_{1,3}$   and $(y_{l+3},y_{l})\in \Gamma_{1,3}$   for   $0\leq i\leq s$ and  $0\leq l\leq s-2$. It follows from
	 Lemma \ref{34lemma0} that
	$k_{1,3}=1$ and  $k_{1,1},k_{1,2}\leq 1$.
	
	If there exists $i_0\in \{0,1,\ldots,s\}$ such that $(y_{i_0},y_{i_0+1})\in \Gamma_{1,3}$, then $0\leq i_0+4\leq s+1$ or $0\leq i_0-3\leq s+1$  since $i_0\leq s-3$ or $i_0\geq s-2\geq5$, which implies $y_{i_0+4},y_{i_0}\in \Gamma_{3,1}(y_{i_0+1})$ or $y_{i_0-3},y_{i_0+1}\in \Gamma_{1,3}(y_{i_0})$, contrary to the fact that $k_{1,3}=1$. Thus,
	$(y_i,y_{i+1})\in \Gamma_{1,1}\cup\Gamma_{1,2}$ for any $0\leq i\leq s$.
	
	Suppose $(y_i,y_{i+1})\in \Gamma_{1,2}$ for any $0\leq i\leq s$. Since $y_1\in P_{(1,2),(1,2)}(y_0,y_2)$, from  Lemma \ref{2s1} \ref{2s1-1}, one has either $\Gamma_{2,1}\in \Gamma_{1,2}^2$ or $\Gamma_{2,2}\in \Gamma_{1,2}^2$. Since  $k_{1,2}=1$, from Lemma \ref{jb} \ref{jb-5}, we get $\Gamma_{1,2}^2=\{\Gamma_{2,1}\}$ or $\{\Gamma_{2,2}\}$. If $\Gamma_{1,2}^2=\{\Gamma_{2,1}\}$, then $(y_0,y_2)\in \Gamma_{2,1}$, and so $y_0,y_3\in \Gamma_{1,2}(y_2)$,  contrary to the fact that $k_{1,2}=1$. Then $\Gamma_{1,2}^2=\{\Gamma_{2,2}\}$.  Since $k_{1,2}=1$, from Lemma \ref{jb} \ref{jb-1}, one has $k_{2,2}=1$. Since $s\geq 7$, we have $(y_0,y_2),(y_2,y_4)\in \Gamma_{2,2}$, which implies that
	$y_4,y_0\in \Gamma_{2,2}(y_2)$, contrary to the fact that $k_{2,2}=1$.
	
By the commutativity of $\Gamma$, we may assume that  $(y_{1},y_{2})\in \Gamma_{1,1}$. Since  $k_{1,1}=1$, we get $(y_{0},y_{1}),(y_{2},y_{3})\in \Gamma_{1,2}$. By Lemma \ref{jb} \ref{jb-5} and Lemma \ref{2s1} \ref{2s1-1}, one has $\Gamma_{1,2}\Gamma_{1,1}=\{\Gamma_{2,1}\}$ or $\{\Gamma_{2,2}\}$. If $\Gamma_{1,2}\Gamma_{1,1}=\{\Gamma_{2,1}\}$, then $(y_{0},y_{2})\in \Gamma_{2,1}$, and so $y_{0},y_{3}\in \Gamma_{1,2}(y_{2})$, contrary to the fact that $k_{1,2}=1$. Then $\Gamma_{1,2}\Gamma_{1,1}=\{\Gamma_{2,2}\}$, and so $(y_{0},y_{2})\in \Gamma_{2,2}$. Since $k_{1,2}=1$, from Lemma \ref{comm0}, one has  $\Gamma_{2,2}\Gamma_{1,2}=\{\Gamma_{1,1}\}$, which implies  that $(y_{0},y_{3})\in \Gamma_{1,1}$, contrary to Lemma \ref{2s1} \ref{2s1-1}.
\end{proof}

\begin{lemma}\label{oddarc}
The following hold:
	\begin{enumerate}
		\item\label{oddarc-2} $\Gamma_{1,q}\Gamma_{1,s}\subseteq \{\Gamma_{2,s-1},\Gamma_{2,s+1}\}$ for any $q+1\in T$;

		\item\label{oddarc-1}
		$T\subseteq \{2,4,\ldots,s-1,s+1\}$.
	\end{enumerate}	
\end{lemma}
\begin{proof}
By Lemma \ref{34lemma1}, one has $\Gamma_{2,s+1}\notin \Gamma_{1,s}\Gamma_{1,t}$ for any $t\geq 5$. Since $s\geq 7$, from Lemma \ref{jiougy1}, we get $s\notin T$.

\ref{oddarc-2} By \eqref{34gso-1}, one gets $\Gamma_{1,q}\Gamma_{1,s}\subseteq \{\Gamma_{1,s},\Gamma_{2,s-1},\Gamma_{2,s+1}\}$ for any $q+1\in T$. Assume the contrary, namely, $\Gamma_{1,s}\in \Gamma_{1,q}\Gamma_{1,s}$ for some  $q+1\in T$. Since $\Gamma_{1,s}^2\subseteq \{\Gamma_{1,s},\Gamma_{2,s-1},\Gamma_{2,s+1}\}$,  from Lemma \ref{one}, we have $\Gamma_{2,s-1}$ or $\Gamma_{2,s+1}\in \Gamma_{1,s}^2$. Let $(x,y)\in \Gamma_{2,a}$ with $a\in \{s-1,s+1\}$ and $z\in  P_{(1,s),(1,s)}(x,y)$. Pick a shortest path $(y,y_1,\ldots,y_{a-1},x)$. Since $\partial(z,x)=\partial(y,z)=s$, one gets $z\notin \{y,y_1,\ldots,y_{a-1},x\}$. Since $\Gamma_{1,s}\in \Gamma_{1,q}\Gamma_{1,s}$, there exists $z_0\in P_{(1,q),(1,s)}(z,y)$.

	Suppose $a=s+1$. Since $s\geq 7$, by setting $l=0$ in Lemma \ref{1207}, we obtain  $z_0\in \{y_1,\ldots,y_s\}$. Since $\partial(y,z_0)=s$ and  $(y,y_1,\ldots,y_{s},x)$ is a shortest path, one has $z_0=y_s$, which implies $s=\partial(z,x)\leq \partial(z,z_0)+\partial(y_s,x)=2$, a contradiction. Then $a=s-1$.

Since $\partial(y,z_0)=s$, one has $z_0\notin \{z,y_1,\ldots,y_{s-2}\}$. Since $(z,z_0,y,y_1,\ldots,y_{s-2})$ is a path of length $s$, we get $(y_{s-2},z)$ or $(z,y_{s-2})\in A\Gamma$. Since $s=\partial(z,x)\leq \partial(z,y_{s-2})+1$, one obtains $(y_{s-2},z)\in A\Gamma$, which implies $s=\partial(y,z)\leq \partial(y,y_{s-2})+\partial(y_{s-2},z)=s-1$, a contradiction. Thus,  $\Gamma_{1,s}\notin \Gamma_{1,q}\Gamma_{1,s}$ for any $q+1\in T$, and so \ref{oddarc-2} holds.

\ref{oddarc-1} If $p+1\in T$ for some positive even integer $p$, by \ref{oddarc-2} and setting $r=0$ in Lemma \ref{jichu6}, then $s\in T$ since $s$ is odd, a contradiction. Thus,  \ref{oddarc-1} holds.
\end{proof}

\begin{lemma}\label{oprop4chai-1}
If $q\geq5$, then $\Gamma_{2,s+1}\notin \Gamma_{1,q}\Gamma_{1,s}$.
\end{lemma}
\begin{proof}
	Assume the contrary, namely, $\Gamma_{2,s+1}\in \Gamma_{1,q}\Gamma_{1,s}$. Lemma \ref{odd57} \ref{odd57-1}  and Lemma \ref{oddarc} imply that $\Gamma_{2,s+1}\in \Gamma_{1,s}\Gamma_{a,1}$ for some $a\in \{1,3\}$. By Lemma \ref{jb} \ref{jb-2}, we get $p_{(2,s+1),(1,a)}^{(1,s)}=\frac{k_{2,s+1}}{k_{1,s}}p_{(1,s),(a,1)}^{(2,s+1)}\ne 0$.

By Lemma \ref{oegy4}, we have $\Gamma_{2,s-1}\notin \Gamma_{2,s+1}^2\cup\Gamma_{2,s+1}\Gamma_{2,s-1}$. In view of Lemma \ref{oddarc} \ref{oddarc-2}, one has $\sum_{\tilde{r}\in\wz{\partial}(\Gamma)}p_{(1,a),(1,s)}^{\tilde{r}}p_{(2,s+1),\tilde{r}}^{(2,s-1)}=0$.  By setting $i=(1,a)$, $l=(1,s)$,  $j=(2,s-1)$ and $h=(2,s+1)$ in Lemma \ref{jb} \ref{jb-4}, one gets
	$p_{(2,s+1),(1,a)}^{(1,s)}p_{(1,s),(1,s)}^{(2,s-1)}=0$, and so $p_{(1,s),(1,s)}^{(2,s-1)}=0$. It follows from Lemma \ref{oddarc} \ref{oddarc-2} that $\Gamma_{1,s}^2=\{\Gamma_{2,s+1}\}$. By Lemma \ref{odd57} \ref{odd57-5} and Lemma \ref{oddarc}, one gets  $\Gamma_{2,s-1}\in \Gamma_{1,s-2}\Gamma_{1,s}\cap \Gamma_{2,s-1}^2$.  Since $s\geq 7$  and $\Gamma_{2,s-1}\in \Gamma_{1,s-2}\Gamma_{1,s}$, from Lemma \ref{oegy4}, we have $\Gamma_{2,s-1}\notin \Gamma_{2,s-1}^2$, a contradiction.
\end{proof}

\begin{lemma}\label{oprop4chai-2}
	If $\Gamma_{1,s}^2=\{\Gamma_{2,s-1}\}$ and $\Gamma_{2,s-1}\notin \Gamma_{2,s-1}^2$, then $T\subsetneq\{2,4,\ldots,s-1,s+1\}$.
\end{lemma}
\begin{proof}
	By Lemma \ref{jb} \ref{jb-2}, we get $\sum_{\tilde{r}\in\wz{\partial}(\Gamma)}p_{(1,s),(1,s)}^{\tilde{r}}p_{(s-1,2),\tilde{r}}^{(2,s-1)}=0$. In view of Lemma \ref{oddarc} \ref{oddarc-1}, one obtains $T\subseteq\{2,4,\ldots,s-1,s+1\}$. Assume the contrary, namely, $T=\{2,4,\ldots,s-1,s+1\}$.	

 By setting $i=l=(1,s)$ and $j=h^*=(2,s-1)$ in Lemma \ref{jb} \ref{jb-4}, one has
	$\sum_{ \wz{t}\in\wz{\partial}(\Gamma)}p_{(s-1,2),(1,s)}^{\wz{t}}p_{\wz{t},(1,s)}^{(2,s-1)}=0$, and so $p_{(s-1,2),(1,s)}^{(1,1)}p_{(1,1),(1,s)}^{(2,s-1)}=0$. By Lemma \ref{jb}  \ref{jb-2}, we get $\frac{k_{2,s-1}}{k_{1,1}}p_{(1,1),(s,1)}^{(s-1,2)}p_{(1,1),(1,s)}^{(2,s-1)}=0$.
	It follows from Lemma \ref{jb} \ref{jb-0} that $\frac{k_{2,s-1}}{k_{1,1}}p_{(1,1),(1,s)}^{(2,s-1)}p_{(1,1),(1,s)}^{(2,s-1)}=0$.
	Thus, $p_{(1,1),(1,s)}^{(2,s-1)}=0$, and so $\Gamma_{2,s-1}\notin \Gamma_{1,1}\Gamma_{1,s}$. By Lemma \ref{oddarc} \ref{oddarc-2}, one has  $\Gamma_{1,1}\Gamma_{1,s}=\{\Gamma_{2,s+1}\}$. 	
	
Since $\Gamma_{a,b}\in \Gamma_{s,1}\Gamma_{1,s-2}$ implies $b\leq1+(s-2)$, we have $\Gamma_{2,s+1}\notin \Gamma_{s,1}\Gamma_{1,s-2}$.  Lemma \ref{jb} \ref{jb-2} implies 	$\sum_{ \wz{t}\in\wz{\partial}(\Gamma)}p_{(s,1),(1,s-2)}^{\wz{t}}p_{\wz{t},(1,1)}^{(1,s)}=0$. 	By setting $i=(1,s-2)$, $l=(1,1)$ and $j=h^*=(1,s)$ in Lemma \ref{jb} \ref{jb-4}, one obtains $p_{(1,s-2),(1,1)}^{(2,s-1)}p_{(s,1),(2,s-1)}^{(1,s)}=0$. Since $\Gamma_{1,s}^2=\{\Gamma_{2,s-1}\}$, from Lemma \ref{jb} \ref{jb-2},  we get $p_{(s,1),(2,s-1)}^{(1,s)}\ne 0$, and so $\Gamma_{2,s-1}\notin\Gamma_{1,1}\Gamma_{1,s-2}$.

Since $\Gamma_{2,s-1}\notin\Gamma_{1,1}\Gamma_{1,s-2}\cup\Gamma_{1,1}\Gamma_{1,s}$, from Lemma \ref{jb} \ref{jb-2}, one has $\Gamma_{1,s-2},\Gamma_{1,s}\notin \Gamma_{2,s-1}\Gamma_{1,1}$. 
Let $(x,y)\in \Gamma_{2,s-1}$ and $(y,y_1,\ldots,y_{s-2},x)$ be a shortest path.  Let $v\in \Gamma_{1,1}(y)$.   Suppose $v\notin \{y_1,\ldots,y_{s-2}\}$. Since $(v,y,y_1,\ldots,y_{s-2},x)$ is a path of length $s$, one has $(x,v)$ or $(v,x)\in A\Gamma$. Since $T=\{2,4,\ldots,s-1,s+1\}$ and $s-1=\partial(y,x)\leq 1+\partial(v,x)$, we get $(x,v)\in \Gamma_{1,s}\cup \Gamma_{1,s-2}$. Since $y\in P_{(2,s-1),(1,1)}(x,v)$, one has $\Gamma_{1,s}$ or $\Gamma_{1,s-2}\in \Gamma_{2,s-1}\Gamma_{1,1}$, a contradiction.  Thus,  $v\in \{y_1,\ldots,y_{s-2}\}$. Since $\partial(y,v)=1$, we have $v=y_1$.  Thus,  $(y,y_1)\in \Gamma_{1,1}$.
	
	Since $\Gamma_{1,s}^2=\{\Gamma_{2,s-1}\}$, there exists $w\in P_{(1,s),(1,s)}(x,y)$.  Since $\Gamma_{1,1}\Gamma_{1,s}=\{\Gamma_{2,s+1}\}$ and $y\in P_{(1,s),(1,1)}(w,y_1)$, one obtains $(w,y_1)\in \Gamma_{2,s+1}$, which implies that $s+1=\partial(y_1,w)\leq \partial(y_1,x)+\partial(x,w)=s-1$, a contradiction.
\end{proof}

Now we are ready to give a proof of Proposition \ref{oprop4}.

\begin{proof}[Proof of Proposition \ref{oprop4}]
Assume the contrary, namely, $t+1=\min T<s+1$.
	By Lemma \ref{oddarc} \ref{oddarc-1},  $t$ is  odd.
	By Lemma \ref{oddarc} \ref{oddarc-2} and setting $r=0$ in Lemma \ref{jichu6}, we have $t+3\in T$.
	Repeating this process, since $\max T=s+1$,
	we get $T=\{t+1,t+3,\ldots,s-1,s+1\}$.

	 By Lemma \ref{oddarc} \ref{oddarc-2} and  Lemma \ref{oprop4chai-1}, we get $\Gamma_{1,s}\Gamma_{1,q}=\{\Gamma_{2,s-1}\}$  for $q\geq 5$, which implies $p_{(2,s-1),(q,1)}^{(1,s)}=k_{1,q}$ from Lemma \ref{comm}.	Since $\Gamma_{1,s}^2=\{\Gamma_{2,s-1}\}$, from Lemma \ref{oegy4}, we get $\Gamma_{2,s-1}\notin \Gamma_{2,s-1}^2$, which implies that $\sum_{\tilde{t}\in\wz{\partial}(\Gamma)}p_{(1,s),(1,s)}^{\tilde{t}}p_{\tilde{t},(2,s-1)}^{(2,s-1)}=0$.
	  By setting $i=h=(1,s)$ and $l=j=(2,s-1)$ in Lemma \ref{jb} \ref{jb-4}, one obtains
	$p_{(1,s),(2,s-1)}^{(1,q)}p_{(1,s),(1,q)}^{(2,s-1)}=0$ for any  $q\geq 5$,
	which implies $p_{(1,s),(2,s-1)}^{(1,q)}=0.$
	
If $\Gamma_{l,1}\in\Gamma_{2,s-1}\Gamma_{s-2,1}$ for some $l>1$, from  Lemma \ref{oegy6}, then $\Gamma_{2,s-1}\Gamma_{1,s}\subseteq \{\Gamma_{1,s-2},\Gamma_{1,s}\}$ since $s\geq7$, contrary to the fact that $p_{(1,s),(2,s-1)}^{(1,q)}=0$ for $q\geq5$. Then $\Gamma_{l,1}\notin \Gamma_{2,s-1}\Gamma_{s-2,1}$ for any $l>1$. By setting $r=2$ and $t=s-2$ in Lemma \ref{jichu0}, we get  $\Gamma_{2,s-1}\Gamma_{s-2,1}\subseteq \{\Gamma_{1,t},\Gamma_{1,t+2},\ldots,\Gamma_{1,s}\}$.
	
Suppose $t\geq 5$. By Lemma \ref{oegy5}, one has $\Gamma_{2,s-1}\Gamma_{s-2,1}\subseteq \{\Gamma_{1,t+2},\ldots,\Gamma_{1,s-2},\Gamma_{1,s}\}$. 	By Lemma \ref{jb} \ref{jb-0}, \ref{jb-2} and  \ref{jb-3}, we obtain 	$k_{2,s-1}=p_{(2,s-1),(s,1)}^{(1,s-2)}+\ldots+p_{(2,s-1),(t+2,1)}^{(1,s-2)}\leq k_{1,s}+\dots+k_{1,t+2}$. Since $p_{(2,s-1),(q,1)}^{(1,s)}=k_{1,q}$ for $q\geq 5$, from Lemma \ref{jb} \ref{jb-3}, one has
	$	k_{2,s-1}\geq  p_{(2,s-1),(s,1)}^{(1,s)}+\ldots+p_{(2,s-1),(t,1)}^{(1,s)}=k_{1,s}+\dots+k_{1,t}$,	a contradiction. Since $\Gamma_{1,s}^2=\{\Gamma_{2,s-1}\}$ and $\Gamma_{2,s-1}\notin \Gamma_{2,s-1}^2$, from Lemma \ref{oprop4chai-2},
	we get  $t=3$.

	Suppose $\Gamma_{2,s+1}\in \Gamma_{1,3}\Gamma_{1,s}$. Let $(x,y)\in \Gamma_{2,s+1}$ and $z\in P_{(1,s),(1,3)}(x,y)$. Pick a shortest path $(y,y_1,\ldots,y_{s},x)$. Since $s$ is odd and $s\geq 7$, from Lemma \ref{2s1} \ref{2s1-2} and \ref{2s1-3}, we get $(x,y_3)\in \Gamma_{1,s-2}$ and $(x,y_1)\in \Gamma_{1,s}$. Since $\min T=4$, from Lemma \ref{2s1} \ref{2s1-1}, we obtain $(y_1,y_2),(y_2,y_3)\in \Gamma_{1,3}$. Since $y_1\in P_{(1,s),(1,3)}(x,y_2)$, from Lemma \ref{oddarc} \ref{oddarc-2}, one has $(x,y_2)\in \Gamma_{2,s-1}\cup \Gamma_{2,s+1}$. Since $(y,y_1,\ldots,y_{s},x)$ is  a shortest path, one obtains $(x,y_2)\in \Gamma_{2,s-1}$. Since $y_3\in P_{(1,s-2),(3,1)}(x,y_2)$, we get $\Gamma_{2,s-1}\in \Gamma_{1,s-2}\Gamma_{3,1}$. By Lemma \ref{jb} \ref{jb-2}, one has $\Gamma_{3,1}\in \Gamma_{2,s-1}\Gamma_{s-2,1}$, a contradiction. In view of Lemma \ref{oddarc} \ref{oddarc-2}, one has $\Gamma_{1,3}\Gamma_{1,s}=\{\Gamma_{2,s-1}\}$. By  Lemma \ref{comm}, we get $p_{(2,s-1),(3,1)}^{(1,s)}=k_{1,3}$.
	
Suppose  $\Gamma_{l,1}\in \Gamma_{2,s-1}\Gamma_{s,1}$ for some $l+1\in T$. Since $\min T=4$, one gets $l\geq 3$. In view of Lemma \ref{jb} \ref{jb-2}, we have  $\Gamma_{2,s-1}\in \Gamma_{l,1}\Gamma_{1,s}$. Let $(u,v)\in \Gamma_{2,s-1}$ and $(v,v_0)\in \Gamma_{1,s}$. Since  $\Gamma_{2,s-1}\in \Gamma_{l,1}\Gamma_{1,s}$, there exists a vertex  $w\in P_{(l,1),(1,s)}(u,v)$. Since $\Gamma_{1,s}^2=\{\Gamma_{2,s-1}\}$, one gets $(w,v_0)\in \Gamma_{2,s-1}$. The fact $s-1=\partial(v,u)\leq 1+\partial(v_0,u)$ implies that $\partial(v_0,u)\geq s-2$. Since $w\in P_{(l,1),(2,s-1)}(u,v_0)$ and $l\geq 3$, by setting $r=2$ and $t=l$ in Lemma \ref{jichu0}, we have $\partial(u,v_0)=1$, and so
	 $(u,v_0)\in \Gamma_{1,s-2}\cup \Gamma_{1,s}$. Since $(u,v)\in \Gamma_{2,s-1}$ and $(v,v_0)\in \Gamma_{1,s}$ were arbitrary, one gets $\Gamma_{2,s-1}\Gamma_{1,s}\subseteq \{\Gamma_{1,s-2},\Gamma_{1,s}\}$, contrary to the fact that $p_{(1,s),(2,s-1)}^{(1,q)}=0$ for $q\geq5$. Thus,  $\Gamma_{l,1}\notin \Gamma_{2,s-1}\Gamma_{s,1}$ for any $l+1\in T$.

By setting $r=2$ and $t\in \{s-2,s\}$ in Lemma \ref{jichu0}, we obtain $\Gamma_{2,s-1}\Gamma_{s-2,1},\Gamma_{2,s-1}\Gamma_{s,1}\subseteq \{\Gamma_{1,3},\Gamma_{1,5},\ldots,\Gamma_{1,s}\}$.	
Note that $p_{(2,s-1),(q,1)}^{(1,s)}=k_{1,q}$ for $q\geq3$. By Lemma \ref{jb} \ref{jb-0}, \ref{jb-2} and \ref{jb-3}, we get
	\begin{align}
	k_{2,s-1}&=p_{(2,s-1),(s,1)}^{(1,s)}+p_{(2,s-1),(s-2,1)}^{(1,s)}+\ldots+p_{(2,s-1),(3,1)}^{(1,s)}=k_{1,s}+k_{1,s-2}+\dots+k_{1,3},\nonumber\\
	k_{2,s-1}&=p_{(2,s-1),(s,1)}^{(1,s-2)}+p_{(2,s-1),(s-2,1)}^{(1,s-2)}+\ldots+p_{(2,s-1),(3,1)}^{(1,s-2)}\leq k_{1,s}+k_{1,s-2}+\dots+k_{1,3},\nonumber
	\end{align}
	which imply that $p_{(2,s-1),(q,1)}^{(1,s-2)}=k_{1,q}$ for any $q+1\in T$. By  Lemma \ref{comm}, we get $\Gamma_{1,s-2}\Gamma_{1,q}=\{\Gamma_{2,s-1}\}$ for any $q+1\in T$, contrary to the fact that $\Gamma_{l,s-3}\in\Gamma_{1,s-2}\Gamma_{1,h}$ for some $h+1\in T$ and $l\in\{1,2\}$.
\end{proof}

Now we are ready to give a proof of Proposition \ref{propxj}.

\begin{proof}[Proof of Proposition \ref{propxj}]
	In view of Propositions \ref{prop4}, \ref{txodd5-2} and \ref{oprop4}, we get $T=\{s+1\}$.
	By Lemma \ref{prop1} \ref{prop1-2}, Lemma \ref{oddarc} \ref{oddarc-2} and \eqref{odd5-3}, one has $\Gamma_{1,s}^2\subseteq \{\Gamma_{2,s-1},\Gamma_{2,s+1}\}$.
	If $\Gamma_{2,s+1}\in \Gamma_{1,s}^2$, from Lemma \ref{2s1} \ref{2s1-1},  then $4\in T$, contrary to the fact that $T=\{s+1\}$. Then  $\Gamma_{1,s}^2=\{\Gamma_{2,s-1}\}$. By Lemma \ref{oegy4}, we have $\Gamma_{2,s-1}\notin \Gamma_{2,s-1}^2$.
	By Lemma \ref{jb} \ref{jb-2}, one gets  $\sum_{\tilde{r}\in\wz{\partial}(\Gamma)}p_{(1,s),(1,s)}^{\tilde{r}}p_{(s-1,2),\tilde{r}}^{(2,s-1)}=0$.
	By setting $i=l=(1,s)$ and $j=h^*=(2,s-1)$ in Lemma \ref{jb} \ref{jb-4}, one obtains  $p_{(1,s),(s-1,2)}^{(1,s)}p_{(1,s),(1,s)}^{(2,s-1)}=0$, and so $p_{(1,s),(s-1,2)}^{(1,s)}=0$. By Lemma \ref{jb} \ref{jb-0}, we get $\Gamma_{s,1}\notin \Gamma_{2,s-1}\Gamma_{s,1}$.  Since $T=\{s+1\}$, by setting $r=2$ and $t=s$ in 
		Lemma \ref{jichu0}, we have $\Gamma_{2,s-1}\Gamma_{s,1}=\{\Gamma_{1,s}\}$. By Lemma \ref{comm}, one obtains $k_{2,s-1}=p_{(s,1),(2,s-1)}^{(1,s)}$. Lemma \ref{jb} \ref{jb-2} implies  $p_{(1,s),(1,s)}^{(2,s-1)}=k_{1,s}$. By \cite[Theorem 2.3]{YF22}, 	the desired result follows.
\end{proof}

We conclude this section with a corollary.

\begin{cor}\label{cor1}
	Let $\Gamma$ be a  commutative weakly distance-regular digraph of valency greater than $3$ and with   diameter at least $s$, where  $s\geq 4$. Then $\Gamma$ is  $s$-quasi-transitive  if and only if
	$\Gamma$ is isomorphic to a $t$-coclique extension of $C_{s+1}$,  where $t>3$.
\end{cor}
\begin{proof}
	If $\Gamma$ has diameter at least $s$,  then $\Gamma$ has cycle polar diameter $s+1$ by   Lemma \ref{even} \ref{even-1} and Lemma \ref{maxt}. In view of Proposition \ref{propxj}, the desired result follows.
\end{proof}

\section{Proofs of Theorems \ref{main thm} and  \ref{main thm1}}

In this section, we prove Theorems \ref{main thm} and  \ref{main thm1}.

\begin{proof}[Proof of Theorem \ref{main thm}]
	Routinely, all digraphs in Theorem \ref{main thm}  are $s$-quasi-transitive. By \cite[Proposition 2.4]{KSW03} and \cite[Lemma 9.1]{lishuang2}, all digraphs in Theorem \ref{main thm} are commutative weakly distance-regular digraphs.
	
	Now we  prove the necessity.    Since $\Gamma$ has cycle polar diameter at least $s$,   from Lemma \ref{maxt},  $\Gamma$ has cycle polar diameter $s$ or $s+1$. If $\Gamma$ has cycle polar diameter $s$, from Proposition \ref{bc4}, then $\Gamma$ is isomorphic to the digraph in Theorem \ref{main thm} \ref{main thm-4}. If $\Gamma$ has cycle polar diameter $s+1$, from Proposition \ref{propxj}, then $\Gamma$ is isomorphic to the digraph in Theorem \ref{main thm} \ref{main thm-3}. 
\end{proof}

In the remainder of this section, $\Gamma$ always denotes  a $4$-quasi-transitive  commutative weakly distance-regular digraph. By Theorem \ref{main thm} and  Corollary \ref{cor1},  we only need to consider the case where $\Gamma$ has cycle polar diameter at most $3$ and diameter at most $3$.
Then
\begin{align}
\Gamma_{1,2}^2,\Gamma_{1,2}\Gamma_{1,1}&\subseteq \{\Gamma_{1,1},\Gamma_{1,2},\Gamma_{2,1},\Gamma_{2,2},\Gamma_{2,3}\},\label{4ni-1}\tag{6.1}\\
\Gamma_{1,1}^2&\subseteq \{\Gamma_{0,0},\Gamma_{1,1},\Gamma_{1,2},\Gamma_{2,1},\Gamma_{2,2}\}.\label{4ni-2}\tag{6.2}
\end{align}

\begin{lemma}\label{4nilemma}
	The following hold:
	\begin{enumerate}
		\item\label{4nilemma-1}
		$\Gamma_{2,1}\in \Gamma_{1,2}^2\cup \Gamma_{1,2}\Gamma_{1,1}$;
		
		\item\label{4nilemma-2} $\Gamma_{1,2}^2\subseteq \{\Gamma_{1,1},\Gamma_{1,2},\Gamma_{2,1}\}$.
	\end{enumerate}
\end{lemma}
\begin{proof}
	(i) Assume the contrary, namely, $\Gamma_{2,1}\notin \Gamma_{1,2}^2\cup \Gamma_{1,2}\Gamma_{1,1}$. Since $T\subseteq \{2,3\}$, we get $\Gamma_{2,1}\in \Gamma_{1,1}^2$.   Lemma \ref{jb} \ref{jb-0} and \ref{jb-2} imply $p_{(1,1),(1,2)}^{(1,1)}=p_{(1,1),(2,1)}^{(1,1)}=\frac{k_{1,2}}{k_{1,1}}p_{(1,1),(1,1)}^{(2,1)}\ne0$, and so
	$\Gamma_{1,1}\in \Gamma_{1,1}\Gamma_{1,2}$.
	
	Let $(x,y),(y,z)\in \Gamma_{1,2}$.
	Since $\Gamma_{2,1}\notin  \Gamma_{1,2}\Gamma_{1,1}$, from Lemma \ref{jb} \ref{jb-0} and  \ref{jb-2}, one has $p_{(1,2),(1,2)}^{(1,1)}=p_{(2,1),(2,1)}^{(1,1)}=\frac{k_{1,2}}{k_{1,1}}p_{(1,2),(1,1)}^{(2,1)}=0$, and so  $\Gamma_{1,1}\notin \Gamma_{1,2}^2$. Combining this with $\Gamma_{2,1}\notin \Gamma_{1,2}^2$, we obtain  $(x,z)\notin \Gamma_{1,1}\cup\Gamma_{2,1}$. Since $\Gamma_{2,1}\in \Gamma_{1,1}^2$,  there exist  $y_1\in P_{(1,1),(1,1)}(y,x)$ and $z_1\in P_{(1,1),(1,1)}(z,y)$. If $y_1\ne z_1$, then $(z,z_1,y,y_1,x)$ is a path of length $4$, which implies  $(x,z)\in \Gamma_{1,2}$ since $(x,z)\notin \Gamma_{1,1}\cup\Gamma_{2,1}$.  If $y_1=z_1$, then $z_1\in P_{(1,1),(1,1)}(x,z)$, which implies  $(x,z)\in  \Gamma_{1,2}\cup \Gamma_{2,2}$ by \eqref{4ni-2}.  Since $(x,y),(y,z)\in\Gamma_{1,2}$ were arbitrary, we get
	$\Gamma_{1,2}^2\subseteq \{\Gamma_{1,2},\Gamma_{2,2}\}$.
	
	By Lemma \ref{one}, one has $\Gamma_{2,2}\in \Gamma_{1,2}^2$. Without loss of generality, we may assume that $(x,z)\in \Gamma_{2,2}$.  Since  $\Gamma_{1,1}\in \Gamma_{1,1}\Gamma_{1,2}$,  there exists $w_1\in P_{(1,1),(1,2)}(y,z_1)$. Since $(x,y,w_1,z_1,z)$ is a path of length $4$, one obtains either $(x,z)\in A\Gamma$ or $(z,x)\in A\Gamma$, a contradiction.
	
(ii)			Let $(x,y),(y,z)\in \Gamma_{1,2}$. By \ref{4nilemma-1},  there exist  $y_1\in P_{(1,2),(1,t)}(y,x)$ and $z_1\in P_{(1,2),(1,t)}(z,y)$ for some $t\in \{1,2\}$. Then $y_1\ne z_1$. If $y_1=z$, then $(x,z)=(x,y_1)\in \Gamma_{t,1}$. If $y_1\ne z$, then  $(z,x)$ or $(x,z)\in A\Gamma$ since $(z,z_1,y,y_1,x)$ is a path of length $4$. Since $(x,y),(y,z)\in \Gamma_{1,2}$ were arbitrary, \ref{4nilemma-2} holds.
	\end{proof}

		\begin{prop}\label{s41}
		Suppose  $T=\{3\}$. Then   $\Gamma$ is isomorphic to one of the digraphs in Theorem \ref{main thm1} \ref{main thm1-7}.
	\end{prop}
	\begin{proof}
Since $T=\{3\}$, from Lemma \ref{4nilemma} \ref{4nilemma-2}, one obtains  $\Gamma_{1,2}^2\subseteq \{\Gamma_{1,2},\Gamma_{2,1}\}$, which implies that $(2,2),(2,3)\notin \wz{\partial}(\Gamma)$. Since $\Gamma$ has diameter at most $3$, we have $\wz{\partial}(\Gamma)\subseteq\{(0,0),(1,2),(2,1),(3,3)\}$.

We claim that $\Gamma_{3,3}\Gamma_{1,2}=\{\Gamma_{1,2}\}$ if $(3,3)\in\wz{\partial}(\Gamma)$. Let $(x,y)\in \Gamma_{3,3}$ and $(x,x_1,x_2,y)$ be a shortest path. Pick a vertex $z\in \Gamma_{1,2}(y)$. If $z\notin \{x_1,x_2\}$,  then  $(x,z)\in \Gamma_{1,2}$  since  $3=\partial(y,x)\leq 1+\partial(z,x)$ and $(x,x_1,x_2,y,z)$ is a path of length $4$. If $z\in \{x_1,x_2\}$, then $z=x_1$ since $\partial(z,y)=2$, which implies $(x,z)=(x,x_1)\in \Gamma_{1,2}$. Since $(x,y)\in \Gamma_{3,3}$ and $z\in \Gamma_{1,2}(y)$ were arbitrary, our claim is valid.

Let $t=\max\{i\mid(i,i)\in\tilde{\partial}(\Gamma)\}$. Then $t\in\{0,3\}$. If $t=0$, then $\langle \Gamma_{0,0}\rangle=\{\Gamma_{0,0}\}$. If $t=3$, from the claim and Lemma \ref{jb} \ref{jb-0}, \ref{jb-2}, then $p_{(3,3),(3,3)}^{(2,1)}=p_{(3,3),(3,3)}^{(1,2)}=\frac{k_{3,3}}{k_{1,2}}p_{(3,3),(1,2)}^{(3,3)}=0$, and so $\Gamma_{3,3}^2\subseteq\{\Gamma_{0,0},\Gamma_{3,3}\}$ since $\wz{\partial}(\Gamma)\subseteq\{(0,0),(1,2),(2,1),(3,3)\}$, which implies $\langle \Gamma_{3,3}\rangle=\{\Gamma_{0,0},\Gamma_{3,3}\}$. We conclude that $F:=\langle \Gamma_{t,t}\rangle=\{\Gamma_{0,0},\Gamma_{t,t}\}$.

Let $(u,v)\in\Gamma_{1,2}$. Pick $u'\in F(u)$ and $v'\in F(v)$. Since $F=\{\Gamma_{0,0},\Gamma_{t,t}\}$, we have $(u,u')\in\Gamma_{0,0}\cup\Gamma_{t,t}$. By the claim, one gets $(u',v)\in \Gamma_{1,2}$. If $v=v'$, then $(u',v')\in\Gamma_{1,2}$. If $v\neq v'$, then $t=3$ and $(v,v')\in\Gamma_{3,3}$, which imply $(u',v')\in \Gamma_{1,2}$ from the claim. We conclude that $(u',v')\in\Gamma_{1,2}$ for $u'\in F(u)$ and $v'\in F(v)$ with $(u,v)\in\Gamma_{1,2}$.	
		
Denote  $\Sigma:=\Gamma/F$. Since $(u',v')\in\Gamma_{1,2}$ for $u'\in F(u)$ and $v'\in F(v)$ with $(u,v)\in\Gamma_{1,2}$, $\Gamma$ is isomorphic to $\Sigma\circ \overline{{K}}_{k_{t,t}+1}$, which implies
	$\tilde{\partial}_{\Gamma}(x_0,y_0)=\tilde{\partial}_{\Sigma}(F(x_0),F(y_0))$ for vertices $x_0,y_0$ with  $y_0\notin F(x_0)$. Since $\wz{\partial}(\Gamma)=\{(0,0),(1,2),(2,1),(t,t)\}$, one has 	$\tilde{\partial}(\Sigma)=\{(0,0),(1,2),(2,1)\}$.

	Let $\wz{h}\in\wz{\partial}(\Sigma)$ and $(F(u_1),F(v_1))\in\Sigma_{\wz{h}}$. Without loss of generality, we may assume that $(u_1,v_1)\in \Gamma_{\wz{h}}$. Let $\wz{i},\wz{j}\in\wz{\partial}(\Sigma)$. If $\wz{i}$ or $\wz{j}=(0,0)$, then $|\Sigma_{\wz{i}}(F(u_1))\cap\Sigma_{\wz{j}^*}(F(v_1))|=\delta_{\wz{h},\wz{j}}$ or $\delta_{\wz{h},\wz{i}}$.
	If $\wz{i},\wz{j}\ne(0,0)$, then
		$$P_{\wz{i},\wz{j}}(u_1,v_1)=\cup_{w\in P_{\wz{i},\wz{j}}(u_1,v_1)}F(w)
		=\cup_{F(w)\in \Sigma_{\wz{i}}(F(u_1))\cap\Sigma_{\wz{j}^*}(F(v_1))}F(w).$$
		Since  $|F(w)|=k_{t,t}+1$ for all  $w\in V\Gamma$, one gets $|\Sigma_{\wz{i}}(F(u_1))\cap\Sigma_{\wz{j}^*}(F(v_1))|=\frac{p_{\wz{i},\wz{j}}^{\wz{h}}}{k_{t,t}+1}.$ Since $\Gamma$ is a commutative weakly distance-regular digraph, $\Sigma$ is a  commutative weakly distance-regular digraph. Thus, the desired result follows.
			\end{proof}

Note that $\Gamma$ has cycle polar diameter at most $3$. In the remainder of this section, we always assume that  $T=\{2,3\}$.

	\begin{lemma}\label{4nilemma4}
	We have $\Gamma_{3,3}\notin \Gamma_{1,2}\Gamma_{2,1}$.
\end{lemma}
\begin{proof}
	Assume the contrary, namely, $\Gamma_{3,3}\in \Gamma_{1,2}\Gamma_{2,1}$.
	Let $(x,y)\in \Gamma_{3,3}$ and $z\in P_{(1,2),(2,1)}(x,y)$.
	Pick a vertex $w\in \Gamma_{1,1}(y)$.
	By  Lemma \ref{4nilemma} \ref{4nilemma-1}, there exists $u\in P_{(1,2),(1,t)}(z,x)$ for some $t\in \{1,2\}$. Since $3=\partial(y,x)\leq 1+\partial(w,x)$ and $\partial(u,x)=1$, one has $w\ne u$.  The fact $(w,y,z,u,x)$ is a path of length $4$ implies $(w,x)$ or $(x,w)\in A\Gamma$. Then $\partial(y,x)=\partial(y,w)+\partial(w,x)=2$ or $\partial(x,y)=\partial(x,w)+\partial(w,y)=2$, contrary to the fact that $(x,y)\in \Gamma_{3,3}$.
\end{proof}

	\begin{lemma}\label{4nilemma2}
	Suppose $\Gamma_{2,2}\in \Gamma_{1,1}^2$. Then $\Gamma_{1,1}\notin \Gamma_{1,2}^2\cup \Gamma_{1,2}\Gamma_{1,1}$.
\end{lemma}
\begin{proof}
	Let $(x,y)\in \Gamma_{2,2}$ and $z\in P_{(1,1),(1,1)}(x,y)$. Assume the contrary, namely, $\Gamma_{1,1}\in \Gamma_{1,2}^2\cup \Gamma_{1,2}\Gamma_{1,1}$.
	 Then  there exist  $z_1\in P_{(1,2),(1,t)}(x,z)$ and $z_2\in P_{(1,2),(1,t)}(z,y)$ for some $t\in \{1,2\}$. Since  $(x,z_1,z,z_2,y)$ is a path of length $4$, one has either $(x,y)\in A\Gamma$ or $(y,x)\in A\Gamma$, a contradiction.
\end{proof}

	\begin{lemma}\label{4nilemma3}
	Suppose $\{\Gamma_{2,2},\Gamma_{2,3}\}\cap \Gamma_{1,2}\Gamma_{1,1}\ne \emptyset$. Then $\Gamma_{1,2}\Gamma_{1,1}\subseteq \{\Gamma_{2,2},\Gamma_{2,3}\}$, $\Gamma_{1,1}^2\subseteq \{\Gamma_{0,0},\Gamma_{2,2}\}$, $\Gamma_{1,2}^2=\{\Gamma_{2,1}\}$ and $\Gamma_{1,2}\Gamma_{2,1}\subseteq \{\Gamma_{0,0},\Gamma_{2,2},\Gamma_{2,3},\Gamma_{3,2}\}$.
\end{lemma}
\begin{proof}
Let $(x,y)\in \Gamma_{2,m}$ for some $m\in \{2,3\}$. Pick a vertex $z\in P_{(1,2),(1,1)}(x,y)$. By Lemma \ref{4nilemma} \ref{4nilemma-1}, there exists $w\in P_{(1,2),(1,t)}(z,x)$ for some $t\in \{1,2\}. $

Suppose $\Gamma_{1,1}\in \Gamma_{1,a}\Gamma_{1,b}$ for some $a,b\in \{1,2\}$. Then there exists $u\in P_{(1,a),(1,b)}(y,z)$. Since $(y,u,z,w,x)$ is a path of length $4$, one has either $(y,x)\in A\Gamma$ or $(x,y)\in A\Gamma$,  a contradiction. Thus, $\Gamma_{1,1}\notin \Gamma_{1,a}\Gamma_{1,b}$ for any $a,b\in \{1,2\}$.

If $\Gamma_{1,2}\in \Gamma_{1,2}\Gamma_{1,1}$, then there exists $v\in P_{(1,2),(1,1)}(z,w)$, which implies $(y,x)$ or $(x,y)\in A\Gamma$ since $(y,z,v,w,x)$ is a path of length $4$, a contradiction. Then $\Gamma_{1,2}\notin \Gamma_{1,2}\Gamma_{1,1}$. 	Since $\{\Gamma_{2,2},\Gamma_{2,3}\}\cap \Gamma_{1,2}\Gamma_{1,1}\ne \emptyset$, from Lemma \ref{jichu3}, we get $\Gamma_{1,2}\notin \Gamma_{1,2}^2$. Thus, $\Gamma_{1,2}\notin \Gamma_{1,2}\Gamma_{1,s}$ for any $s\in \{1,2\}$.

Since $\Gamma_{1,1},\Gamma_{1,2}\notin\Gamma_{1,2}^2$, from Lemma \ref{4nilemma} \ref{4nilemma-2}, we get $\Gamma_{1,2}^2=\{\Gamma_{2,1}\}$. In view of Lemma \ref{jb} \ref{jb-0} and \ref{jb-2}, one gets $p_{(1,2),(1,1)}^{(2,1)}=p_{(2,1),(1,1)}^{(1,2)}=\frac{k_{1,1}}{k_{1,2}}p_{(1,2),(1,2)}^{(1,1)}=0$. Since $\Gamma_{1,1},\Gamma_{1,2}\notin\Gamma_{1,2}\Gamma_{1,1}$, from \eqref{4ni-1}, we obtain $\Gamma_{1,2}\Gamma_{1,1}\subseteq \{\Gamma_{2,2},\Gamma_{2,3}\}$. In view of Lemma \ref{jb} \ref{jb-0} and \ref{jb-2} again, one has $p_{(1,1),(1,1)}^{(2,1)}=p_{(1,1),(1,1)}^{(1,2)}=\frac{k_{1,1}}{k_{1,2}}p_{(1,1),(1,2)}^{(1,1)}=0$. Since $\Gamma_{1,1}\notin\Gamma_{1,1}^2$, one obtains $\Gamma_{1,1}^2\subseteq \{\Gamma_{0,0},\Gamma_{2,2}\}$.

Since $\Gamma_{1,2}\notin \Gamma_{1,2}^2\cup \Gamma_{1,2}\Gamma_{1,1}$, from Lemma \ref{jb} \ref{jb-0} and \ref{jb-2}, one has $p_{(2,1),(1,2)}^{(2,1)}=p_{(1,2),(2,1)}^{(1,2)}=p_{(1,2),(1,2)}^{(1,2)}=0$ and $p_{(1,2),(2,1)}^{(1,1)}=\frac{k_{1,2}}{k_{1,1}}p_{(1,1),(1,2)}^{(1,2)}=0$, and so $\Gamma_{2,1},\Gamma_{1,2},\Gamma_{1,1}\notin \Gamma_{1,2}\Gamma_{2,1}$. Since $\Gamma$ has diameter at most $3$, from Lemma \ref{4nilemma4}, we get $\Gamma_{1,2}\Gamma_{2,1}\subseteq \{\Gamma_{0,0},\Gamma_{2,2},\Gamma_{2,3},\Gamma_{3,2}\}$.
\end{proof}

\begin{lemma}\label{4nilemma5}
Suppose $\Gamma_{1,2}\Gamma_{1,1}=\{\Gamma_{2,3}\}$. Then $\Gamma_{2,3}\Gamma_{1,1}=\{\Gamma_{1,2}\}$ and $\Gamma_{2,3}\Gamma_{2,1}\subseteq\{\Gamma_{1,1},\Gamma_{2,1}\}$.
\end{lemma}
\begin{proof}
	Let $(x,y)\in \Gamma_{2,3}$ and $z\in P_{(1,1),(1,2)}(x,y)$. Since $\Gamma_{2,1}\in \Gamma_{1,2}^2$ by Lemma \ref{4nilemma3}, there exists  $z_1\in P_{(1,2),(1,2)}(y,z)$.
	
	Let  $u\in \Gamma_{1,1}(y)$.
Since $3=\partial(y,x)\leq 1+\partial(u,x)$ and  $(u,y,z_1,z,x)$ is a path of length $4$,  we get  $(x,u)\in \Gamma_{1,2}$. Since $u\in \Gamma_{1,1}(y)$ was arbitrary, we get $\Gamma_{2,3}\Gamma_{1,1}=\{\Gamma_{1,2}\}$.

Let  $v\in \Gamma_{2,1}(y)$. Since $\Gamma_{2,3}\notin \Gamma_{1,2}^2$, one has $(x,v)\notin \Gamma_{1,2}$. If $v=z$, then $(x,v)=(x,z)\in \Gamma_{1,1}$. If $v\ne z$, then $(v,y,z_1,z,x)$ is a path of length $4$, which implies $(x,v)\in\Gamma_{1,1}\cup\Gamma_{2,1}$. Since $v\in \Gamma_{2,1}(y)$ was arbitrary, we get $\Gamma_{2,3}\Gamma_{2,1}\subseteq\{\Gamma_{1,1},\Gamma_{2,1}\}$.
\end{proof}

	\begin{lemma}\label{s42}
	We have $\Gamma_{2,2}\notin \Gamma_{1,1}^2$.
\end{lemma}
\begin{proof}
	Assume the contrary, namely, $\Gamma_{2,2}\in \Gamma_{1,1}^2$.
By Lemma \ref{4nilemma2}, one has $\Gamma_{1,1}\notin \Gamma_{1,2}^2\cup \Gamma_{1,2}\Gamma_{1,1}$.    Lemma \ref{jb} \ref{jb-0} and \ref{jb-2} imply
$p_{(1,2),(1,1)}^{(2,1)}=p_{(2,1),(1,1)}^{(1,2)}=\frac{k_{1,1}}{k_{1,2}}p_{(1,2),(1,2)}^{(1,1)}=0$.  It follows from \eqref{4ni-1} and Lemma \ref{4nilemma} \ref{4nilemma-2} that $\Gamma_{1,2}\Gamma_{1,1}\subseteq \{\Gamma_{1,2},\Gamma_{2,2},\Gamma_{2,3}\}$ and  $\Gamma_{1,2}^2\subseteq \{\Gamma_{1,2},\Gamma_{2,1}\}$. Lemma \ref{one} implies $\Gamma_{2,1}\in \Gamma_{1,2}^2$.

Let $(x,y)\in \Gamma_{2,2}$ and $z\in P_{(1,1),(1,1)}(x,y)$. Suppose $\Gamma_{2,2},\Gamma_{2,3}\notin \Gamma_{1,2}\Gamma_{1,1}$. Then $\Gamma_{1,2}\Gamma_{1,1}=\{\Gamma_{1,2}\}$. By Lemma \ref{jb} \ref{jb-2}, one has
$p_{(1,2),(2,1)}^{(1,1)}=\frac{k_{1,2}}{k_{1,1}}p_{(1,1),(1,2)}^{(1,2)}\ne0$, and so $\Gamma_{1,1}\in \Gamma_{1,2}\Gamma_{2,1}$. Combining this with $\Gamma_{2,1}\in \Gamma_{1,2}^2$, there exist  $w_1\in P_{(1,2),(2,1)}(x,z)$ and $w_2\in P_{(1,2),(1,2)}(w_1,z)$. Since  $(x,w_1,w_2,z,y)$ is a path of length $4$, we get $(x,y)$ or $(y,x)\in A\Gamma$, a contradiction. Thus, $\{\Gamma_{2,2},\Gamma_{2,3}\}\cap \Gamma_{1,2}\Gamma_{1,1}\ne \emptyset$.

By  Lemma \ref{4nilemma3}, we get $\Gamma_{1,2}^2=\{\Gamma_{2,1}\}$, $\Gamma_{1,2}\Gamma_{1,1}\subseteq \{\Gamma_{2,2},\Gamma_{2,3}\}$ and $
\Gamma_{1,1}^2=\{\Gamma_{0,0},\Gamma_{2,2}\}.$ 	 Lemma \ref{jb} \ref{jb-2} implies that
$\sum_{\tilde{r}\in\wz{\partial}(\Gamma)}p_{(1,1),(1,1)}^{\tilde{r}}p_{(2,1),\tilde{r}}^{(1,2)}=0$ and $\sum_{\tilde{t}\in\wz{\partial}(\Gamma)}p_{(1,2),(1,2)}^{\tilde{t}}p_{\tilde{t},(1,1)}^{(1,2)}=0$. By setting $i=l=(1,1)$ and $j=h^*=(1,2)$ in Lemma \ref{jb} \ref{jb-4},  we get $p_{(1,1),(2,1)}^{(2,2)}p_{(2,2),(1,1)}^{(1,2)}=0$. It follows from Lemma \ref{jb} \ref{jb-0} and \ref{jb-2} that $\Gamma_{2,2}\notin \Gamma_{1,2}\Gamma_{1,1}$, and so $\Gamma_{1,2}\Gamma_{1,1}=\{\Gamma_{2,3}\}$.

By Lemma \ref{4nilemma5}, one has $\Gamma_{2,3}\Gamma_{1,1}=\{\Gamma_{1,2}\}$ and $\Gamma_{2,3}\Gamma_{2,1}\subseteq\{\Gamma_{1,1},\Gamma_{2,1}\}$. Note that $\sum_{\tilde{t}\in\wz{\partial}(\Gamma)}p_{(1,2),(1,2)}^{\tilde{t}}p_{\tilde{t},(1,1)}^{(1,2)}=0$.
By setting $i=j=h=(1,2)$ and $l=(1,1)$ in Lemma \ref{jb} \ref{jb-4}, one has  $p_{(1,2),(1,1)}^{(2,3)}p_{(1,2),(2,3)}^{(1,2)}=0$, and so $p_{(1,2),(2,3)}^{(1,2)}=0$. Lemma \ref{jb} \ref{jb-2} implies $p_{(2,1),(1,2)}^{(2,3)}=p_{(2,3),(2,1)}^{(2,1)}=0$, and so $\Gamma_{2,3}\notin \Gamma_{1,2}\Gamma_{2,1}$ and $\Gamma_{2,1}\notin \Gamma_{2,3}\Gamma_{2,1}$. Then $\Gamma_{2,3}\Gamma_{2,1}=\{\Gamma_{1,1}\}$. Combining this with $\Gamma_{2,3}\Gamma_{1,1}=\{\Gamma_{1,2}\}$, we obtain $k_{1,1}=p_{(1,1),(1,2)}^{(2,3)}=p_{(1,2),(1,1)}^{(2,3)}=k_{1,2}$ by   Lemma \ref{jb} \ref{jb-2} and \ref{jb-3}. Since $\Gamma_{1,2}\Gamma_{1,1}=\{\Gamma_{2,3}\}$, from Lemma \ref{jb} \ref{jb-1}, we get $p_{(1,1),(1,2)}^{(2,3)}=\frac{k_{1,1}k_{1,2}}{k_{2,3}}$. It follows that $k_{1,1}=k_{1,2}=k_{2,3}$.

Since $\Gamma_{2,3}\notin \Gamma_{1,2}\Gamma_{2,1}$, from Lemma \ref{jb} \ref{jb-0} and Lemma  \ref{4nilemma3}, one gets $\Gamma_{1,2}\Gamma_{2,1}\subseteq \{\Gamma_{0,0},\Gamma_{2,2}\}$. Since $
\Gamma_{1,1}^2=\{\Gamma_{0,0},\Gamma_{2,2}\}$, we get $k_{1,2}=k_{1,1}>1$, which implies $\Gamma_{1,2}\Gamma_{2,1}=\{\Gamma_{0,0},\Gamma_{2,2}\}$.

Since $k_{1,1}=k_{1,2}$, from Lemma \ref{jb} \ref{jb-2} and \ref{jb-3}, we have $p_{(1,2),(2,2)}^{(1,2)}=p_{(1,1),(2,2)}^{(1,1)}=k_{1,1}-1$. By Lemma \ref{jb} \ref{jb-2}, one obtains $p_{(2,1),(1,2)}^{(2,2)}=\frac{k_{1,2}}{k_{2,2}}p_{(1,2),(2,2)}^{(1,2)}=\frac{k_{1,1}}{k_{2,2}}p_{(1,2),(2,2)}^{(1,2)}=\frac{k_{1,1}(k_{1,1}-1)}{k_{2,2}}$. By setting $i=l^*=(1,2)$ and $j=h=(1,1)$ in Lemma \ref{jb} \ref{jb-4}, we get $k_{1,1}+p_{(1,2),(2,1)}^{(2,2)}p_{(1,1),(2,2)}^{(1,1)}=p_{(1,1),(1,2)}^{(2,3)}p_{(2,3),(2,1)}^{(1,1)}=(p_{(1,1),(1,2)}^{(2,3)})^2$. By Lemma \ref{jb} \ref{jb-2}, one has $k_{2,2}=k_{1,1}-1$, and so $p_{(1,2),(2,2)}^{(1,2)}=p_{(1,1),(2,2)}^{(1,1)}=k_{2,2}$. It follows from Lemma \ref{comm} that $\Gamma_{1,2}\Gamma_{2,2}=\{\Gamma_{1,2}\}$ and $\Gamma_{1,1}\Gamma_{2,2}=\{\Gamma_{1,1}\}$. Then $p_{(1,2),(2,2)}^{(2,3)}=p_{(1,2),(2,2)}^{(3,2)}=p_{(1,1),(2,2)}^{(2,3)}=0$. By Lemma \ref{jb} \ref{jb-0} and \ref{jb-2}, one has $\Gamma_{1,2},\Gamma_{2,1},\Gamma_{1,1}\notin \Gamma_{2,3}\Gamma_{2,2}$.

Let  $u\in \Gamma_{2,3}(y)$. Since $\Gamma_{2,3}\in \Gamma_{1,2}\Gamma_{1,1}$, there exists $y_1\in P_{(1,2),(1,1)}(y,u)$. Since $(x,z,y,y_1,u)$ is a path of length $4$, we obtain $(x,u)\in A\Gamma$ or $(u,x)\in A\Gamma$, contrary to the fact that $\Gamma_{1,1},\Gamma_{1,2},\Gamma_{2,1}\notin \Gamma_{2,2}\Gamma_{2,3}$.
\end{proof}

	\begin{prop}\label{s43}
	Suppose  $T=\{2,3\}$. Then $\Gamma$ is isomorphic to one of the digraphs in Theorem \ref{main thm1} \ref{main thm1-2}.
\end{prop}
\begin{proof}
Suppose $\{\Gamma_{2,2},\Gamma_{2,3}\}\cap \Gamma_{1,2}\Gamma_{1,1}\ne \emptyset$. By  Lemmas \ref{4nilemma3} and  \ref{s42},  we have $\Gamma_{1,2}^2=\{\Gamma_{2,1}\}$, $\Gamma_{1,2}\Gamma_{1,1}\subseteq \{\Gamma_{2,2},\Gamma_{2,3}\}$, $\Gamma_{1,1}^2=\{\Gamma_{0,0}\}$ and $\Gamma_{1,2}\Gamma_{2,1}\subseteq \{\Gamma_{0,0},\Gamma_{2,2},\Gamma_{2,3},\Gamma_{3,2}\}$.
	
	Since $\Gamma_{1,1}^2=\{\Gamma_{0,0}\}$, one gets $k_{1,1}=1$. In view of Lemma \ref{jb} \ref{jb-5}, we obtain $\Gamma_{1,2}\Gamma_{1,1}=\{\Gamma_{2,m}\}$ for some $m\in \{2,3\}$.
	Since  $\Gamma_{1,2}^2=\{\Gamma_{2,1}\}$, one has $\Gamma_{1,2}\Gamma_{2,1}\subseteq \{\Gamma_{0,0},\Gamma_{2,m},\Gamma_{m,2}\}$.

	Since $\Gamma_{1,2}^2=\{\Gamma_{2,1}\}$ and  $\Gamma_{1,2}\Gamma_{1,1}=\{\Gamma_{2,m}\}$, from Lemma \ref{jb} \ref{jb-2}, one obtains $\sum_{\tilde{t}\in\wz{\partial}(\Gamma)}p_{(1,2),(1,2)}^{\tilde{t}}p_{\tilde{t},(1,1)}^{(1,2)}=0$.
	By setting $i=j=h=(1,2)$ and $l=(1,1)$ in Lemma \ref{jb} \ref{jb-4}, we get   $p_{(1,2),(1,1)}^{(2,m)}p_{(1,2),(2,m)}^{(1,2)}=0$, and so $p_{(1,2),(2,m)}^{(1,2)}=0$. It follows from Lemma \ref{jb} \ref{jb-0} and \ref{jb-2} that $\Gamma_{1,2}\Gamma_{2,1}=\{\Gamma_{0,0}\}$, and so $k_{1,2}=1$.
	Since $T=\{2,3\}$ and $k_{1,1}=1$, $\Gamma$ has valency $2$, a contradiction. Thus,  $\Gamma_{2,2},\Gamma_{2,3}\notin \Gamma_{1,2}\Gamma_{1,1}$.

	In view of  \eqref{4ni-1}, \eqref{4ni-2}, Lemma \ref{4nilemma} \ref{4nilemma-2} and Lemma \ref{s42}, one has $\Gamma_{1,2}\Gamma_{1,1},\Gamma_{1,2}^2,\Gamma_{1,1}^2\subseteq \{\Gamma_{0,0},\Gamma_{1,1},\Gamma_{1,2},\Gamma_{2,1}\}$. Then $\Gamma_{2,2},\Gamma_{2,3}\notin \tilde{\partial}(\Gamma)$. By Lemma \ref{4nilemma4}, we obtain $\Gamma_{1,2}\Gamma_{2,1}\subseteq \{\Gamma_{0,0},\Gamma_{1,1},\Gamma_{1,2},\Gamma_{2,1}\}$. It follows that
	$\tilde{\partial}(\Gamma)\subseteq \{(0,0),(1,1),(1,2),(2,1)\}$. Then $\Gamma$ is isomorphic to the digraph in Theorem \ref{main thm1} \ref{main thm1-2}.
\end{proof}

\begin{proof}[Proof of Theorem \ref{main thm1}]
	Routinely, all digraphs in Theorem \ref{main thm1} are $4$-quasi-transitive. By \cite[Proposition 2.4, Theorem 1.1]{KSW03}  and \cite[Lemma 9.1]{lishuang2}, all digraphs in  Theorem \ref{main thm1} are commutative weakly distance-regular digraphs.
	
	Now we prove the necessity. If $\Gamma$ has cycle polar diameter $4$ or  diameter at least $4$, from Theorem \ref{main thm} and Corollary \ref{cor1}, then $\Gamma$ is isomorphic to one of the digraphs in Theorem \ref{main thm1} \ref{main thm1-3} and \ref{main thm1-6}.
	If $\Gamma$ has cycle polar diameter at most $3$ and diameter at most $3$, from Propositions \ref{s41} and \ref{s43}, then $\Gamma$ is isomorphic to one of the digraphs in Theorem \ref{main thm1} \ref{main thm1-2} and \ref{main thm1-7}.
\end{proof}

\section*{Acknowledgements}

Y.~Yang is supported by NSFC (12101575, 12671382) and the Fundamental Research Funds for the Central Universities (2652019319). 

\section*{Data Availability Statement}

No data was used for the research described in the article.

\section*{Conflict of Interest}

The authors declare that they have no conflict of interest.

\end{document}